\documentclass{scrartcl}

\usepackage{url} 

\usepackage{soul}

\usepackage{lmodern,lscape,textcomp,placeins}
\usepackage[T1]{fontenc}

\usepackage{authblk}

\usepackage{amsmath,amssymb,dsfont,units,booktabs}
\usepackage{tabularx}
\usepackage{appendix}
\usepackage{graphicx, color}
\usepackage{color}

\usepackage{natbib}
\usepackage[utf8]{inputenc} 

\usepackage{algorithm}
\usepackage{algpseudocode}


\newcommand{\vt}{\mathbf{v}}

\newcommand{\taut}{\pmb{\tau}}

\newcommand{\sigmat}{\pmb{\sigma}}

\newcommand{\epsilont}{\pmb{\epsilon}}

\newcommand{\te}{\pmb{e}}

\renewcommand{\div}{\operatorname{div}}
\newcommand{\tr}{\operatorname{tr}}

\begin{document}

\date{\small\textbf{Abstract}}
\title{Simulating linear kinematic features in viscous-plastic sea ice models on quadrilateral and triangular grids}

\author[1]{\small C. Mehlmann}
\author[2]{S. Danilov}
\author[2]{M. Losch}
\author[3]{J.F. Lemieux}
\author[2]{N. Hutter}
\author[4]{T. Richter}
\author[5]{P. Blain}
\author[6]{E.C. Hunke}
\author[1]{P. Korn}

\affil[1]{Max-Planck Institute for Meteorology, Hamburg, Germany } 
\affil[2]{Alfred-Wegener-Institut, Helmholtz Zentrum für Polar- und Meeresforschung, Bremerhaven, Germany}
\affil[3]{Recherche en Prévision Numérique Environnementale, Environnement et Changement Climatique Canada, Dorval, QC, Canada} 
\affil[4]{Otto-von-Guericke Universit\"at Magdeburg, Magdeburg, Germany} 
\affil[5]{Service  Météorologique du Canada, Environnement et Changement Climatique Canada, Dorval, QC, Canada}
\affil[6]{Los Alamos National Laboratory, Los Alamos, NM, USA}



\maketitle

\begin{abstract}%

Linear Kinematic Features (LKFs) are found everywhere in the Arctic sea ice cover. They are strongly localized deformations often associated with the formation of leads and pressure ridges.
Viscous-plastic sea ice models start to generate LKFs at high spatial grid resolution, typically with a grid spacing below 5\,km. Besides grid spacing, other aspects of a numerical implementation, such as discretization details, may affect the number and definition of simulated LKFs.
To explore these effects, the solutions of sea ice models with different grid spacings, mesh types, and numerical discretization techniques are compared in an idealized configuration, which could also serve as a benchmark problem in the future.
The A, B, and C-grid discretizations of sea ice dynamics on quadrilateral meshes leads to a similar number of LKFs as the  A-grid approximation on triangular meshes (with the same number of vertices).
The discretization on an Arakawa CD-grid on both structured quadrilateral and triangular meshes  resolves the same number of LKFs as conventional Arakawa A-grid, B-grid, and C-grid approaches, but on two times coarser meshes. This is due to the fact that the CD-grid approach has a higher number of degrees of freedom to discretize the velocity field.
 Due to its enhanced resolving properties, the CD-grid discretization is an attractive alternative to conventional discretizations.

\end{abstract}

\section{Introduction}
Sea ice is an {important} component of the climate system. As such, it is {crucial} for climate models to accurately represent interactions between the atmosphere, the sea ice cover and the ocean. 
Observations from Synthetic Aperture Radar (SAR) show that the sea ice cover is characterized by highly localized linear deformations (i.e.\ strain rates of the sea ice velocity field) \citep{Kwok2008}. These deformations are referred to as Linear Kinematic Features (LKFs). As LKFs give rise to leads and pressure ridges, they have a strong effect on sea ice production, salt rejection and the associated mixing, and on ocean-ice-atmosphere exchanges of momentum, heat and moisture.

For realistic simulations of LKFs, sea ice models require  an accurate representation of rheology. The rheology is the relation between applied stresses, material properties and the resulting deformations. 
In almost all climate models, sea ice is represented as a continuous viscous-plastic (VP) material \citep{Blockley2020} either in the classical VP framework \citep{Hibler1979} or following the  elastic-viscous-plastic (EVP) approximation \citep{Hunke1997}, although recently, alternative rheologies have been suggested \citep{Girard2011, Tsamados2013, Rampal2016}. 

All of these rheologies rely on the continuum assumption, 
which implies that statistical averages can be taken over a large number of floes \citep{Gray1994, Feltham2008}. Thus, the application of continuum rheological models at or below the scale of individual floes is only appropriate if the mode of failure of a single floe is the same as the mode of failure of an aggregate of floes \citep{Feltham2008}.
The validity of the continuum assumption and hence of the VP rheology at grid spacings of the size of ice floes is unclear  \citep{Coon2007, Feltham2008}, and this prompted the development of sea ice models based on discrete particles \citep[i.e., ice floes,][]{Wilchinsky2012, Herman2013}. In spite of the questionable validity of the VP model at grid spacings of the size of ice floes, 
the majority of practical applications still use the (E)VP model and will do so in the foreseeable future  \citep{Blockley2020}.

The widespread use of the (E)VP approach is also explained by the sparsity of available high frequency deformation observations, which makes it difficult to give preference to one rheology over others. Past research has therefore often focused on the development of efficient numerical methods \citep{Hunke1997,Hibler1997, Lemieux2009,Lemieux2010} for (E)VP models. These improvements, combined with the increased computer power, allow the simulation of well defined deformations and associated features such as leads. In the (E)VP model, LKFs are an emergent property at high spatial grid resolution of 5\,km and  finer \citep{Hutter2018}. The number of simulated LKFs also depends on the numerical convergence of the solving procedures \citep{Lemieux2012, Koldunov}.  To improve the numerical convergence, a number of methods have been proposed \citep{Lemieux2012, Losch2014, Kimmritz2015, MehlmannRichter2016newton}.

Recent analyses 
indicate that at high spatial resolution the VP rheology simulates many scalings and statistics of observed LKFs \citep{Hutter2020}. This result motivates our interest in sea ice models based on the VP rheology and leads to the following questions: 
(1) Which factors besides the grid size and numerical convergence influence the formation of LKFs?
 (2) How do current sea ice models differ in resolving these features?

To address these questions, we use a simple test case to compare several spatial discretizations on triangular and quadrilateral meshes with respect to their ability to resolve LKFs. The results are obtained with sea ice packages, such as the \emph{Los Alamos Sea Ice Model} \citep[CICE,][]{Hunke2015}, the sea ice module of the \emph{Massachusetts Institute of Technology
general circulation model} \citep[MITgcm,][]{LOSCH2010}, the \emph{Finite-Volume Sea Ice–Ocean model} \citep[FESOM,][]{Danilov2015} and sea ice module of the \emph{Icosahedral Nonhydrostatic Weather and Climate Model} \citep[ICON,][]{Mehlmann2020}, but also with the academic software library \emph{Gascoigne} \citep{Gascoigne}. We indeed show that the spatial discretization used for velocities plays an important role in resolving LKFs, in terms of both definition and number.\\

The paper is structured as follows. Section \ref{sec:models} presents the  viscous-plastic sea ice model. Section \ref{sec:testcase} describes the idealized test case. The Arakawa A, B, C, and CD-grid staggering and the corresponding discretization are introduced in Section \ref{sec:dis} and further described in the \ref{app:conf}. The numerical results are discussed in Section \ref{sec:num}. We conclude in Section \ref{sec:con}. 
\section{The viscous-plastic sea ice model}\label{sec:models}

For brevity, only a simplified model for three prognostic variables is presented here: the sea ice concentration $A$, the mean sea ice thickness $H$ and the sea ice velocity $\vt$. The sea ice dynamics are described by a system of coupled partial differential equations:

  \begin{align}
    \rho_\text{ice}H \partial_t \vt 
    + f_c \te_z\times \vt&=
    \div\,\sigmat + 
    A\taut(\vt)  -\rho_\text{ice}H g\nabla \tilde H_g,\label{eq:mom}\\
    \partial_t A + \div\,(\vt A) &= 0, \quad  A\le 1,\label{eq:A}  \\
    \partial_t H + \div\,(\vt H) & = 0, \label{eq:h}
  \end{align}
where $\rho_{ice}$ is the ice density, $f_c$ is the Coriolis
parameter, $g$ is the gravitational acceleration, $\tilde H_g$ is the sea surface height and $\te_z$ is the vertical ($z$-direction) unit vector.
The forcing term $\taut(\vt)$ is the sum of the ocean and atmospheric surface stresses, that is 

\begin{align}
 \taut(\vt) =
  C_{w}\rho_{w}|\vt_{w}-\vt|_2(\vt_{w}-\vt) +
  C_{a}\rho_{a}|\vt_{a}|_2\vt_{a},   
  \label{eq:forcing}
\end{align}
with the ocean velocity $\vt_{w}$ and the wind velocity
$\vt_{a}$. 
The internal stresses in the ice $\sigmat$ are modeled by the viscous-plastic (VP) sea ice rheology
\citep{Hibler1979}.

\begin{table}[t]
  \begin{center}
    \begin{tabular}{l|l|l}
      \toprule
      \text{Parameter} & \text{Definition} & \text{Value}\\
      \midrule
     $\rho_\text{ice}$ & sea ice density &$\unit[900]{kg/m^{3}}$\\
      $ \rho_\text{a}$& air density & $\unit[1.3]{kg/m^{3}}$\\
      $\rho_\text{w}$& water density & $\unit[1026]{kg/m^{3}}$\\
      $C_\text{a}$&air drag coefficient &$\unit[1.2] \cdot {10^{-3}}$\\
      $C_\text{w}$&water drag coefficient& $\unit[5.5] \cdot {10^{-3}}$\\
     $f_c$ &Coriolis parameter&$\unit[1.46] \cdot \unit[10^{-4}]{s^{-1}}$ \\
     $P^{\star}$&ice strength parameter&$\unit[27.5]\cdot \unit[ 10^{3}]{N/m^2}$\\
     $C$&ice concentration parameter&$20$\\
      $e$& ellipse aspect ratio&$2$\\
     \bottomrule
    \end{tabular}

 \caption{Physical parameters of the momentum equation\label{Con}.}
 \end{center}
\end{table}
The nonlinear viscous-plastic rheology
relates  the strain rate tensor

\begin{align}\label{strain}
 \dot\epsilont=\frac{1}{2}\Big(\nabla\vt+\nabla\vt^T\Big),\quad
\dot\epsilont':=\dot\epsilont-\frac{1}{2}\operatorname{tr}(\dot\epsilont)I ,   
\end{align}
where $\operatorname{tr}(\cdot)$ is the trace, to the stress tensor $\sigmat$. The relationship is given by

\begin{equation}\label{model:stress}
  \begin{aligned}
    \sigmat &= 2\eta \dot\epsilont' + \zeta \tr(\dot\epsilont)I -
    \frac{P}{2} I,%
    \end{aligned} 
\end{equation}
with the viscosities $\eta$ and $\zeta$, given by $\eta=e^{-2}\zeta$ and

\begin{equation}\label{deltalimit}
  \zeta=\frac{P_0}{2\Delta(\dot\epsilont)},\quad
  \Delta(\dot\epsilont):= 
  \sqrt{\frac{2}{e^2} \dot\epsilont':\dot\epsilont'+\tr(\dot\epsilont)^2 +
    \Delta_\text{min}^2}.
\end{equation}
$\Delta_\text{min}= 2\cdot 10^{-9} \, \unit{s^{-1}}$  is the threshold that describes
the transition between the viscous and the plastic regimes. 
The replacement pressure $P$ and the ice strength $P_0$ in (\ref{model:stress}) are respectively expressed as

\begin{equation}\label{icestrength}
  P=P_0\frac{\Delta}{(\Delta+\Delta_\text{min})},\quad 
  P_0(H,A)= P^\star H \exp\big(-C(1-A)\big).
\end{equation}
All parameters are collected in
Table~\ref{Con}. 

\section{Test case}\label{sec:testcase}
The test case corresponds to the initial phase of the deformation of sea ice caused by a moving cyclone.
We consider the quadratic domain
$\Omega=(0,\unit[512]{km})^2$  and measure the time $t $ in $\unit{days}$. The simulation is run for $T=[0,2]$ days.
We prescribe a
circular steady ocean current

        \begin{align}
  \vt_{ocean}&=\bar v_{ocean}^{max}\begin{pmatrix}
 \phantom{-}(2y-L)/L\\
 -(2x-L)/L
  \end{pmatrix},
\end{align}
with $L=512000$ m and

\[
\bar v_{ocean}^{max} = \unit[0.01]{m\, s ^{-1}}. 
\]
The wind field is described by a cyclone which moves from the center of the domain to the upper right edge. The center of the cyclone moves in time as

\begin{align}
m_x(t)=m_y(t) =   \unit[51200]{m}+\unit[51200]{m\,day^{-1}} \cdot t.
\end{align}
The maximal wind speed is set to $v_\text{max}=\unit[15]{m\,s^{-1}}$. To reduce the wind speed away from the center, it is multiplied by the factor

\begin{align}
  s =& \frac{1}{50}\exp(-0.01r), \quad  r=\sqrt{(m_x-x)^2+(m_y-y)^2}.
\end{align}
Given the 
  convergence angle $\alpha=72^\circ$ for the cyclone, 
the wind vector is finally expressed as

\begin{align}
  \vt_{a} &= 
  - s \cdot {v_{\max}}
  \begin{pmatrix} 
  \phantom{-}\cos(\alpha)(x-m_x) + \sin(\alpha)(y-m_y) \\
  -\sin(\alpha)(x-m_x) + \cos(\alpha)(y-m_y)
  \end{pmatrix}.
\end{align}

We initialize the simulation with sea ice at rest and assume a constant ice concentration of 1.0 and a small perturbation of the ice thickness around a mean of $H^0(x,y)=\unit[0.3]{m}$. These initial conditions are

\begin{align}\label{initial}
\vt(0,x,y)&=\vt^0(x,y):=\unit[0]{m\,s^{-1}},\\\label{initial:A} A(0,x,y)&=A^0:=1,\\\label{initial:H}
H(0,x,y)& = H^0(x,y):=\unit[0.3]{m} + \unit[0.005]{m}\left(
  \sin\left(6\cdot 10^{-5} x\right)
  +\sin\left(3\cdot 10^{-5}y\right)\right),
\end{align}
where $x,y$ are given in meters.
At the boundary of the domain we apply a no-slip condition for the velocity, i.e.

\begin{align}
 \vt=0\text{ on }\partial\Omega.    
\end{align}

\section{Discretization}\label{sec:dis}
\begin{figure}
\begin{center}
 \begin{tabular}{c c c c c}
 & A-grid& B-grid& CD-grid& C-grid\\
  \hline
FE:&Q1-Q1 &Q1-Q0& CR-Q0& - \\
&\includegraphics[scale=0.4]{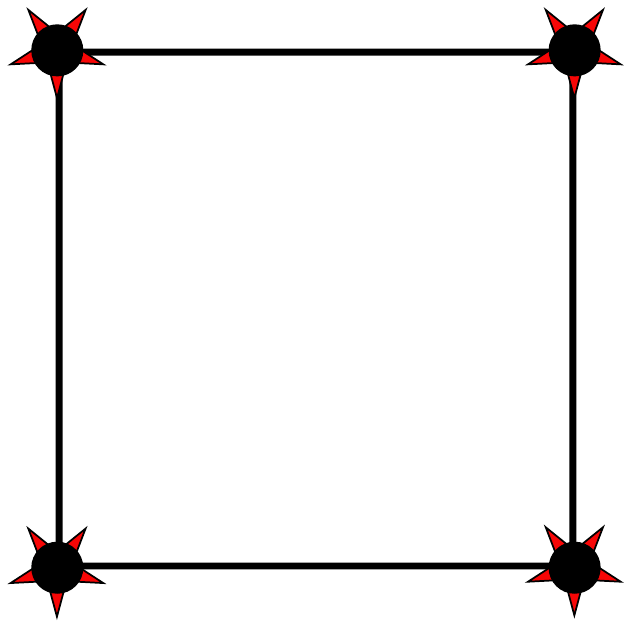}
&\includegraphics[scale=0.4]{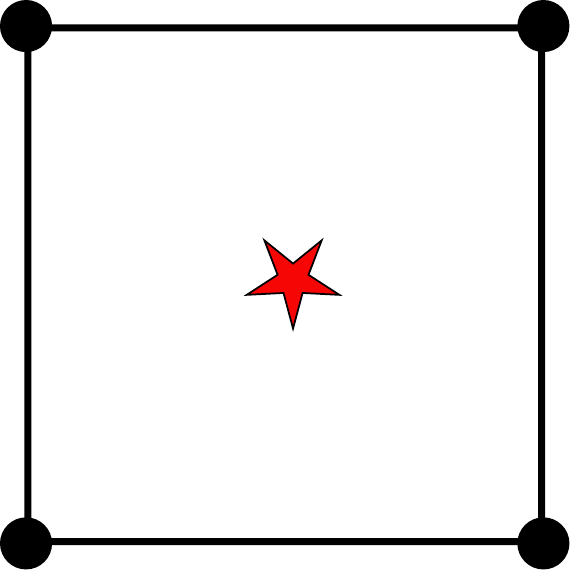}& \includegraphics[scale=0.4]{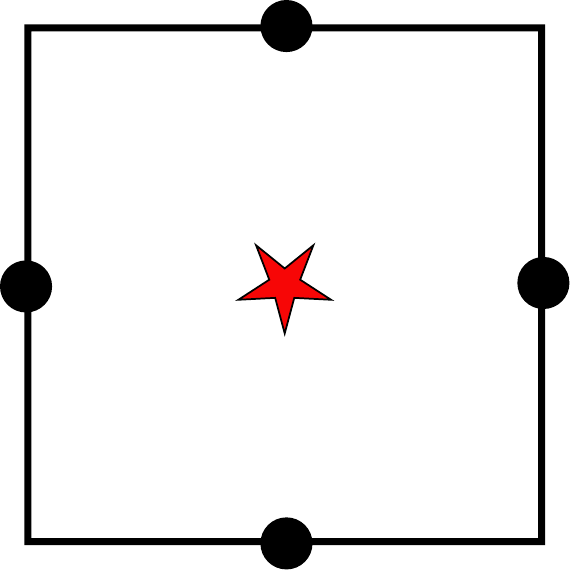}&\includegraphics[scale=0.4]{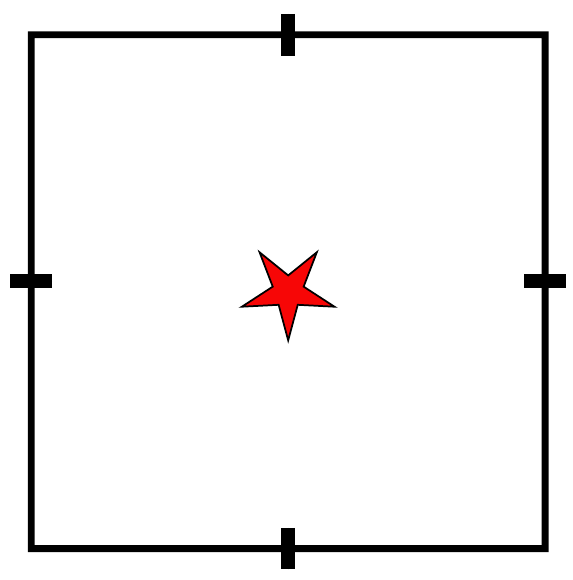}\\
FE:&P1-P1 &P0-P1& CR-P0& -\\
&\includegraphics[scale=0.4]{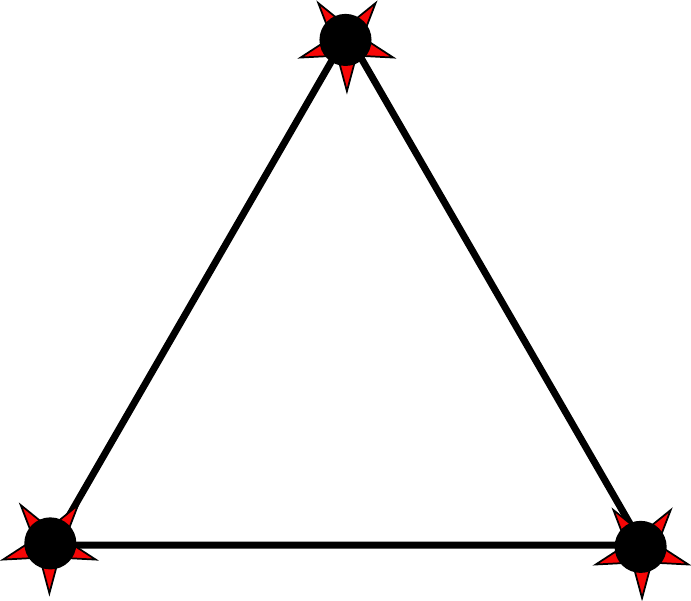}&\includegraphics[scale=0.4]{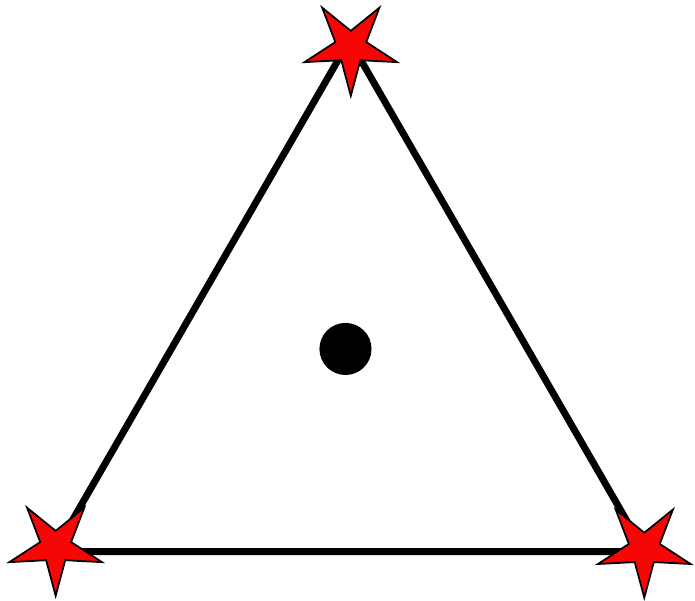}& \includegraphics[scale=0.4]{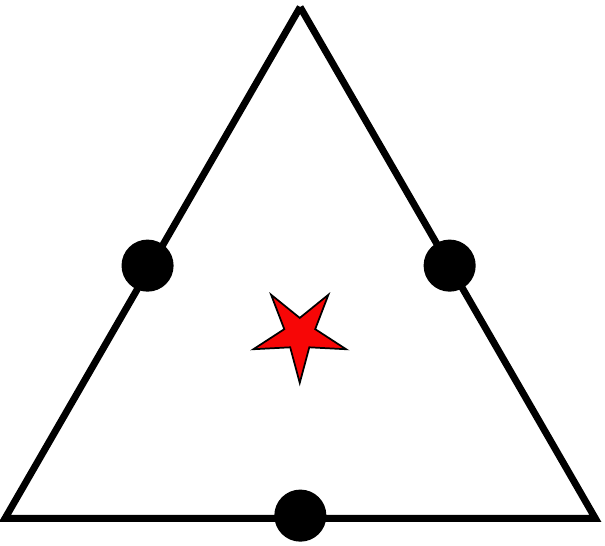}\\
\end{tabular}
\caption{Different  staggering explored in this manuscript and their finite element (FE) equivalents. 
We indicate the placement of the horizontal velocity $\vt=(v^1,v^2)$, the normal component of the velocity and the staggering of the tracers by $\bullet$, \textbf{-}, and {\color{red}$\star$}, respectively.\label{fig:FE:quad_tri}}   
\end{center}
\end{figure}
Most sea ice models use structured (quadrilateral) meshes and either an Arakawa B-grid or Arakawa C-grid discretization. We are not aware of the use of Arakawa A-grids or CD-grids on quadrilateral meshes in sea ice models. In this paper, we compare Arakawa B-grids and C-grids to Arakawa A-grids and CD-grids taking into account structured quadrilateral grids (upper row in
Figure~\ref{fig:FE:quad_tri}) and triangular meshes (lower row in Figure~\ref{fig:FE:quad_tri}). 


We use the sea ice model CICE \citep{Hunke2015} for Arakawa B-grid simulations and the sea ice module of the MITgcm \citep{LOSCH2010} for Arakawa C-grid simulations on quadrilateral grids.
The discretizations based on the structured A-grid and CD-grid type staggering are carried out in the finite element software library Gascoigne \citep{Gascoigne}. The Arakawa A and CD-grids  are equivalent to the Q1-Q1 and CR-Q0 finite element pairs, respectively. The first component of the pair refers to the discretization of the sea ice velocity, and the second component addresses the discretization of the sea ice thickness and concentration. Q0 denotes the piecewise constant element and  Q1 refers to the bilinear quadrilateral element. On quadrilaterals 
CR is the noncoforming rotated bilinear element \citep{RannacherTurek1992}, which is a variant of the triangular nonconforming piecewise linear Crouzeix-Raviart element \citep{CrouzeixRaviart1973}.  In fact, the B-grid discretization of CICE uses some aspects  of a finite-element approach through a variational formulation \citep{Hunke2002}. 
The MITgcm implements a finite volume discretization on a C-grid. Only strain-rates are approximated by finite differences \citep{LOSCH2010}.

  Recent developments on unstructured meshes was accompanied by the use of triangular grids. This includes the A-grid like P1-P1 finite element which is used in the sea ice component of FESOM \citep{Danilov2015}, the CD-grid like CR-P0 finite element pair  applied in the sea ice module of ICON \citep{Mehlmann2020} and a quasi-B grid discretization used in sea ice module of the \emph{Finite‐Volume Community Ocean Model}
\citep[FVCOM,][]{Gao2011}. The latter is similar to the quasi-B grid staggering realized on hexagonal meshes used in  the \emph{Model for Prediction Across Scales} \citep[MPAS,][]{Petersen2019}. This staggering corresponds to a P0-P1 finite element pair. In analogy to the quadrilateral case, P0, P1, and CR refer to the piecewise constant element, the piecewise linear element, and the nonconforming linear element on triangles. On triangular meshes we compare the A-grid, B-grid, and CD-grid type approximations implemented in the sea ice module of FESOM to the CD-grid type discretization in ICON. The Arakawa C-grid is not considered on triangular grids as the corresponding finite element space is not large enough to approximate the full strain rate tensor \citep{Acosta2011} of the VP sea ice rheology. Instead, this issue motivates the use of the Arakawa CD-grid on triangular meshes.
Discretizing the sea ice momentum equation on quadrilateral and triangular meshes with the CD-grid like CR-element causes oscillations in the velocity field. To reduce these oscillations, a stabilization approach is used. A short description of the stabilization can be found in \ref{app:conf}. 

The standard approach for solving the coupled sea ice system (\ref{eq:mom})-(\ref{eq:h}) is a time splitting method \citep{Lemieux2014}. 
We first  compute the solution of the sea ice  momentum equation (\ref{eq:mom}), followed by the solution of the  transport equations (\ref{eq:A}) and (\ref{eq:h}). For stability reasons, a fully explicit time
stepping scheme for the momentum equation with a VP rheology would require a time step smaller than a second - even on a grid resolution as coarse as  $\unit[100]{km}$ \citep{Hibler1991}.
Therefore the authors recommended an implicit treatment in time \citep{Hibler1991}. An implicit time discretization requires iterative methods such as a Picard solver \citep{Hibler1997,Lemieux2009} or Newton like methods \citep{Lemieux2010,Losch2014, MehlmannRichter2016newton}. To avoid an implicit discretization, the EVP model was introduced \citep{Hunke1997}\, where an artificial elastic term added to the VP rheology allows an explicit discretization of the momentum equation with relatively large time steps. However, {the EVP model does not lead to the same deformations as simulated with a VP model.} Thus, a 
modified version of the EVP method was developed that ensures convergence to the solution of the VP model \citep{Lemieux2012, Boullion2013, Kimmritz2015}.  We refer to this pseudo-time stepping method as the mEVP solver. In this manuscript, we consider an implicit discretization of the VP model and explicit discretization with the mEVP solver. 

In the VP case, the sea ice momentum equation is approximated with an implicit Euler scheme in time. We solve the resulting nonlinear system with a Newton scheme, which is converged to a certain tolerance, given in \ref{app:conf}.
As there is no Newton-type method available in CICE, the solution is computed with a Picard solver. The nonlinear solution is approximated by performing 100 Picard iterations. A higher number of Picard iterations does not affect the number of resolved LKFs. 
Similarly, a higher solver tolerance for the Newton methods does not influence the number of resolved LKFs.


In the triangular configuration the discretization of the sea ice momentum equation is based on the explicit mEVP solver. For a large number of sub-cycles per time step the method converges to 
the VP solution \citep{Kimmritz2015}. 
However, in practice only a limited number of iterations of the  mEVP are applied to reduce the numerical cost \citep{Koldunov}. Thus we focus on an efficient approximation and use  $100$ sub-cycling steps per time step. 
In the test case considered, the number of LKFs does not significantly grow if the number of sub-iterations is increased. 


The transport scheme in the tracer advection differs between the A-grid, B-grid, C-grid and CD-grid like approximations. An overview is given in  \ref{app:conf}.  In the test case, which is run for a short period of 2~days, the influence of the advection scheme on the number of resolved LKFs is small. This allows us to focus on the role of the velocity staggering.

The test case is analyzed on quadrilateral and triangular meshes with grid spacings of 8\,km, 4\,km, and 2\,km.
The quadrilateral grid has respectively 4096, 16384, and 65536 cells on the 8\,km, 4\,km, and 2\,km meshes. The triangular mesh in FESOM contains respectively 9490 cells (8\,km), 37926 cells (4\,km), and 151630 cells (2\,km).
Due to the use of a uniform mesh in this test case, the triangular 8\,km grid in ICON has 9070 cells, the 4\,km grid has 37082 cells, and the 2\,km mesh contains 149938 cells.
To compare the approximation on quadrilateral grids to the discretization on triangular meshes we chose grids with the same number of vertices. Therefore the edge length in the triangular case is slightly larger, 8.6 km, 4.3 km and 2.15 km. Keeping the same number of vertices in the triangular grid results in twice as many cells and 1.5 times more edges. 

We solve the system  on all mesh levels with a time step of $\Delta t=2$ min. The choice of time step is determined by the use of the explicit mEVP solver on the 2\,km grid. The CR-P0 discretization requires slightly smaller time steps than the P0-P1 or the P1-P1 discretization. 
Using the implicit solvers, we observed that a larger time step ($\Delta t=30$\,min) slightly affects the position of the resolved LKFs, but did not significantly influence the number of the resolved LKFs (not shown).

\section{Numerical evaluation}\label{sec:num}
In this section, we examine the effect
of the grid staggering on the formation of LKFs.
The number of LKFs is affected by different factors such as the solver convergence \citep{Lemieux2009,Koldunov}, the mesh resolution \citep{Hutter2018} or the ice strength parameterization \citep{Hutter2020}. 
To quantify the effect of velocity staggering on the formation of LKFs, we analyze the sea ice concentration $A$ and shear deformation
\begin{align*}
     \epsilon_{II}=\sqrt{(\dot\epsilont_{11}-\dot\epsilont_{22})^2 -4 \dot\epsilont_{12}^2},
\end{align*}
where $\dot \epsilont_{11}, \dot \epsilont_{22}$ and $\dot \epsilont_{12}$  are the entries of the strain rate tensor given in (\ref{strain}). The number of LKFs is determined with a recently developed LKF detection algorithm \citep{Hutter2019}. 

As the detection algorithm requires data on a regular grid, all model outputs are interpolated on a 2 km cartesian mesh.
We adapt the original version of the algorithm \citep{Hutter2019} for our idealized experiments by the following minor changes: (1) We do not use the histogram equalization in the LKF pixel filtering process, as all simulated shear fields have the same range of magnitudes. (2) The maximum and minimum kernel sizes of the Difference of Gaussian (DoG) filter is set to $6\cdot\frac{\Delta x}{2\,km}$ and $1.2\cdot\frac{\Delta x}{2\,km}$, to ensure that the detected LKFs are wider than one pixel and no grid-scale noise is detected. Here, $\Delta x$ is the length of the grid edge. (3) We use a filter threshold for the DoG of  $\unit[0.1]{\log_{10}\Big(\frac{1}{s}\Big)}$. (4) The smallest length of detected LKFs is set to $4.8\cdot\Delta x$.

\subsection{Quadrilaterals}\label{sec:quad}
 \begin{figure}[t]
  \begin{center}
    \begin{tabular}{c c c c | c   }
     & A-grid& B-grid & CD-grid & C-grid \\
     &(Gascoigne)& (Gascoigne) &(Gascoigne)& (MITgcm) \\
     \rotatebox{90}{\phantom{abcde}4 km}&  \includegraphics[scale=0.04]{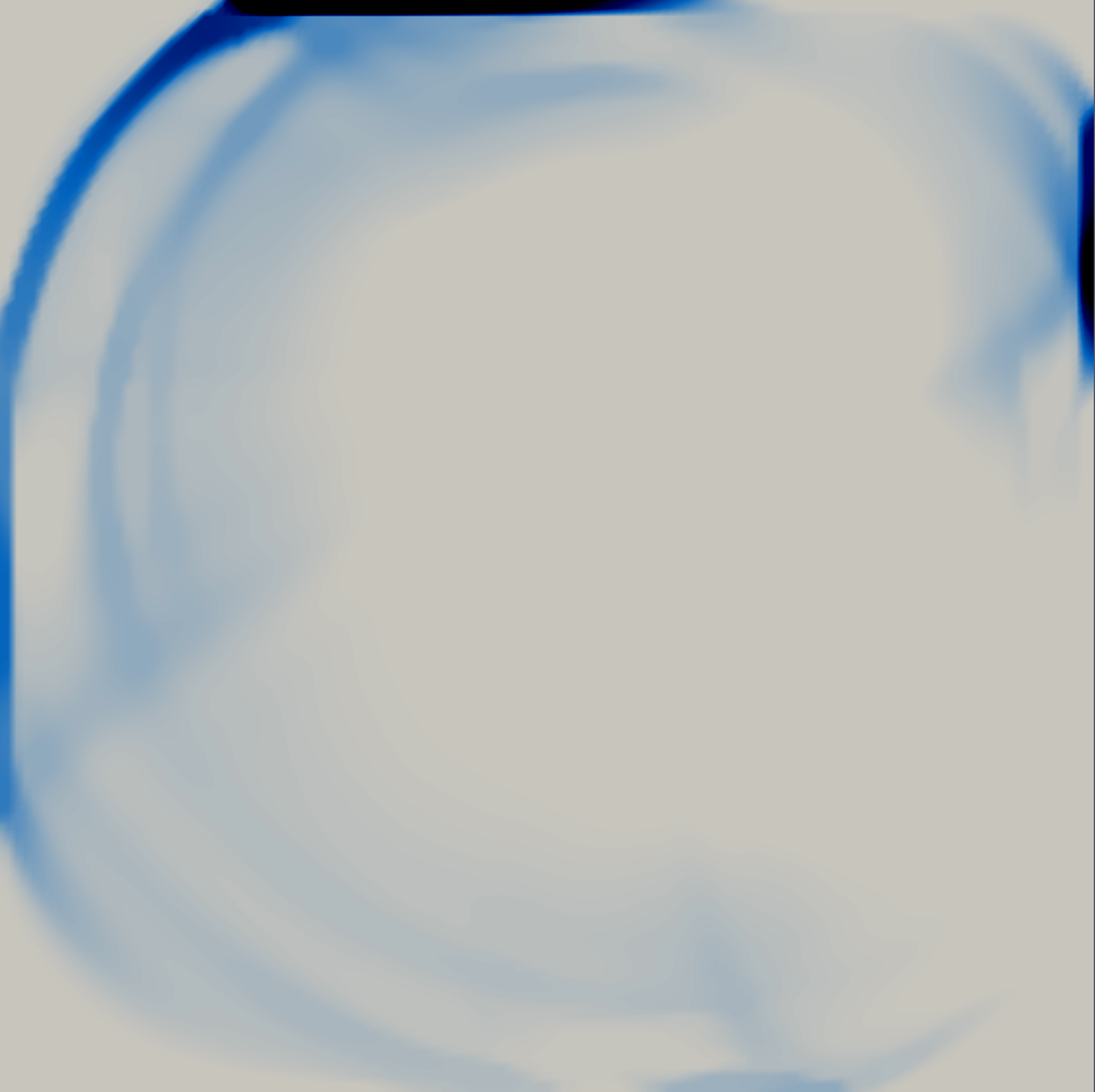}& 
       \includegraphics[scale=0.04]{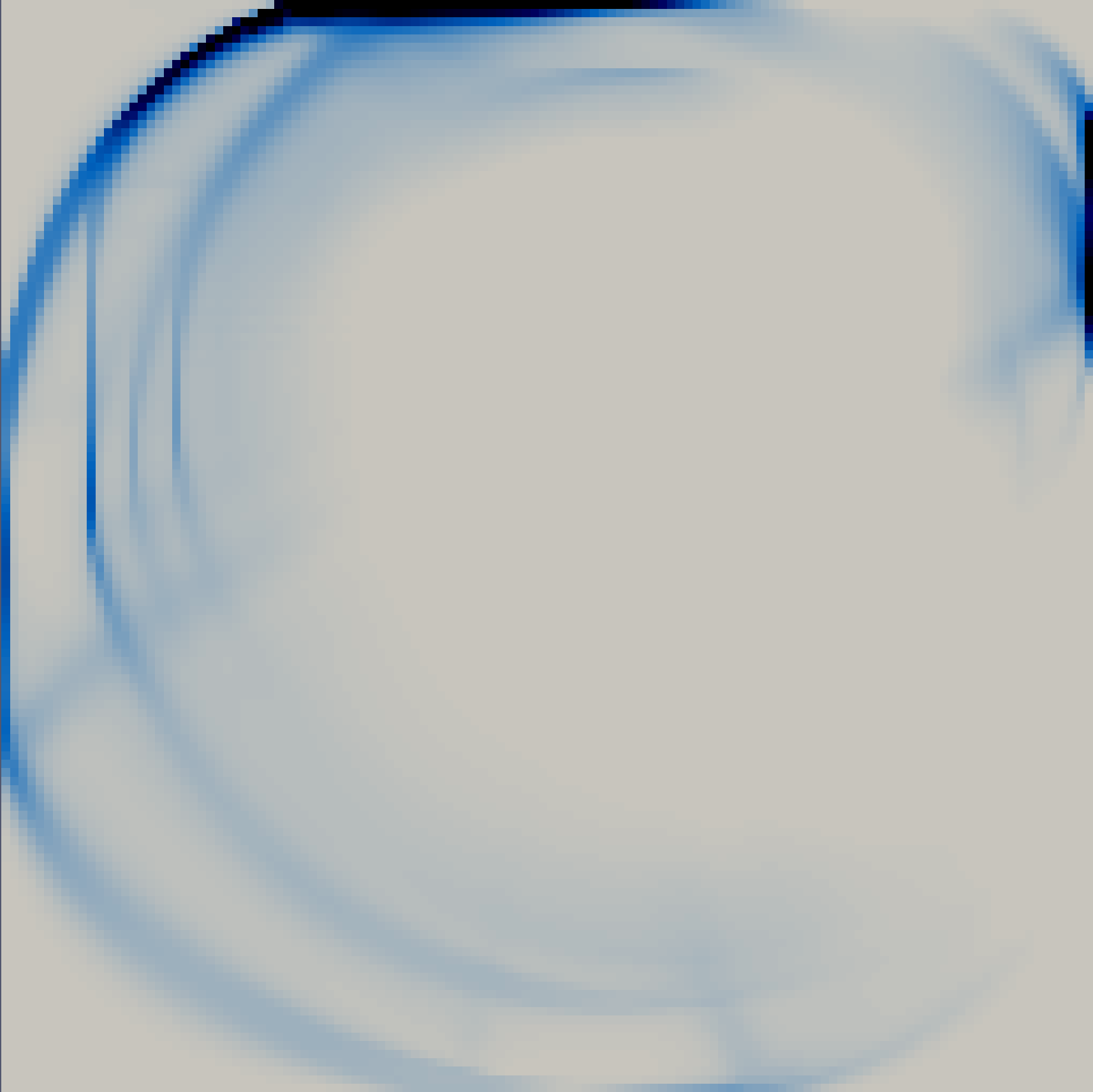}& 
      \includegraphics[scale=0.04]{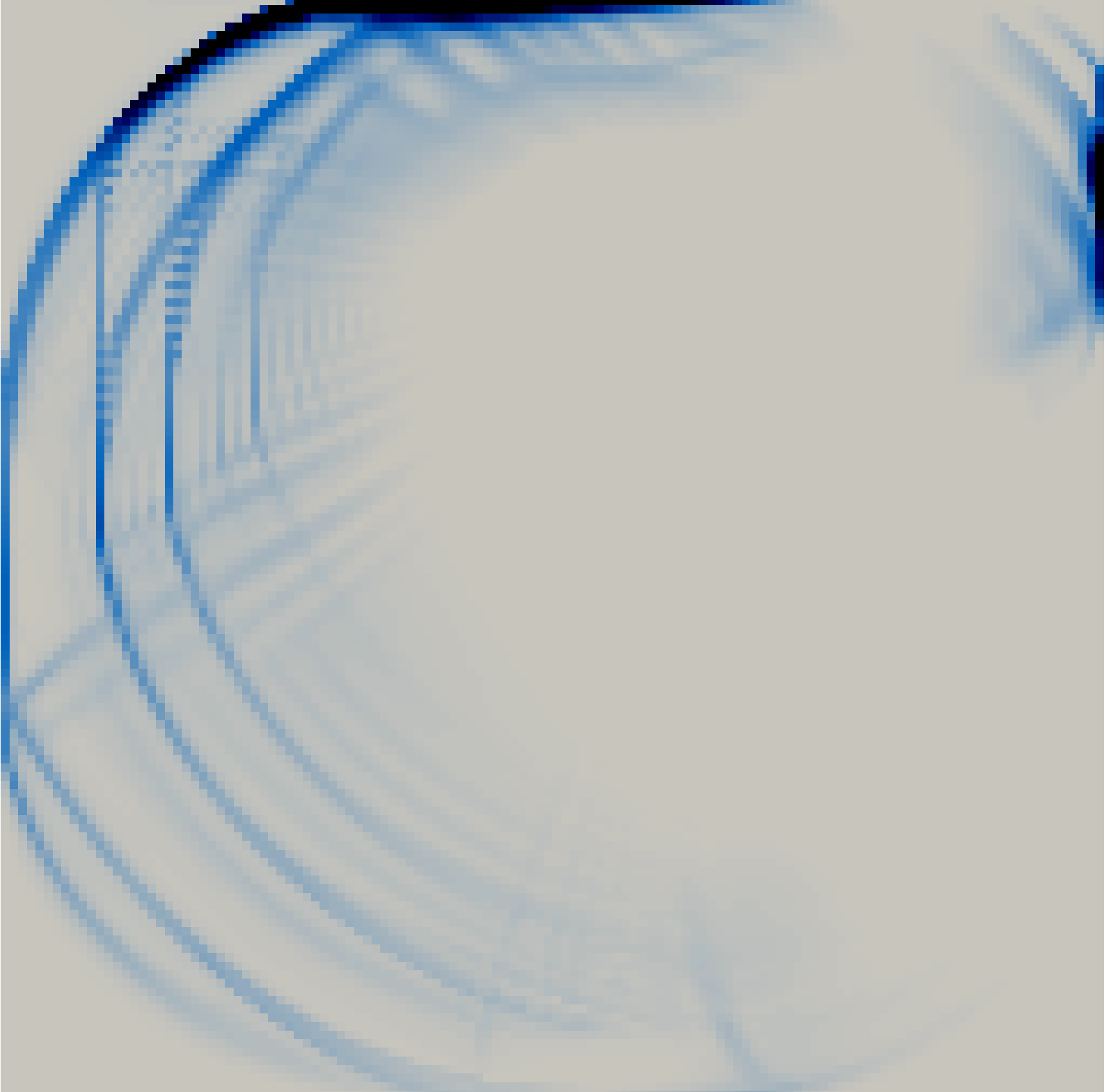}&  
\includegraphics[scale=0.04]{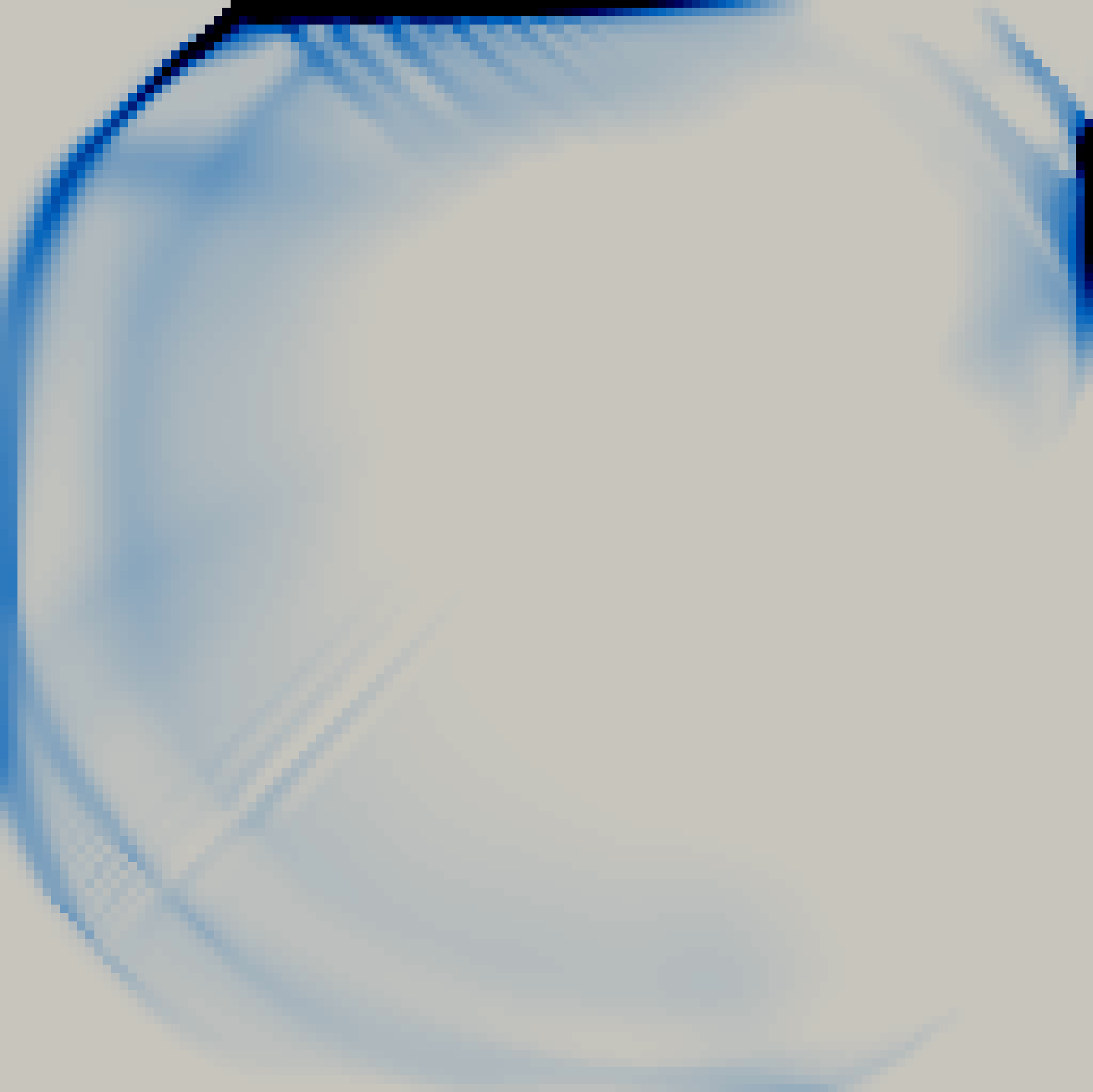}\\
     \rotatebox{90}{\phantom{abcde}2 km}&        \includegraphics[scale=0.054]{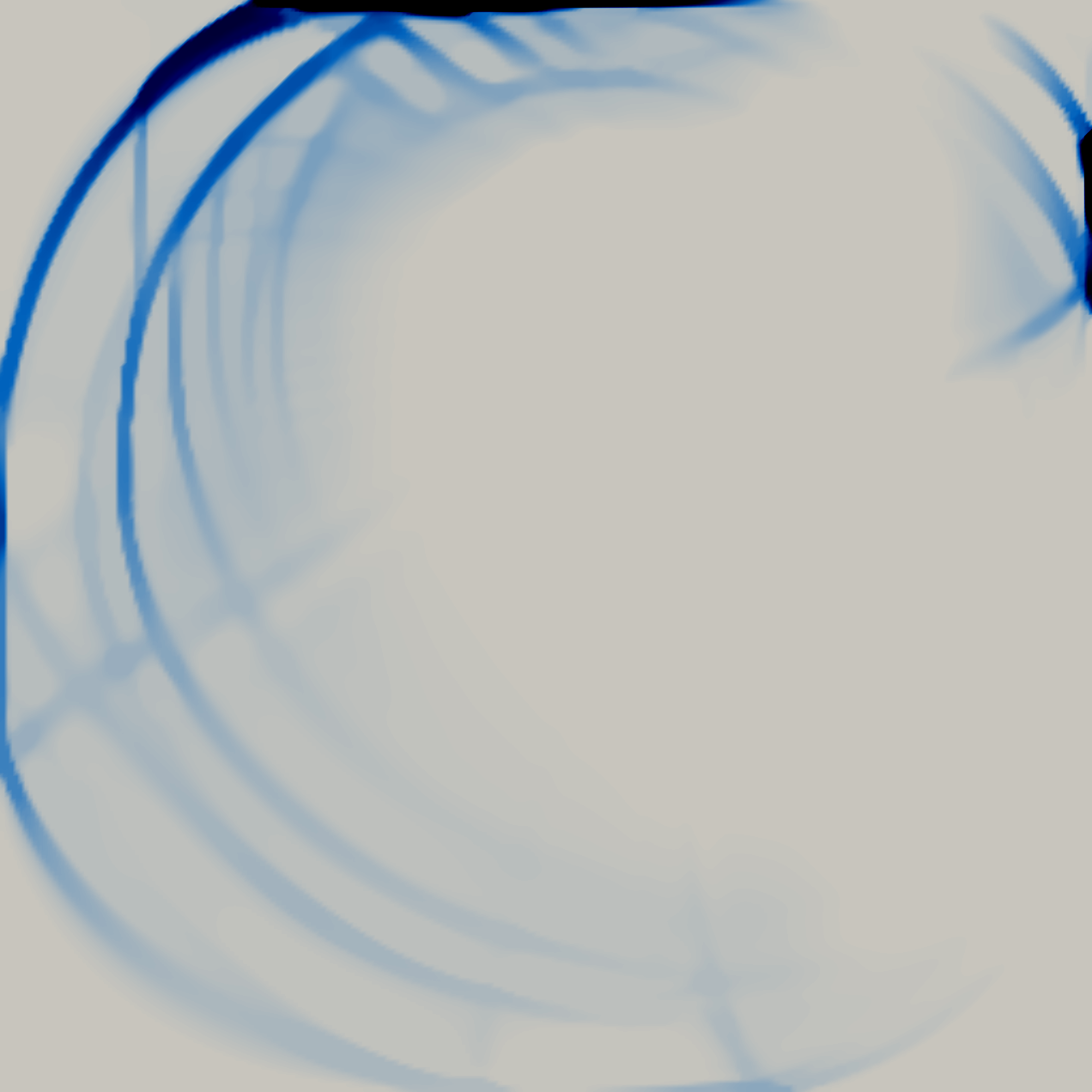}&
      \includegraphics[scale=0.04]{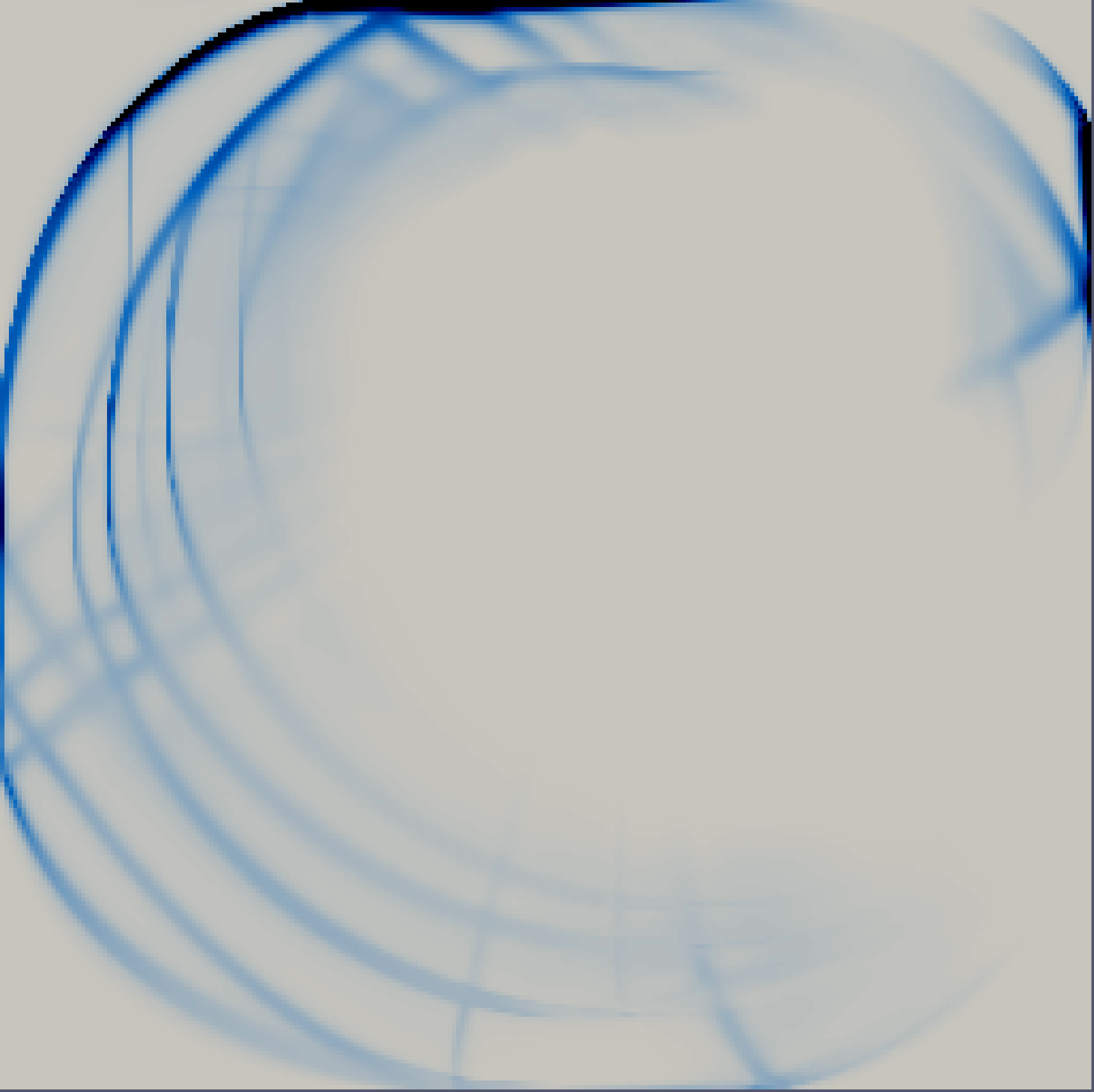}& 
         \includegraphics[scale=0.04]{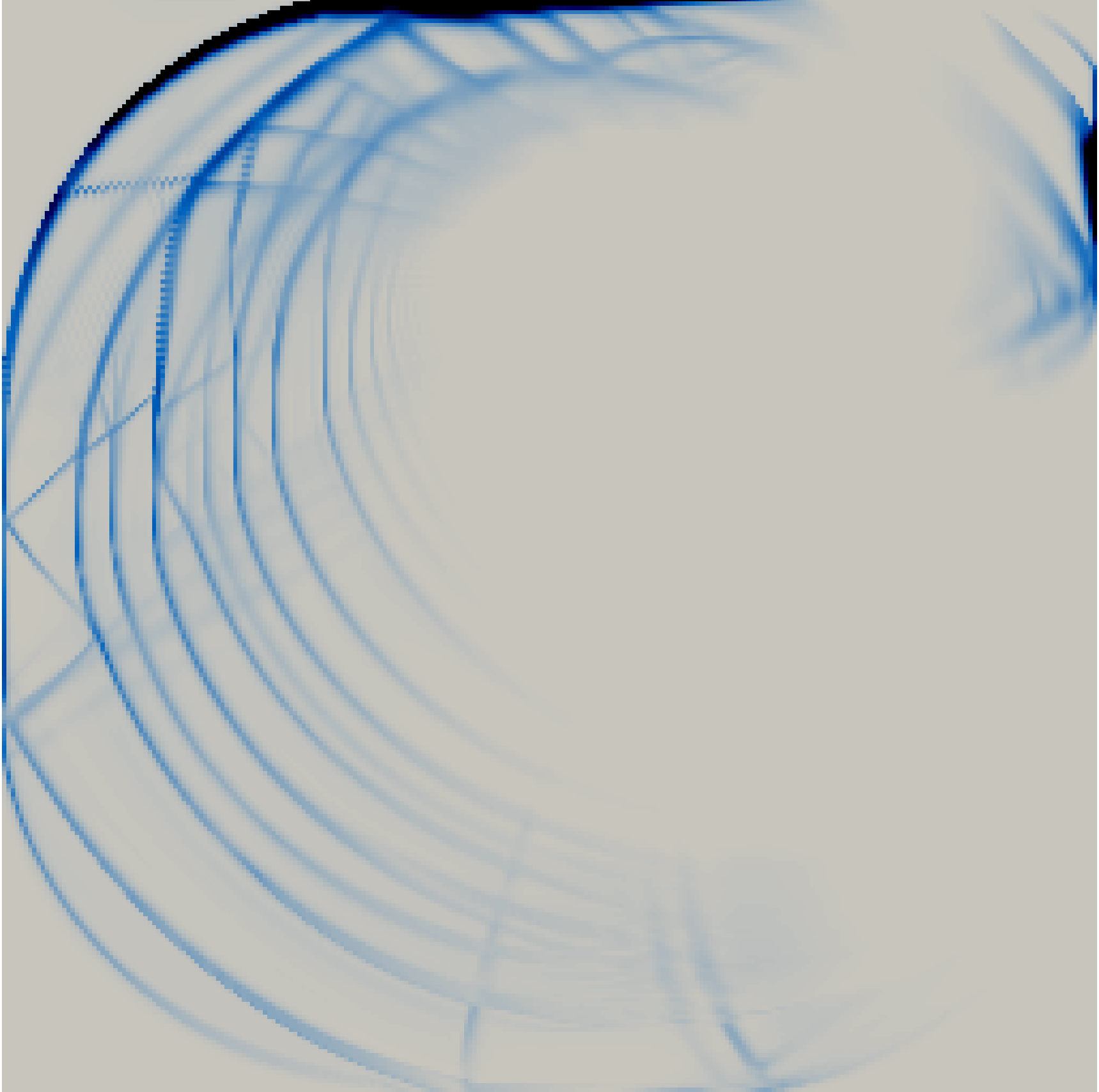}& 
        \includegraphics[scale=0.04]{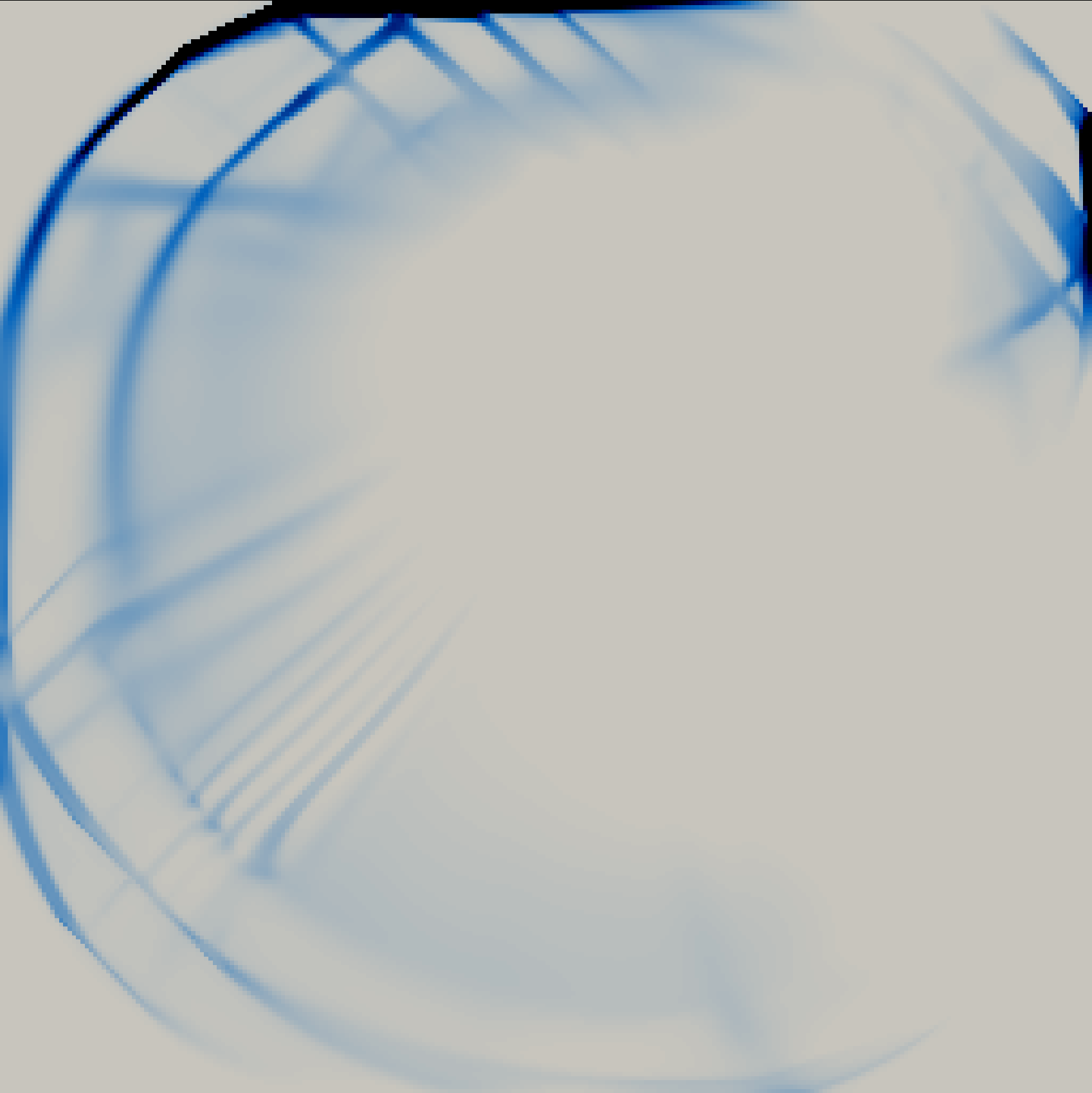}
       \\
&         \includegraphics[scale=0.16]{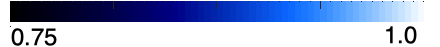}&  \includegraphics[scale=0.16]{pics_CR/ICON_A_sc}&
         \includegraphics[scale=0.16]{pics_CR/ICON_A_sc}&
         \includegraphics[scale=0.16]{pics_CR/ICON_A_sc}\\
     \end{tabular}
  \end{center}
  \begin{center}
   \caption{Sea ice concentration computed on a quadrilateral mesh. 
   \label{fig:Cgrid}}
  \end{center}
\end{figure}

We start by considering the discrete solutions obtained with the software library Gascoigne \citep{Gascoigne} using different finite element discretizations. The B-grid  and the CD-grid like finite element discretization differ only by the staggering of the velocity components while the A-grid type and B-grid like finite element discretizations are identical except for the placement of the sea ice concentration and thickness.\\

As both B-grid and CD-grid like approximations are based on an upwind scheme, we can attribute the relatively large differences in the sea ice concentration in Figure~\ref{fig:Cgrid} to the different velocity placement.
In the case of the A-grid discretization, however, the sea ice concentration and sea ice thickness are advected with a second order Taylor-Galerkin flux-correction scheme \citep{Mehlmann2019}, so that the much smaller differences between A-grid and B-grid  (Figure~\ref{fig:Cgrid}) is a consequence of a combination of different tracer point placement and the related difference in advection schemes. Therefore, neither the advection scheme nor the placement of the tracers are important for the evolution of LKFs in this test case.

\begin{figure}[tp]
  \begin{center}
    \begin{tabular}{c c c  |c  c   }
     & B-grid& CD-grid & B-grid & C-grid \\
     & (Gascoigne)& (Gascoigne) &(CICE) & (MITgcm) \\
     \rotatebox{90}{\phantom{abcde}4 km}&  
       \includegraphics[scale=0.04]{pics_CR/Q1_Q0_A_4km_2min.png}& 
      \includegraphics[scale=0.04]{pics_CR/CR_Q0_A_4km_02.png}&
      \includegraphics[scale=0.04]{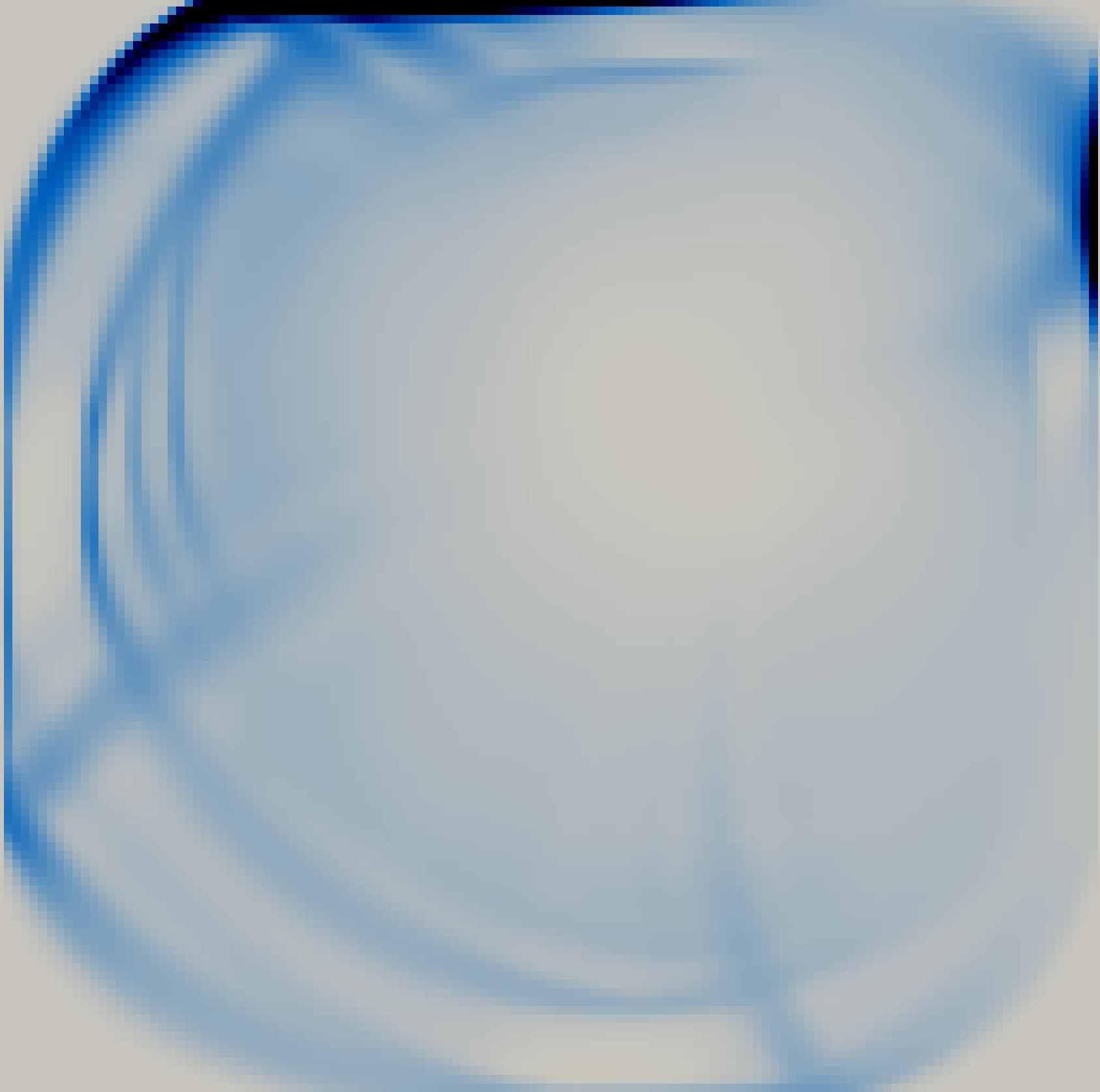}&  
\includegraphics[scale=0.04]{pics_CR/MITgcm_A_4km_2min.png}\\
     \rotatebox{90}{\phantom{abcde}2 km}&        \includegraphics[scale=0.04]{pics_CR/Q1_Q0_A_2km_2min.png}& 
         \includegraphics[scale=0.04]{pics_CR/CR_Q0_A_2km_01.png}& 
         \includegraphics[scale=0.04]{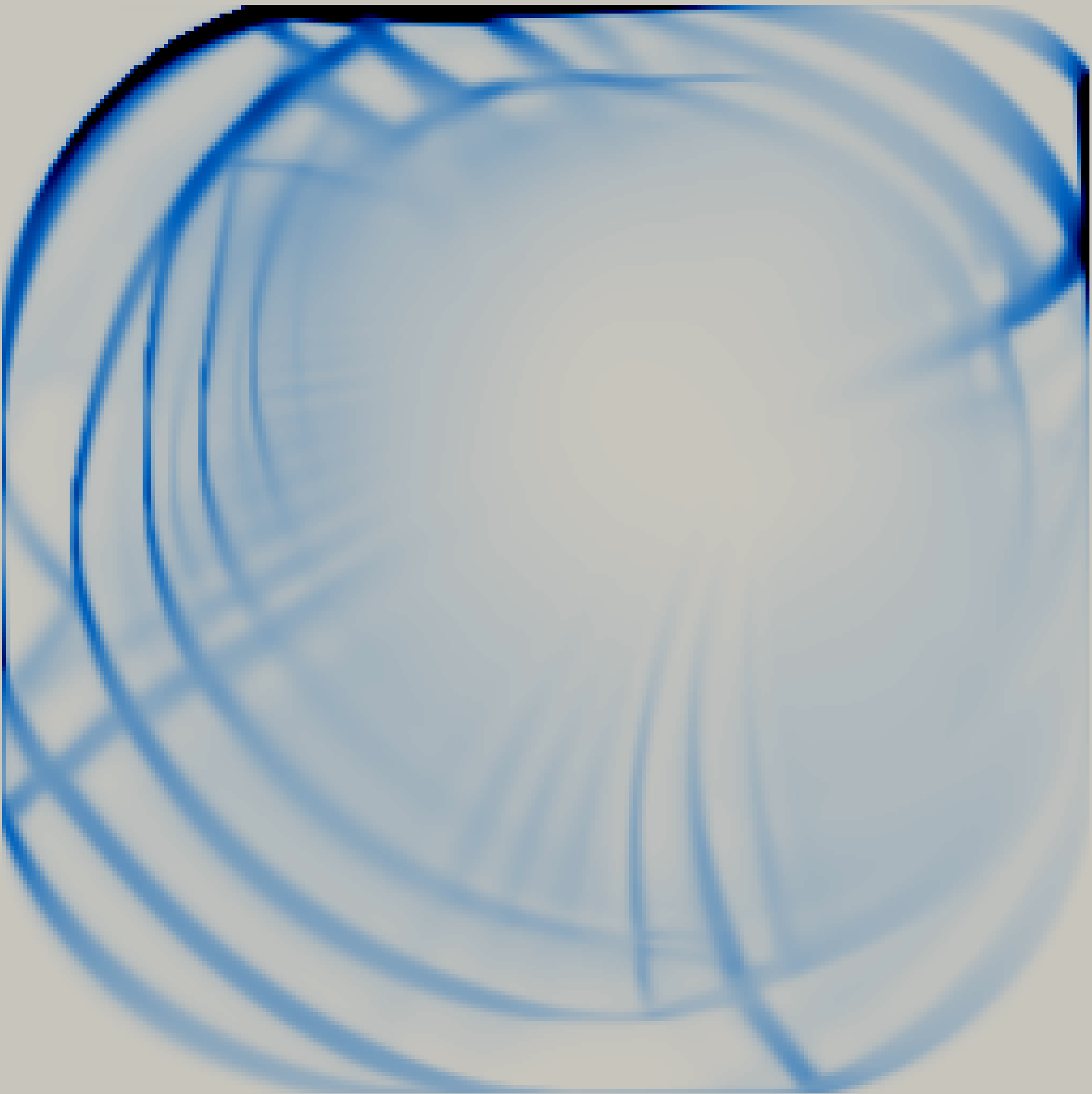}&
        \includegraphics[scale=0.04]{pics_CR/MITgcm_A_2km_2min.png}
       \\
&         \includegraphics[scale=0.16]{pics_CR/ICON_A_sc}& 
          \includegraphics[scale=0.16]{pics_CR/ICON_A_sc}&    
          \includegraphics[scale=0.16]{pics_CR/ICON_A_sc}&    
         \includegraphics[scale=0.16]{pics_CR/ICON_A_sc}\\
         \hline
              \rotatebox{90}{\phantom{abcde}4 km}&  
 \includegraphics[scale=0.04]{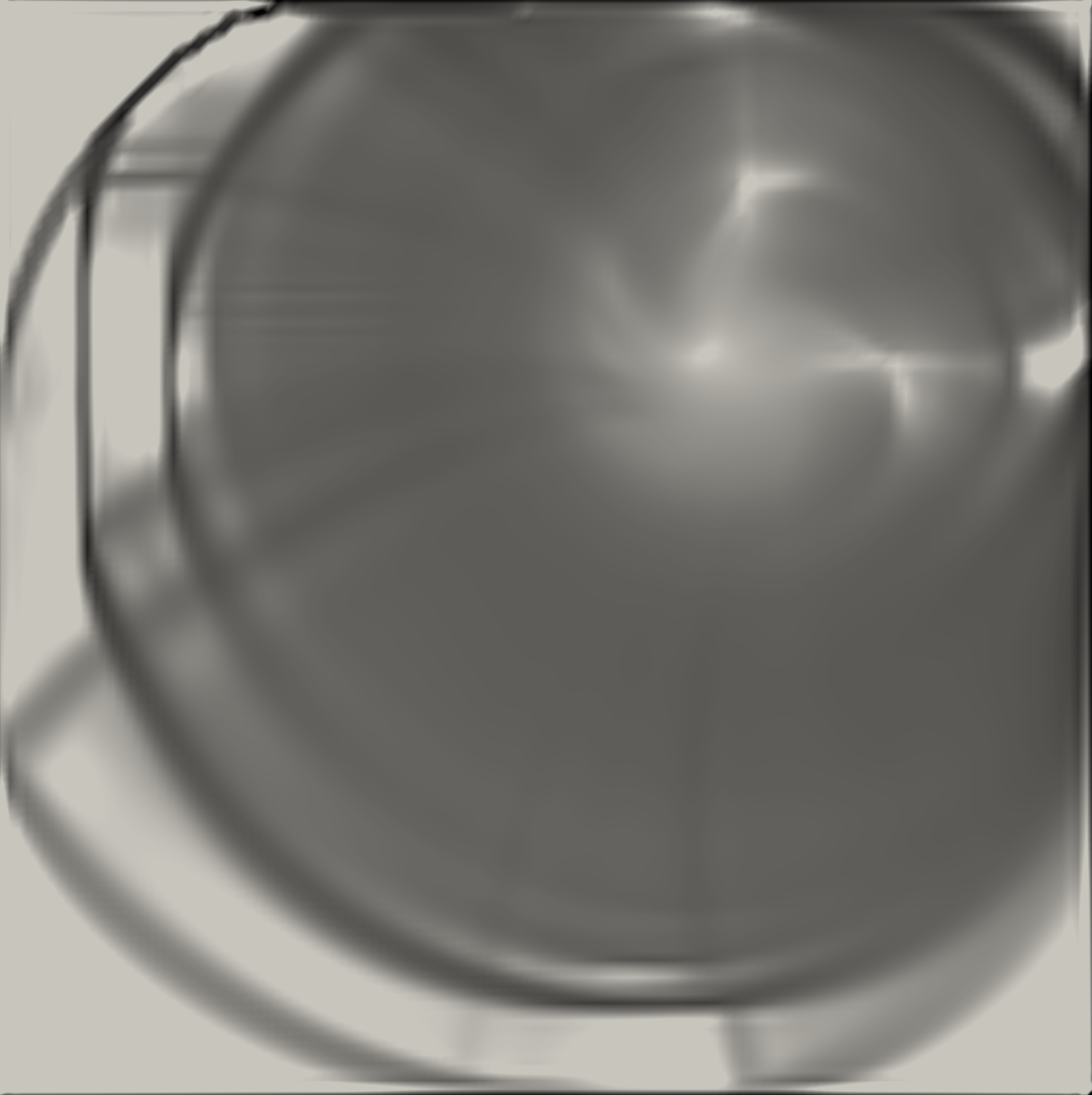}&
  \includegraphics[scale=0.04]{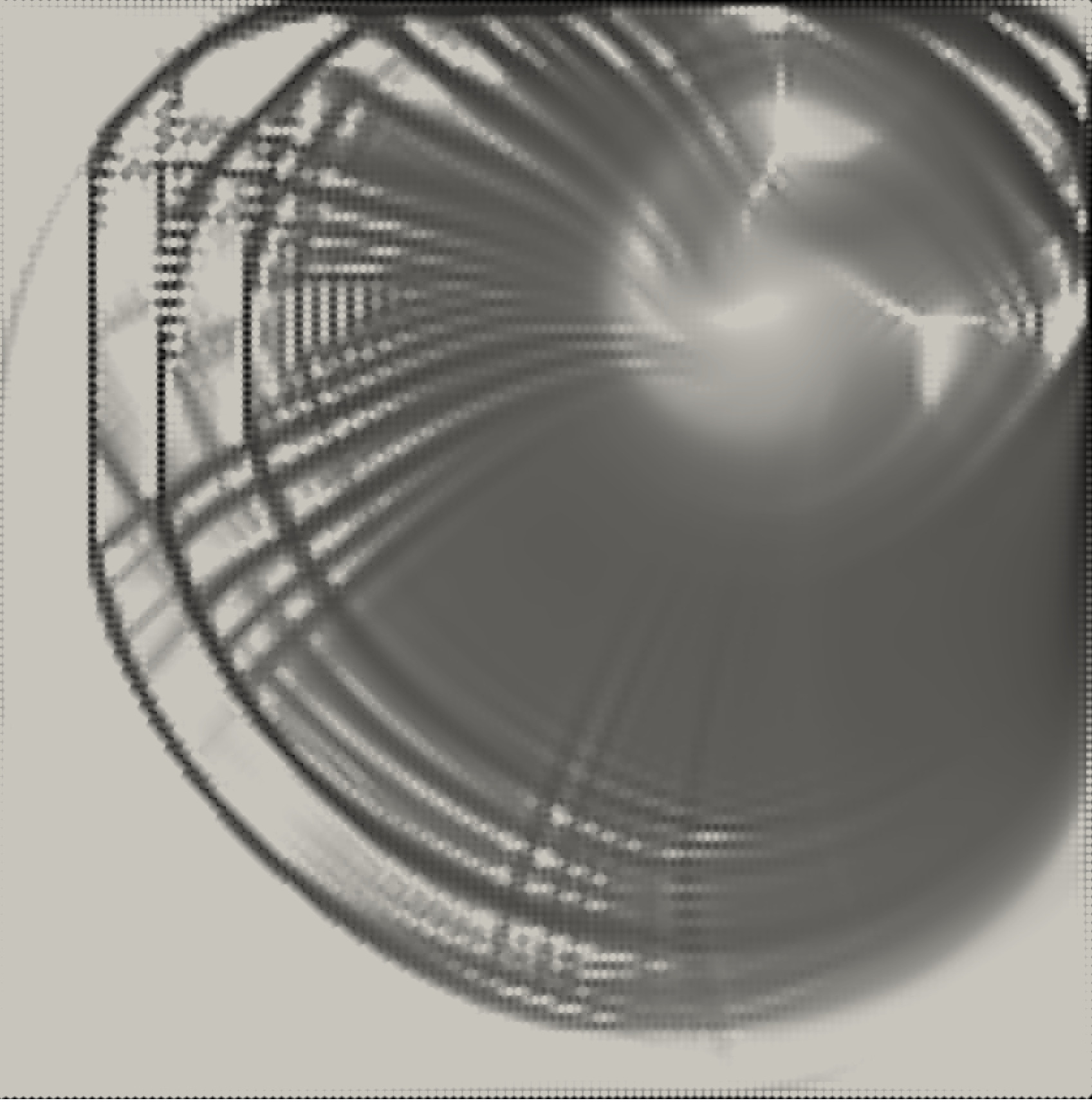}
 &\includegraphics[scale=0.04]{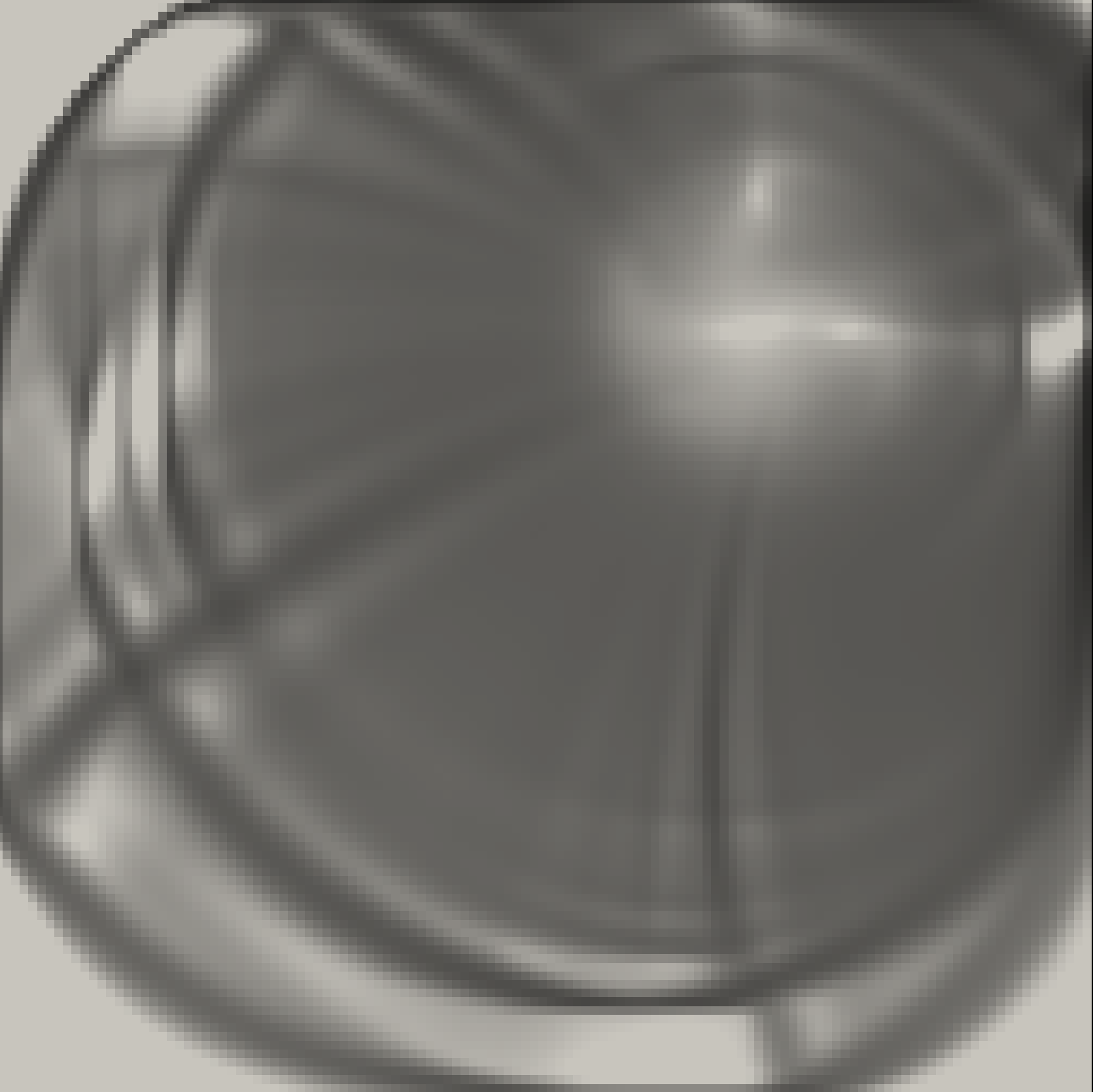}
 &\includegraphics[scale=0.04]{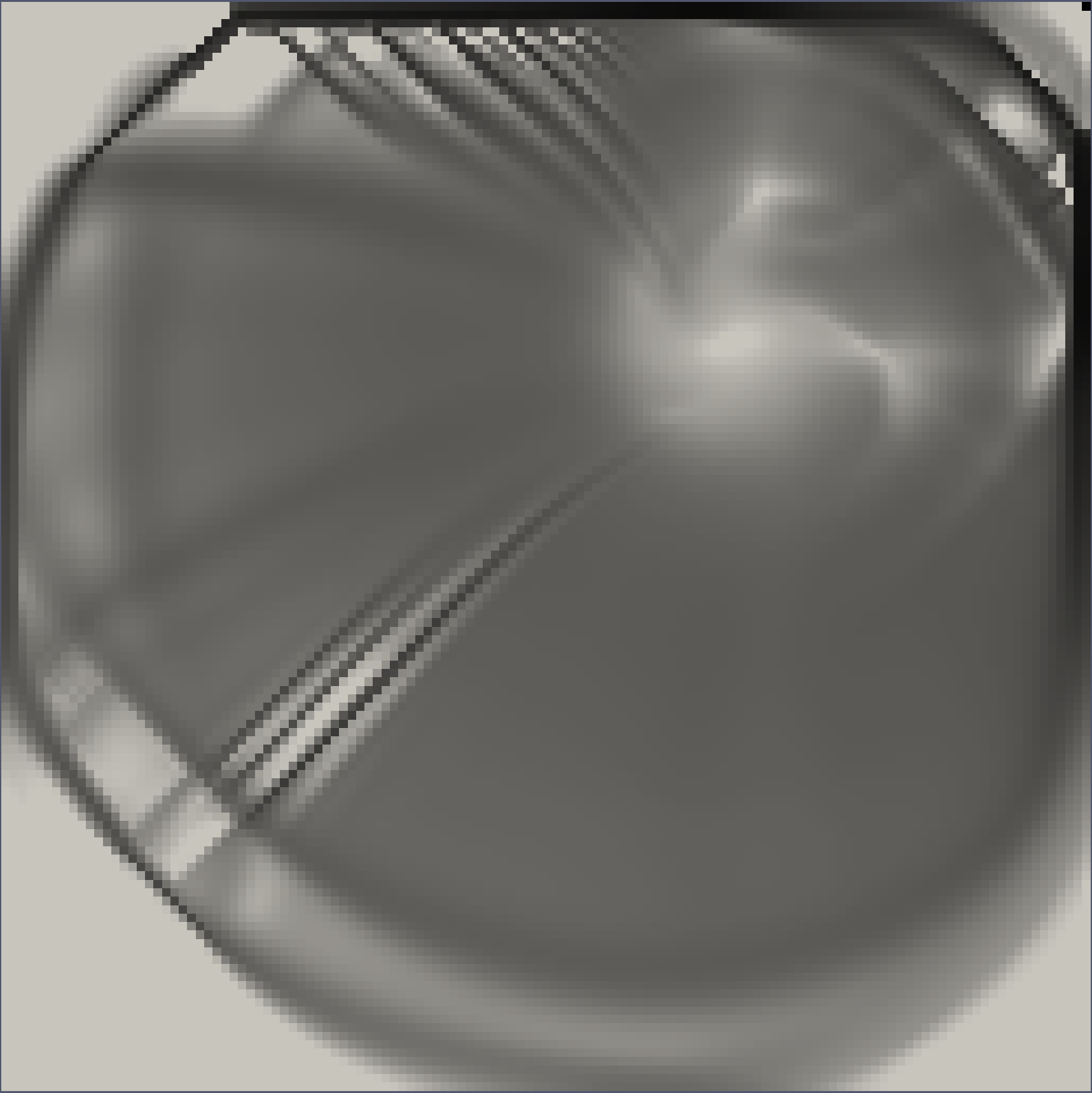} \\
     \rotatebox{90}{\phantom{abcde}2 km}&   \includegraphics[scale=0.04]{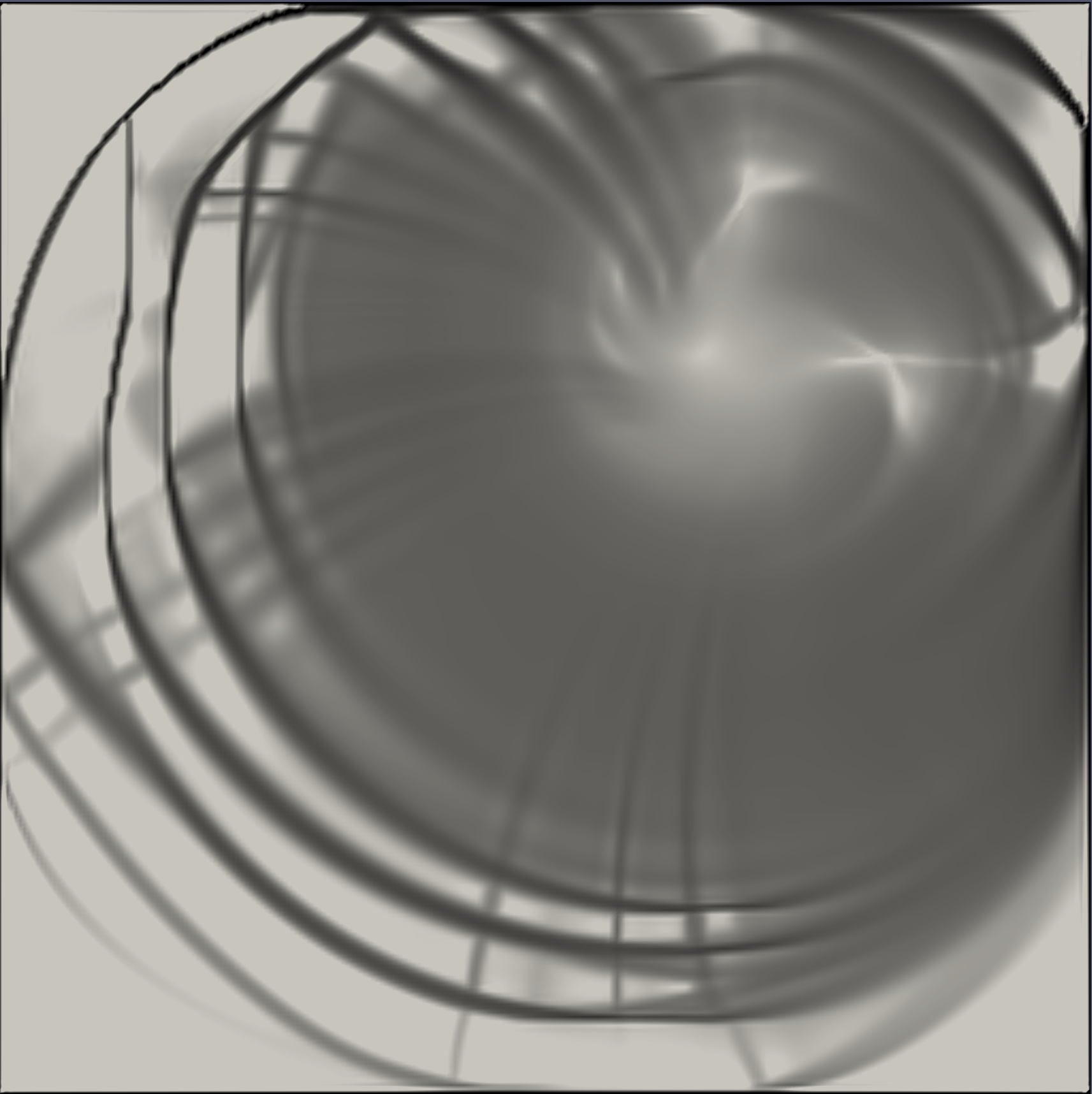}
 &\includegraphics[scale=0.05]{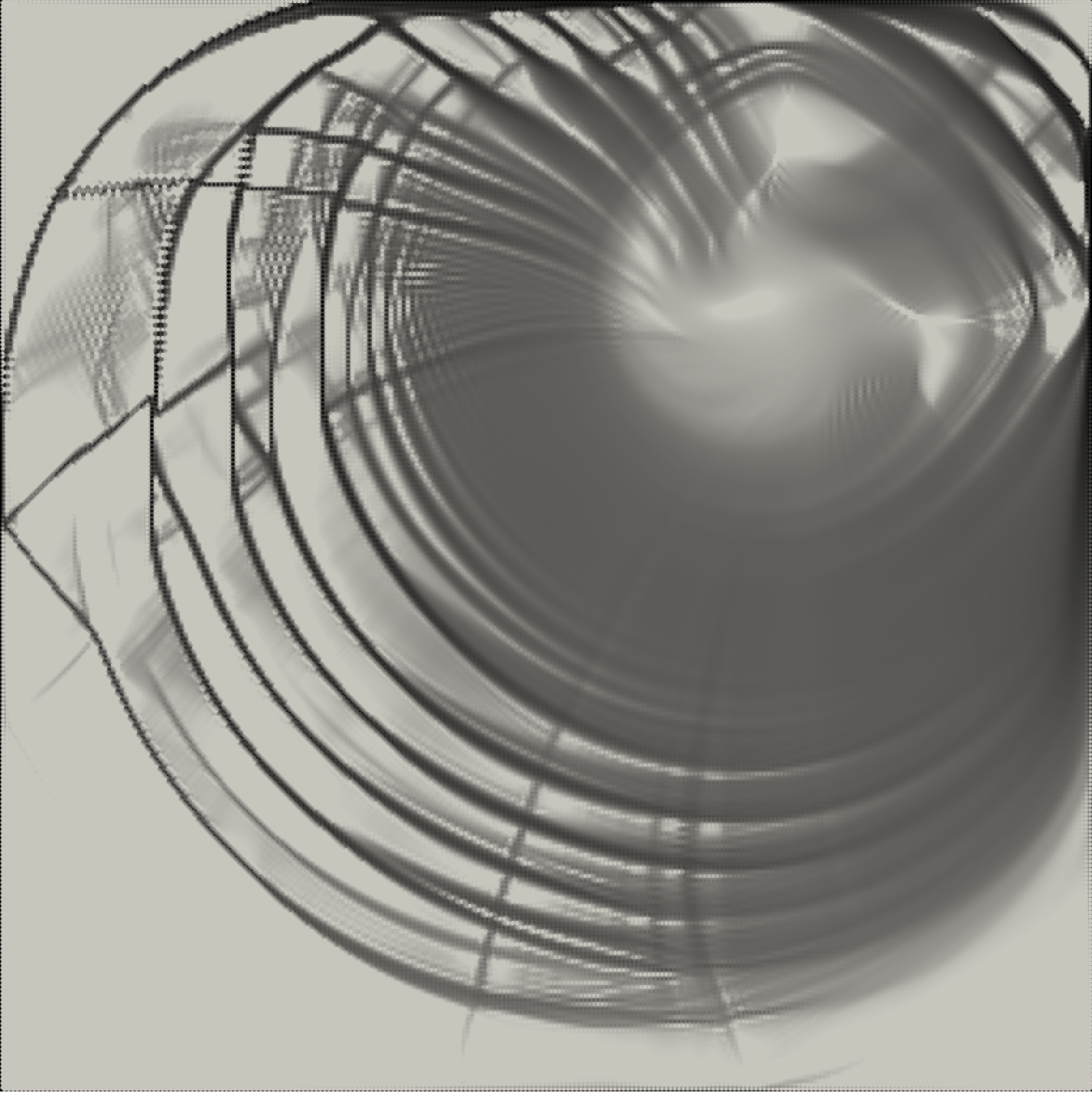}
 &\includegraphics[scale=0.04]{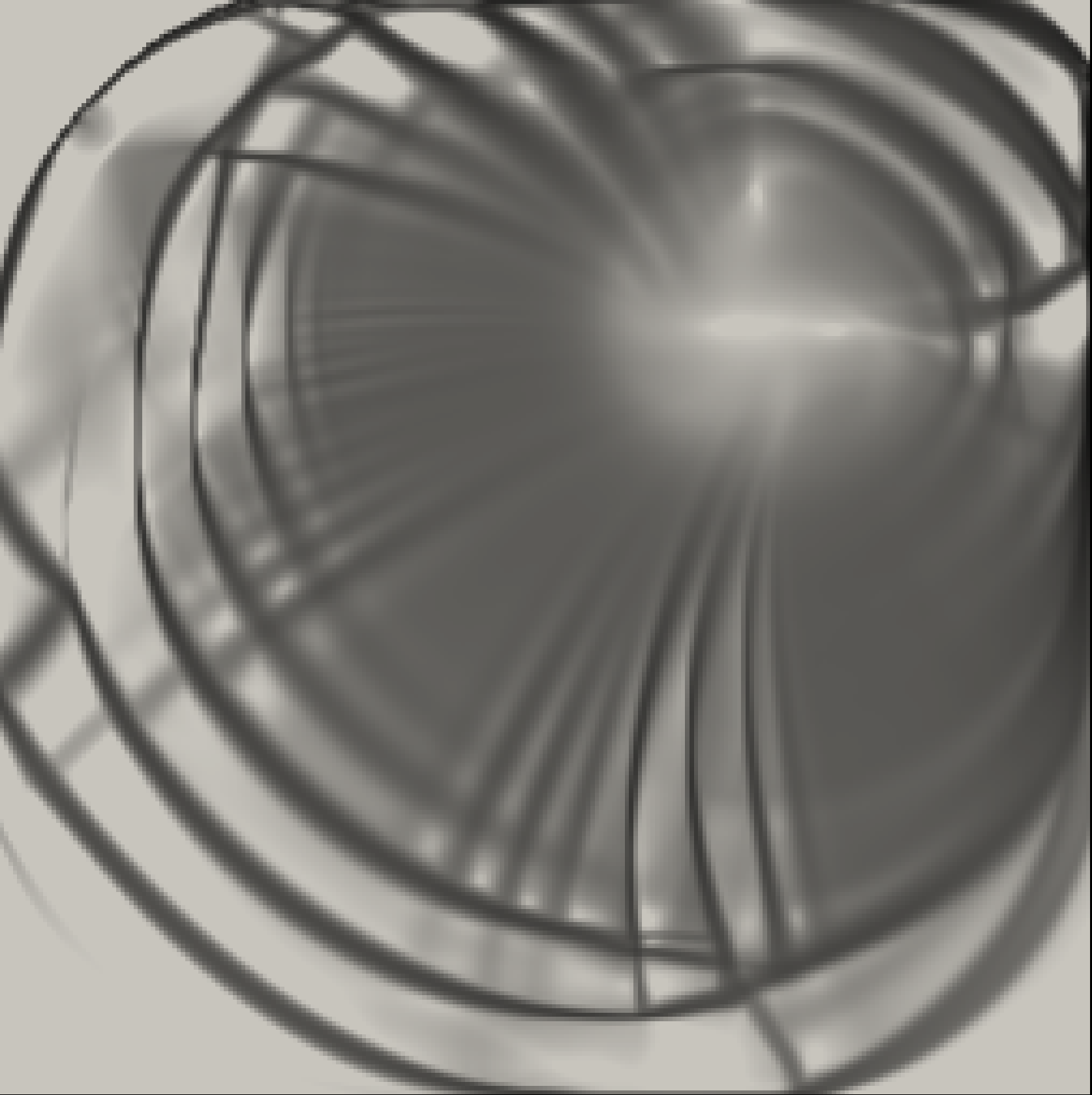}
 &\includegraphics[scale=0.04]{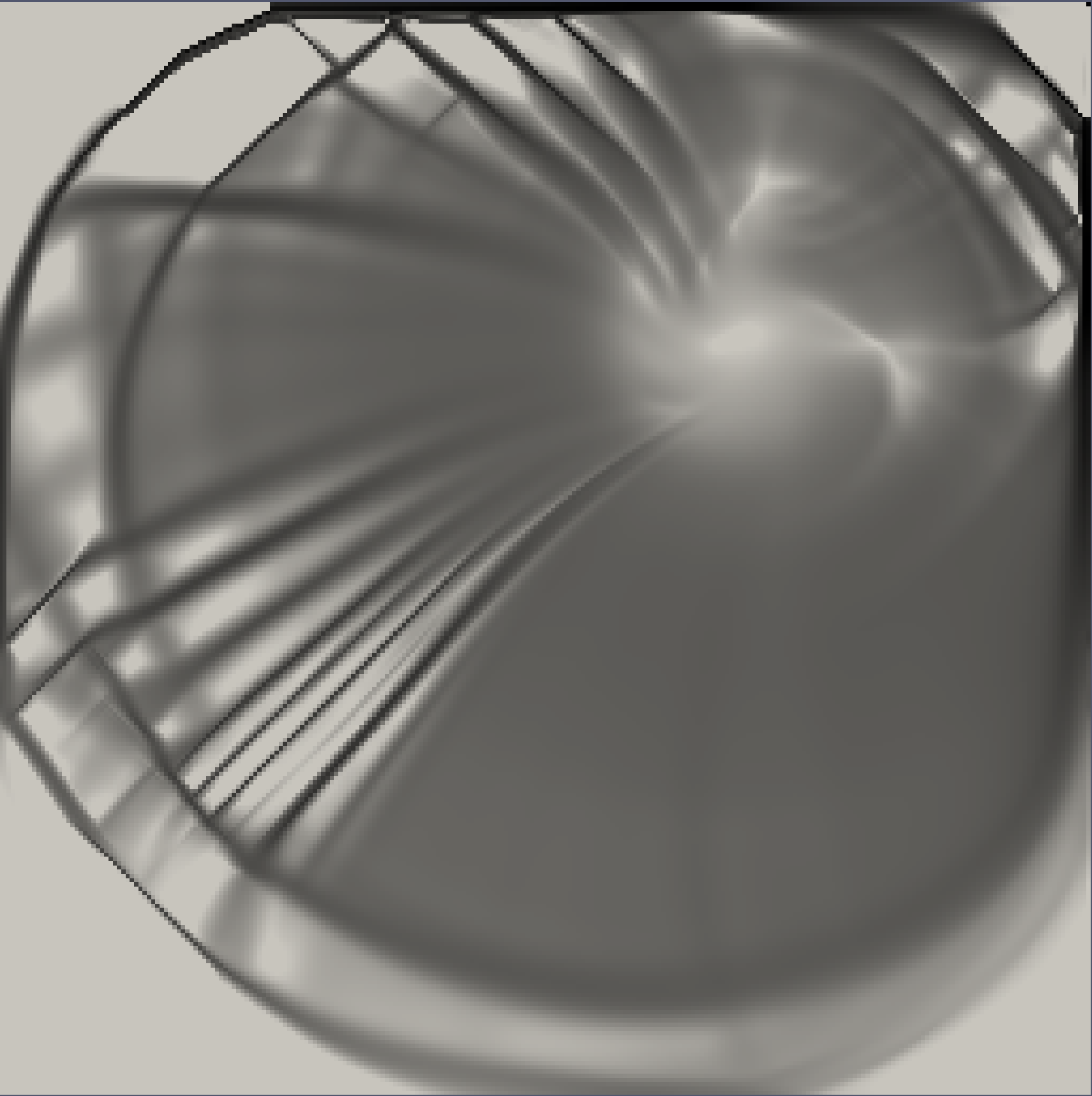} \\
& \includegraphics[scale=0.16]{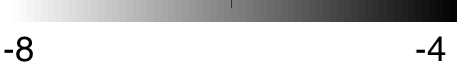}&  \includegraphics[scale=0.16]{pics_CR/ICON_sc}&
 \includegraphics[scale=0.16]{pics_CR/ICON_sc}&
 \includegraphics[scale=0.16]{pics_CR/ICON_sc} \\
     \end{tabular}
  \end{center}
  \begin{center}
   \caption{The first and second row show the sea ice concentration. The third and fourth row present the shear deformation calculated on a quadrilateral mesh. 
   \label{fig:Bgrid}}
  \end{center}
\end{figure}

Next we compare the A-grid, B-grid and CD-grid type approximations obtained with Gascoigne to the C-grid type discretization computed in the MITgcm model. The MITgcm configuration differs from Gascoigne's B-grid and CD-grid type setups by the placement and discretization of the sea ice velocity and the choice of advection scheme (second order central differences with a superbee flux limiter). We repeat the experiment with the MITgcm sea ice module with extreme choices of advection schemes, a first-order upwind scheme and a seventh-order monotonicity preserving advection scheme, and confirm that the advection scheme is not important in our context (results not shown). 
Thus, we attribute the difference in the sea ice concentration (Figure~\ref{fig:Cgrid}) and shear deformation (Figure~\ref{fig:Bgrid}) to the different variable staggering and in part to the finite element vs. finite volume discretizations of the sea ice velocity.

CICE uses a B-grid type staggering such as the Q1-Q0 finite element discretization in Gascoigne. The sea ice thickness and concentration are advected with the incremental remapping scheme \citep{Liscomb2004}. We also ran the test case in CICE with an upwind advection, but did not observe significant changes in the number and width of the resolved LKFs.
Figure~\ref{fig:Bgrid}  shows that the LKFs in the CICE discretization are wider than the LKFs in the B-grid like  Gascoigne approximation. 
Furthermore, the simulated sea ice concentration is more diffusive. We attribute the difference in the sea ice concentration and shear deformation field shown in Figure~\ref{fig:Bgrid} to the discretization of the velocity field and stress tensor. 
\begin{figure}[t]
  \begin{center}
         \includegraphics[scale=0.5]{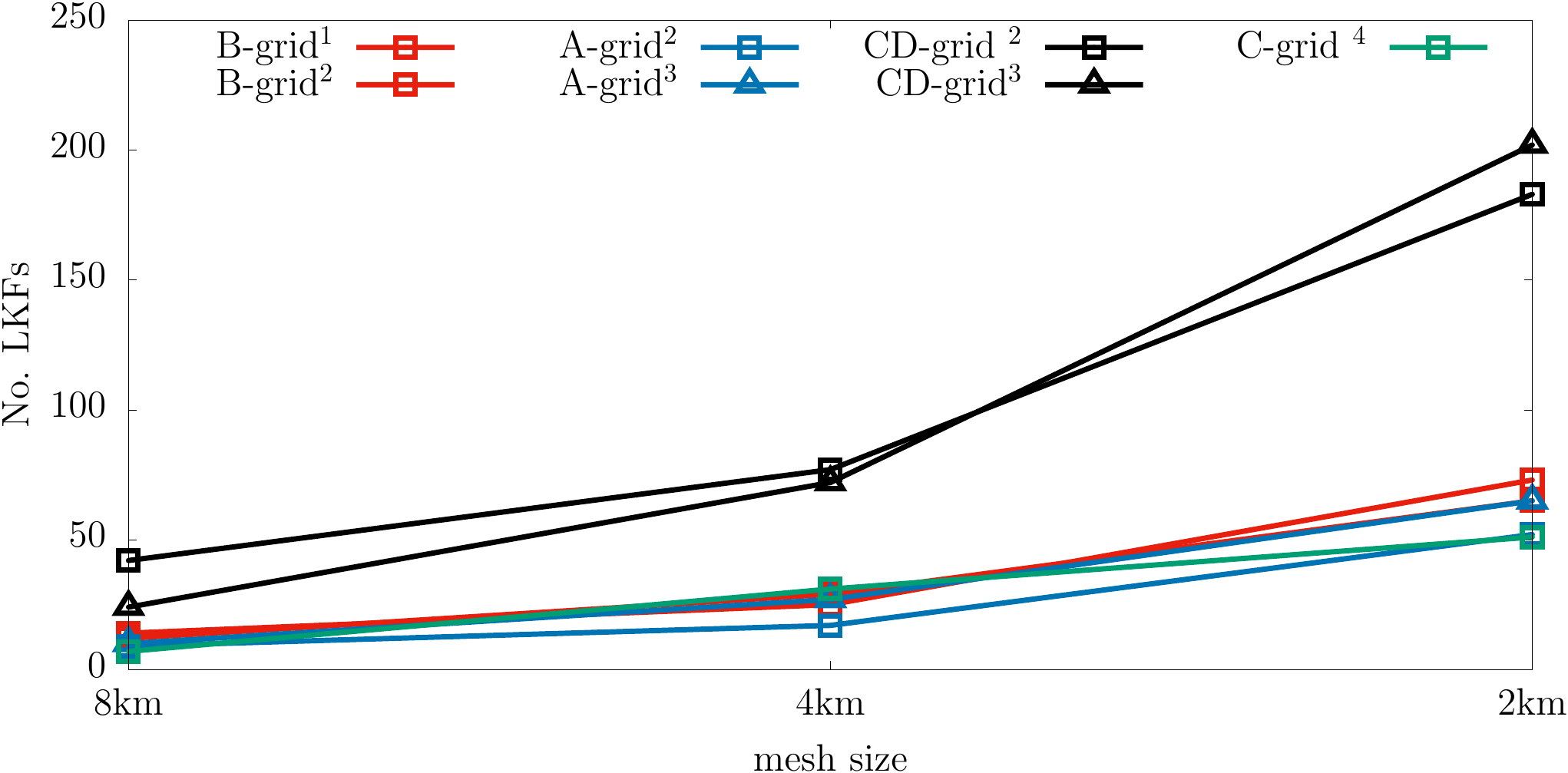}
  \end{center}
  \begin{center}
   \caption{ The number of the detected LKFs from the shear deformation field presented in Figure~\ref{fig:Bgrid} and \ref{fig:comp_quads_tri} on quadrilateral $\Box$ and triangular $\triangle$ meshes. The subscript 1, 2, 3 and 4 refer to the simulations carried out in the framework of CICE, Gascoigne, FESOM and MITgcm respectively. \label{fig:vgl_LKF}}
  \end{center}
\end{figure}
%
\begin{figure}[t]
  \begin{center}
    \begin{tabular}{c c c  |c  c }
     &\multicolumn{2}{c|}{triangular grid}&\multicolumn{2}{c}{quadrilateral grid}\\
     \hline
     &A-grid& CD-grid & A-grid & CD-grid \\
     &(FESOM)& (FESOM) &(Gascoigne) & (Gascoigne)\\
      \rotatebox{90}{\phantom{abcde}4 km}&  \includegraphics[scale=0.05]{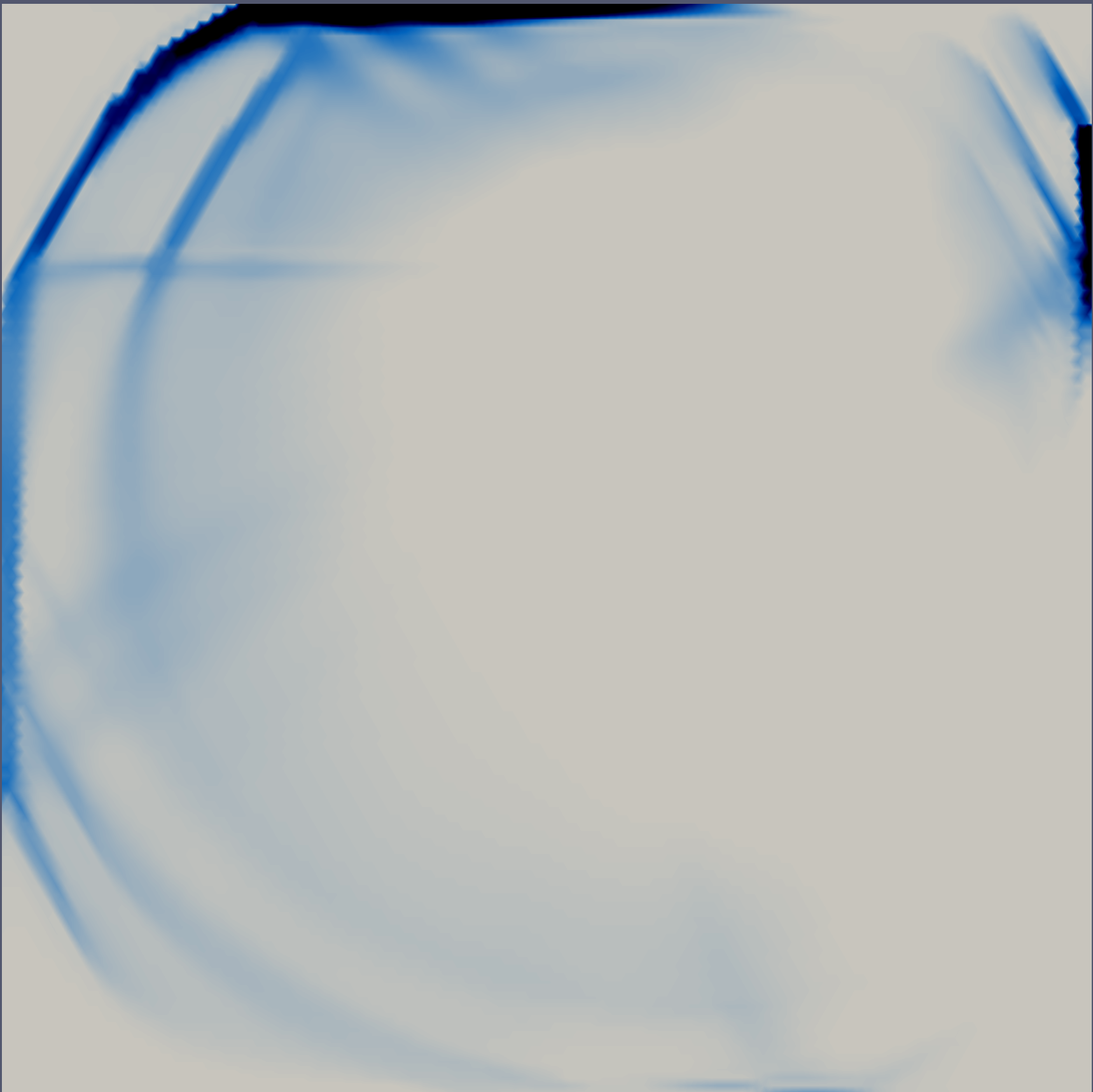}& 
      \includegraphics[scale=0.05]{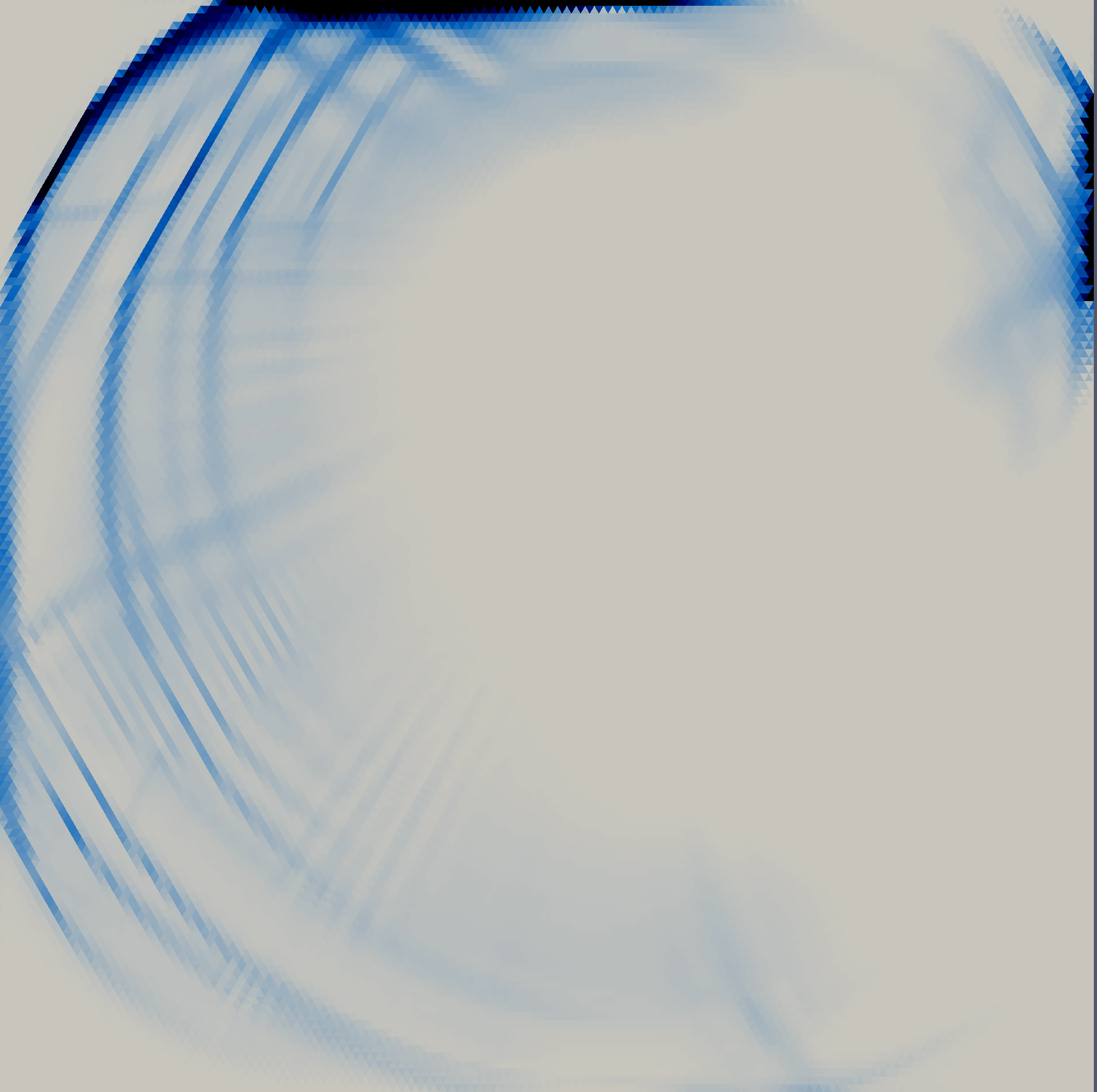}&
      \includegraphics[scale=0.049]{pics_CR/Q1-Q1_4km_A_2min.png}& 
      \includegraphics[scale=0.05]{pics_CR/CR_Q0_A_4km_02.png}\\
       \rotatebox{90}{\phantom{abcde}2 km}&
      \includegraphics[scale=0.05]{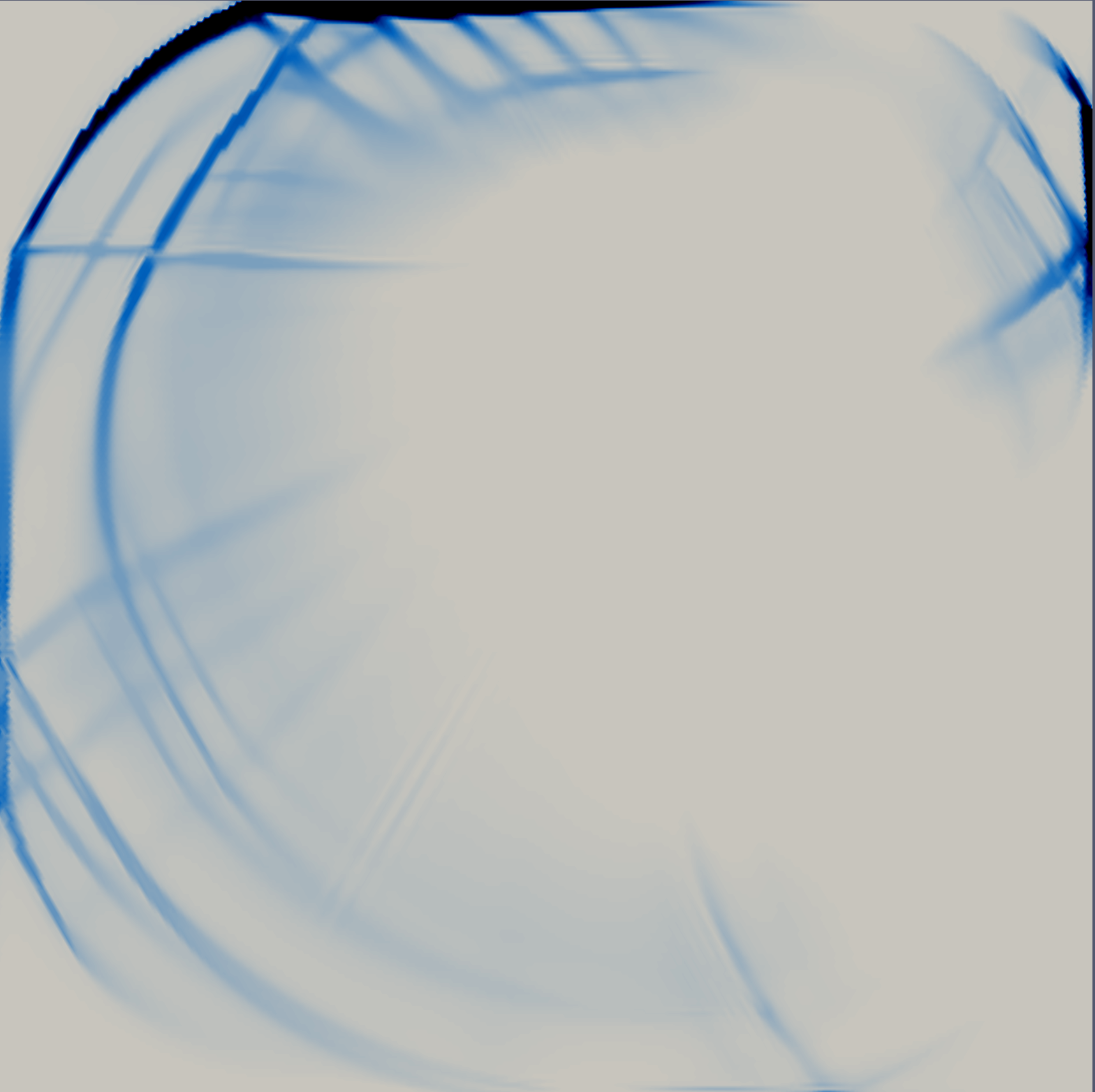}& 
         \includegraphics[scale=0.115]{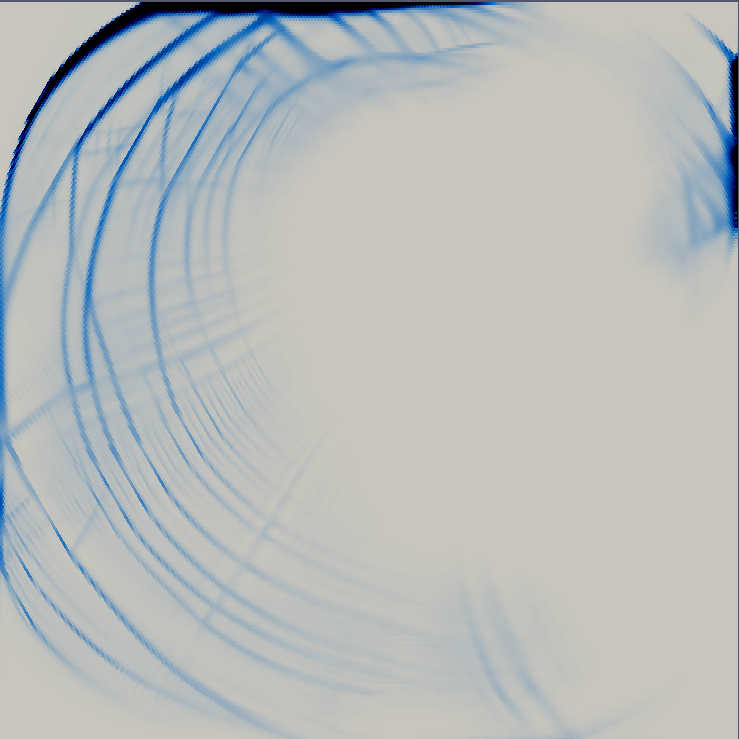}& 
         \includegraphics[scale=0.065]{pics_CR/Q1_Q1_A_2km_neu.png}&
        \includegraphics[scale=0.05]{pics_CR/CR_Q0_A_2km_01.png}
       \\
        & \includegraphics[scale=0.16]{pics_CR/ICON_A_sc}& 
          \includegraphics[scale=0.16]{pics_CR/ICON_A_sc}&    
          \includegraphics[scale=0.16]{pics_CR/ICON_A_sc}&    
         \includegraphics[scale=0.16]{pics_CR/ICON_A_sc}\\
     \end{tabular}
  \end{center}
  \begin{center}
    \caption{Sea ice concentration on triangular and quadrilateral meshes using the same number of vertices. Therefore the edge length in the triangular case is 8.6 km, 4.3 km and 2.15 km. The triangular mesh has around double the number of cells and 1.5 more edges compared to the quadrilateral mesh. 
    \label{fig:comp_quads_tri} }
  \end{center}
\end{figure}

Both sea ice concentration and deformation fields (Figures~\ref{fig:Cgrid} and~\ref{fig:Bgrid}) show that the CD-grid type approximation on the 4\,km mesh resolves more LKFs than the A-grid, the B-grid, and the C-grid like discretization on the 2\,km mesh. The number of detected LKFs in the deformation fields confirms this observation (Figure~\ref{fig:vgl_LKF}). This result can plausibly be attributed to the fact that the CD-type discretization has {twice} as many velocity degrees of freedom (DoF) compared to an A-grid, B-grid or C-grid type staggering.\footnote{The degrees of freedom are given by the nodes per cell divided by the number of neighbouring cells sharing each node. The CD-grid has 4N DoF (with N the number of grid nodes), whereas A, B, C-grid contains 2N DoF.}
Furthermore, the CD-grid approximation resolves almost the same number of LKFs on the 8\,km mesh compared to A, B and C-grids on the 2\,km grid, see Figure~\ref{fig:vgl_LKF}.  


\subsection{Comparison of quadrilateral and triangular grids}\label{sec:comp_quad_tri}

In this section we compare the approximation of sea ice dynamics on quadrilateral and triangular grids with respect to their ability to resolve LKFs. The triangular grids contain twice as many cells as the quadrilateral grids, but the same number of vertices. 
The discretization on the triangular grids is done in the framework of the sea ice model FESOM using an A-grid type and CD-grid type staggering. 

Figure~\ref{fig:comp_quads_tri} compares the approximations of the sea ice concentration with A-grid and CD-grid type discretizations on quadrilateral grids to those on triangular grids. 
 We observe that the A-grid and CD-grid like discretization are qualitatively similar on both grids. 
 
 We conclude that, given the  same number of vertices, the number of resolved LKFs is more sensitive to the velocity staggering than to the use of quadrilateral or triangular grids. This is supported by the results of the detection algorithm (Figure~\ref{fig:vgl_LKF}). The number of LKFs is similar for the CD-type discretization on the structured quadrilateral and triangular meshes. 
 The same is true for the A-grid type discretizations. 
  Figure~\ref{fig:vgl_LKF} also shows that traditionally used Arakawa B-grid and C-grid discretizations on structured quadrilateral meshes resolve quantitatively the same number of LKFs as an A-grid like approximation on triangular meshes.
  
   Furthermore, we compared the probability density function of shear deformation field and found that the magnitude of the shear deformation is nearly identical among the different discretizations (not shown).  

\subsection{Triangles}\label{sec:tri}
 \begin{figure}[t]
  \begin{center}
    \begin{tabular}{c c c c c}
     & A-grid& B-grid & CD-grid & CD-grid\\
    &(FESOM)& (FESOM)& (FESOM) &(ICON)\\
 \rotatebox{90}{\phantom{abcde}4 km}&       \includegraphics[scale=0.05]{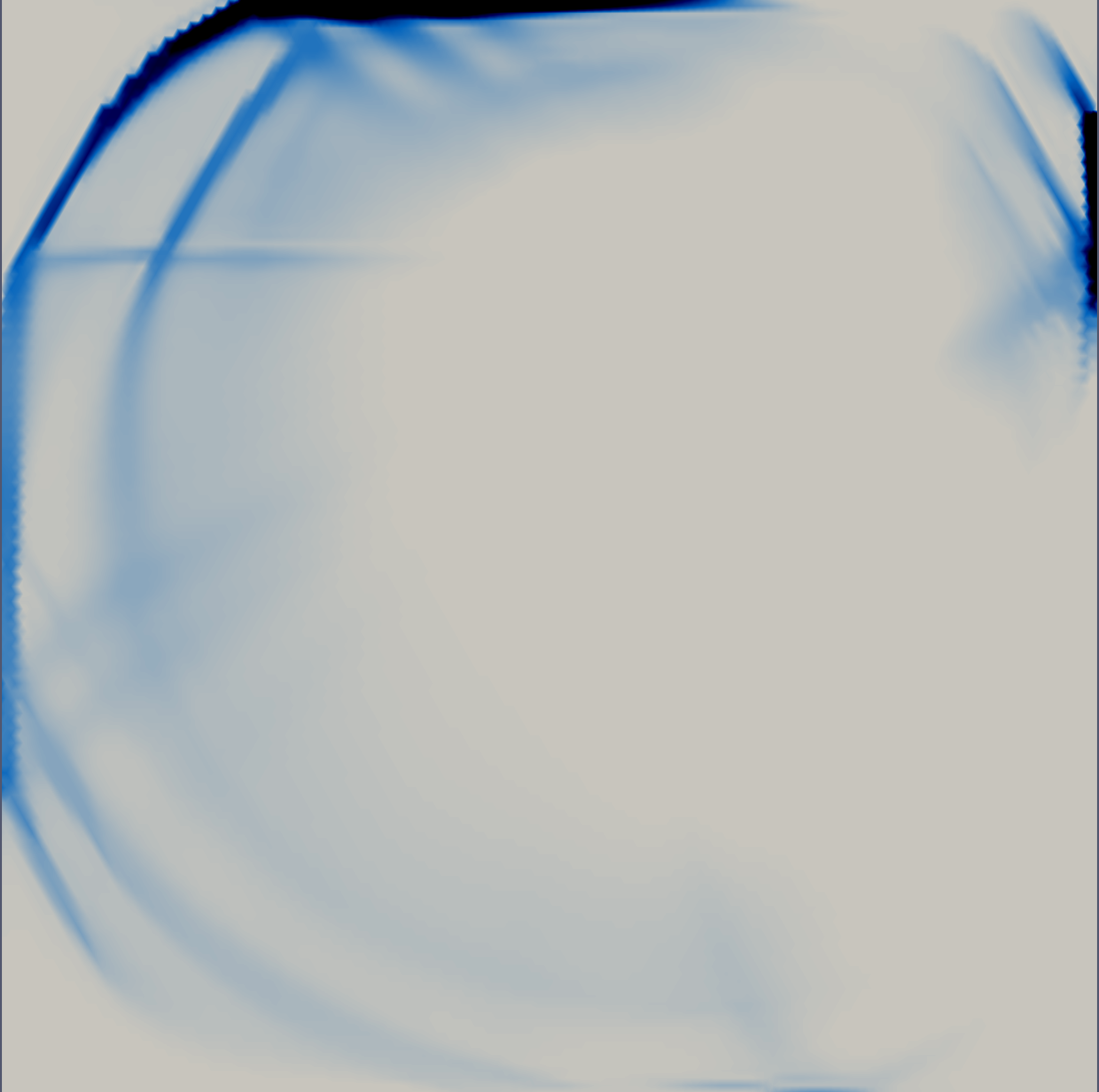}& 
       \includegraphics[scale=0.05]{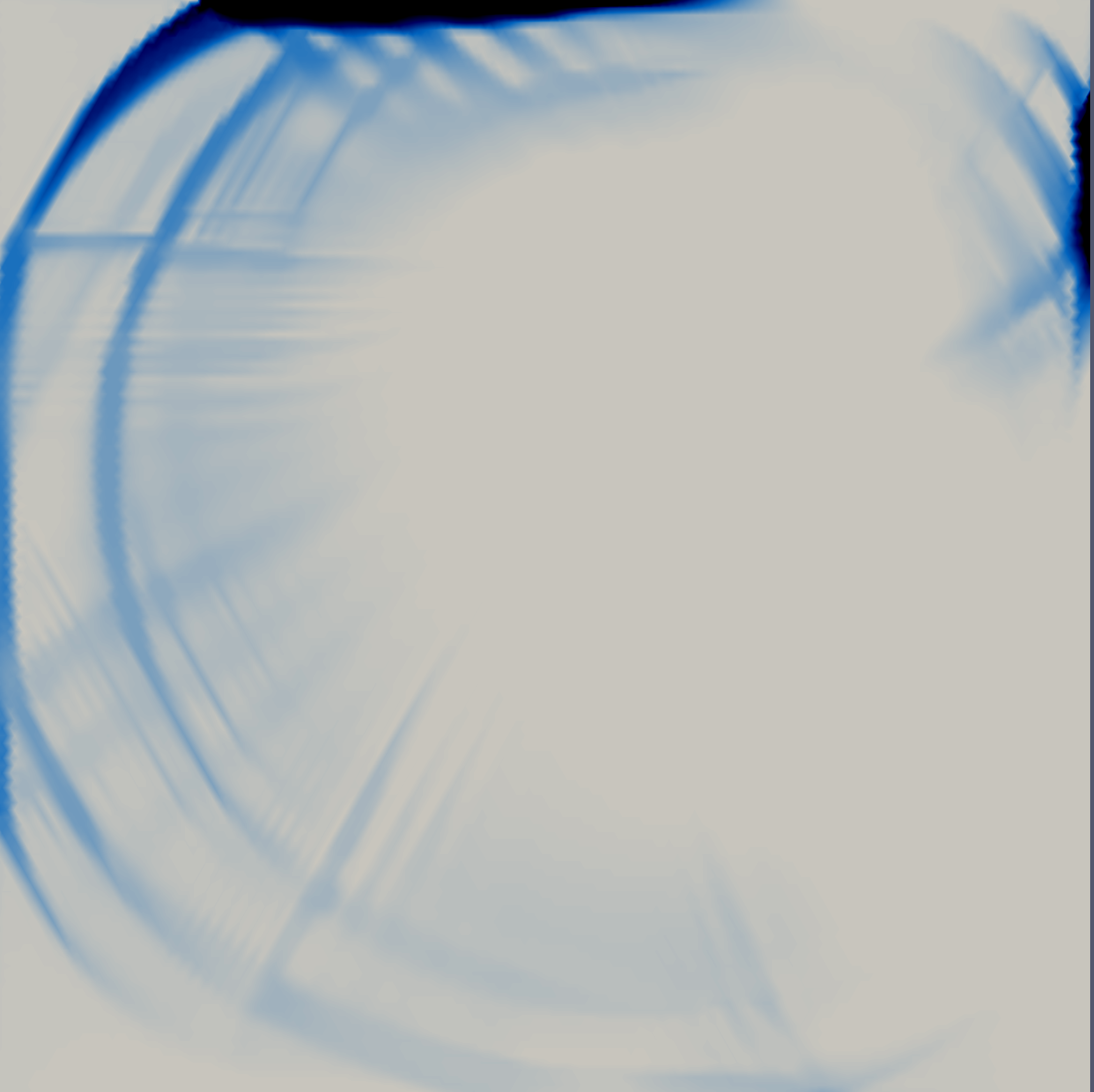}&
       \includegraphics[scale=0.05]{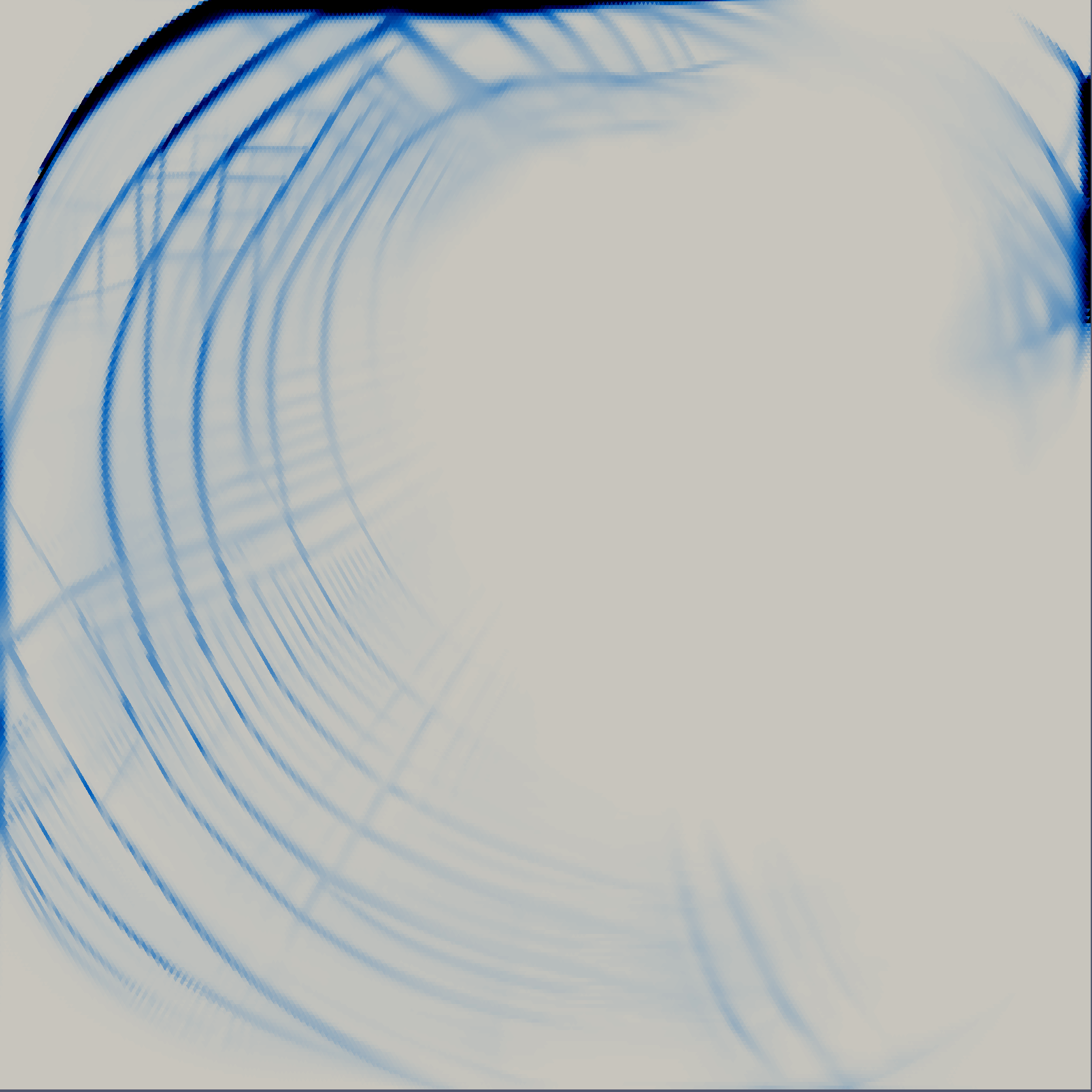}&
      \includegraphics[scale=0.05]{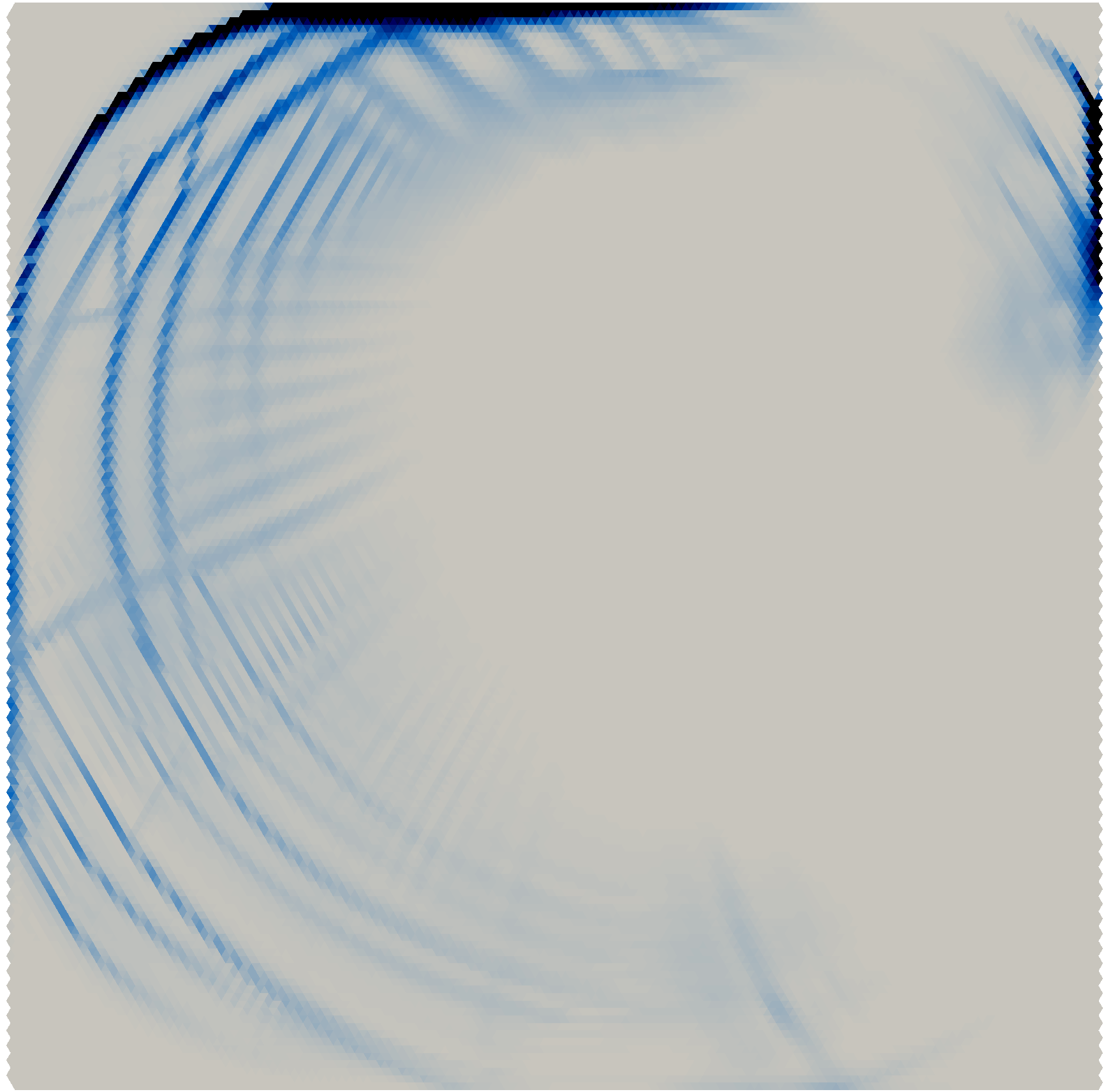}\\
 \rotatebox{90}{\phantom{abcdef}2 km}&
 \includegraphics[scale=0.05]{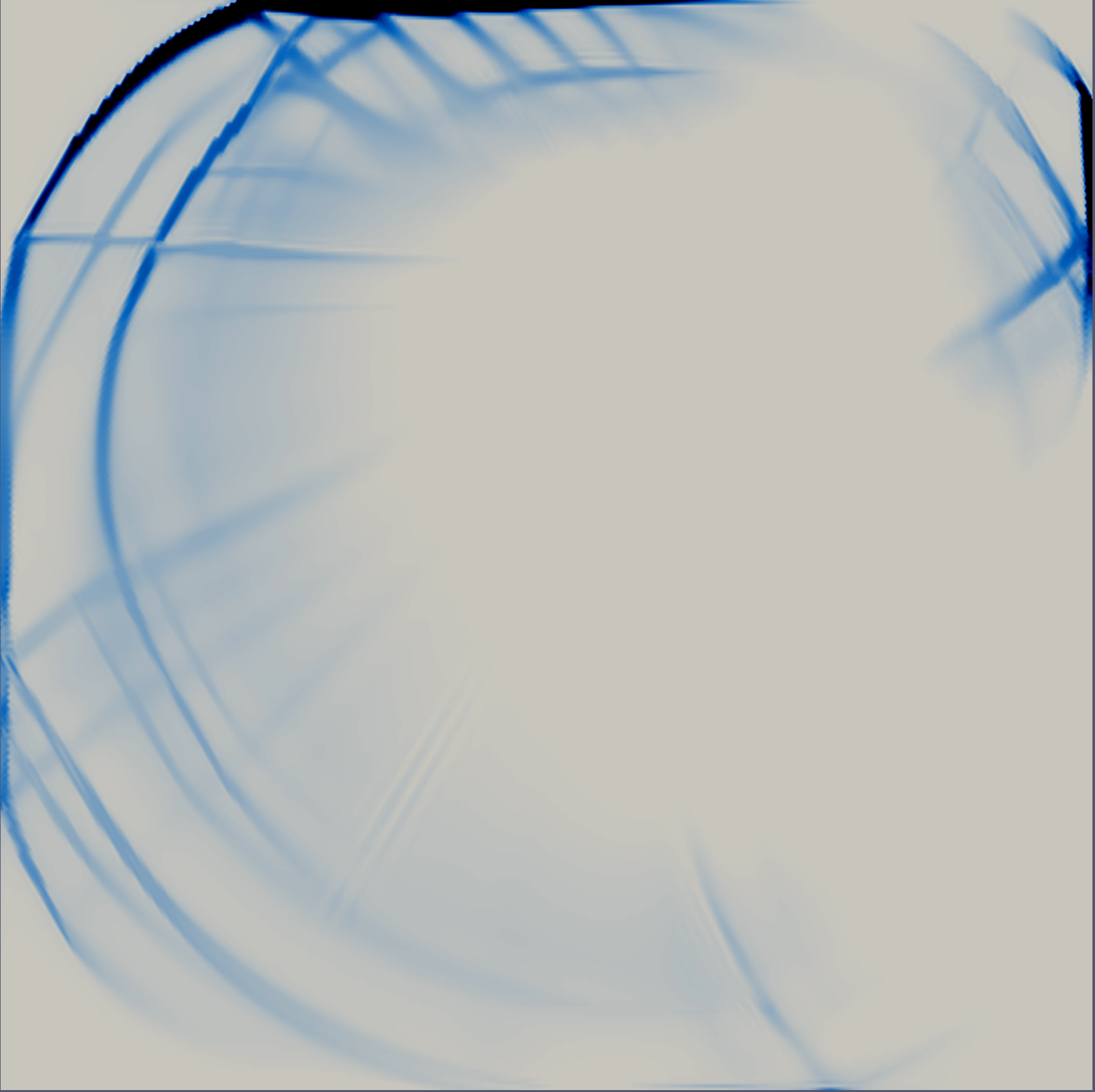}&
    \includegraphics[scale=0.05]{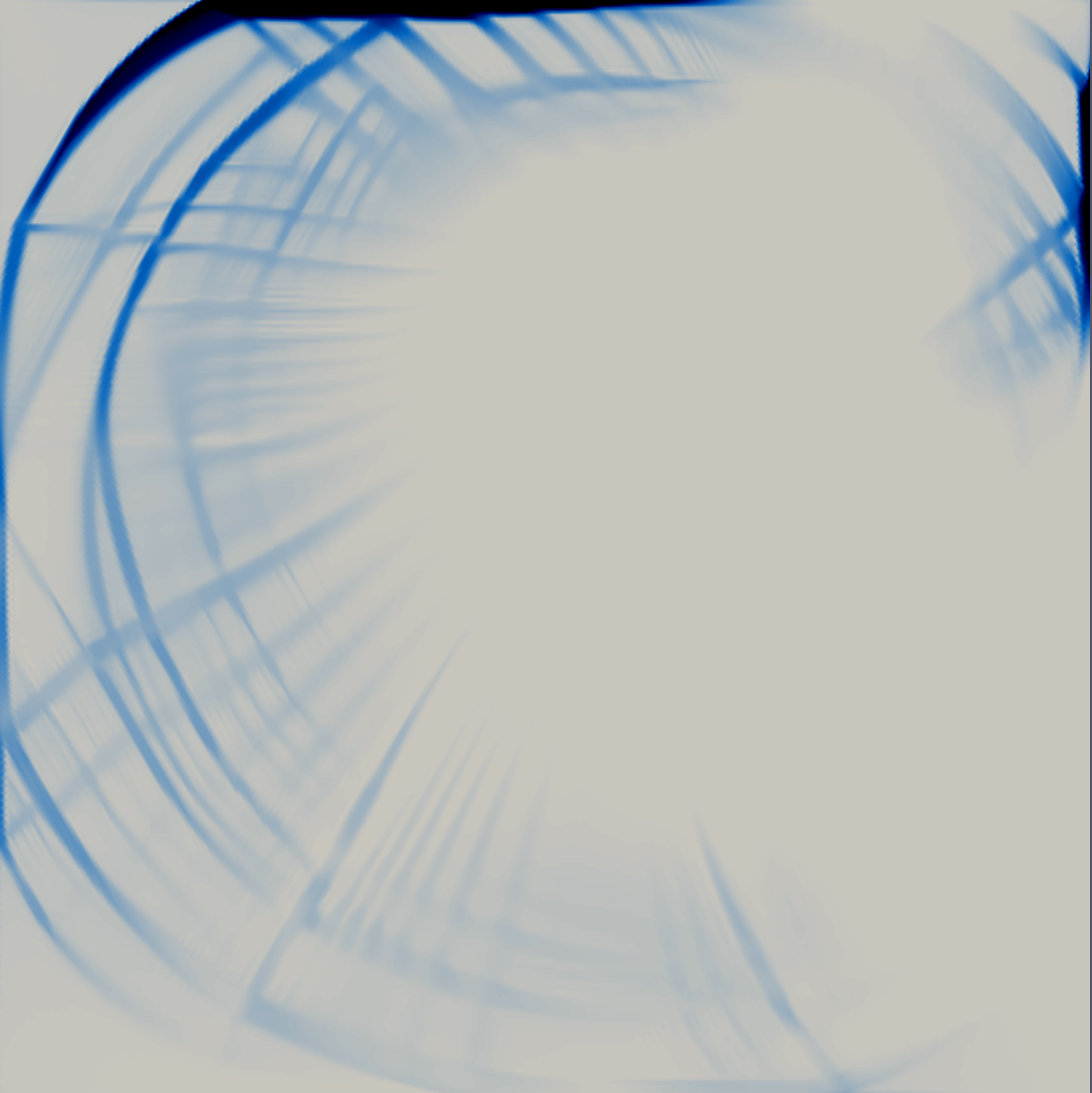}& 
      \includegraphics[scale=0.05]{pics_CR/A_A_CR_2km_2min.png}& 
         \includegraphics[scale=0.05]{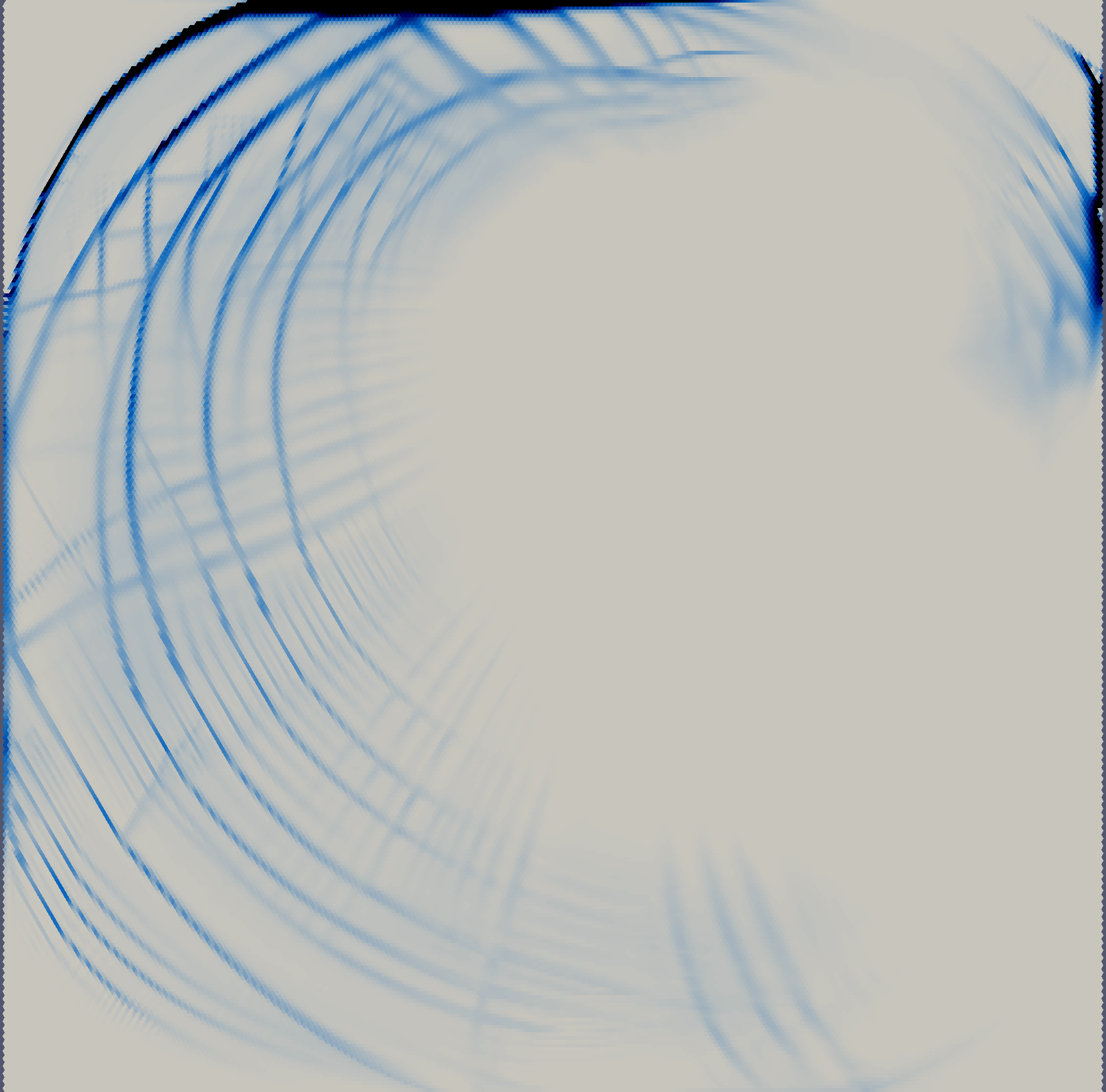}\\
  &       \includegraphics[scale=0.16]{pics_CR/ICON_A_sc}&
             \includegraphics[scale=0.16]{pics_CR/ICON_A_sc}&
         \includegraphics[scale=0.16]{pics_CR/ICON_A_sc}&
         \includegraphics[scale=0.16]{pics_CR/ICON_A_sc}\\
     \end{tabular}
  \end{center}
  \begin{center}
   \caption{Sea ice concentration computed on a triangular grid. The P1-P1 (A-grid), P1-P0 (B-grid), and the CR-P0 pair (CD-grid) is computed in FESOM. The CD-grid discretization in ICON  also uses the CR-P0 finite element. 
   \label{fig:triA} }
  \end{center}
\end{figure}

In this section we compare the discrete solution obtained with the A-grid, B-grid, and CD-grid type discretizations on triangular grids in the framework of the sea ice module of FESOM \citep{Danilov2015} and ICON \citep{Mehlmann2020}. 
To be consistent with the quadrilateral case in Section \ref{sec:quad} we use triangles with side lengths of 8\,km, 4\,km, and 2\,km.  

\begin{figure}[t]
\begin{center}
\begin{tabular}{c c c c c}
 &A-grid& B-grid & CD-grid & CD-grid \\
 &(FESOM)& (FESOM)&(FESOM)& (ICON)\\

%
 \rotatebox{90}{\phantom{abcde}4 km}&\includegraphics[scale=0.05]{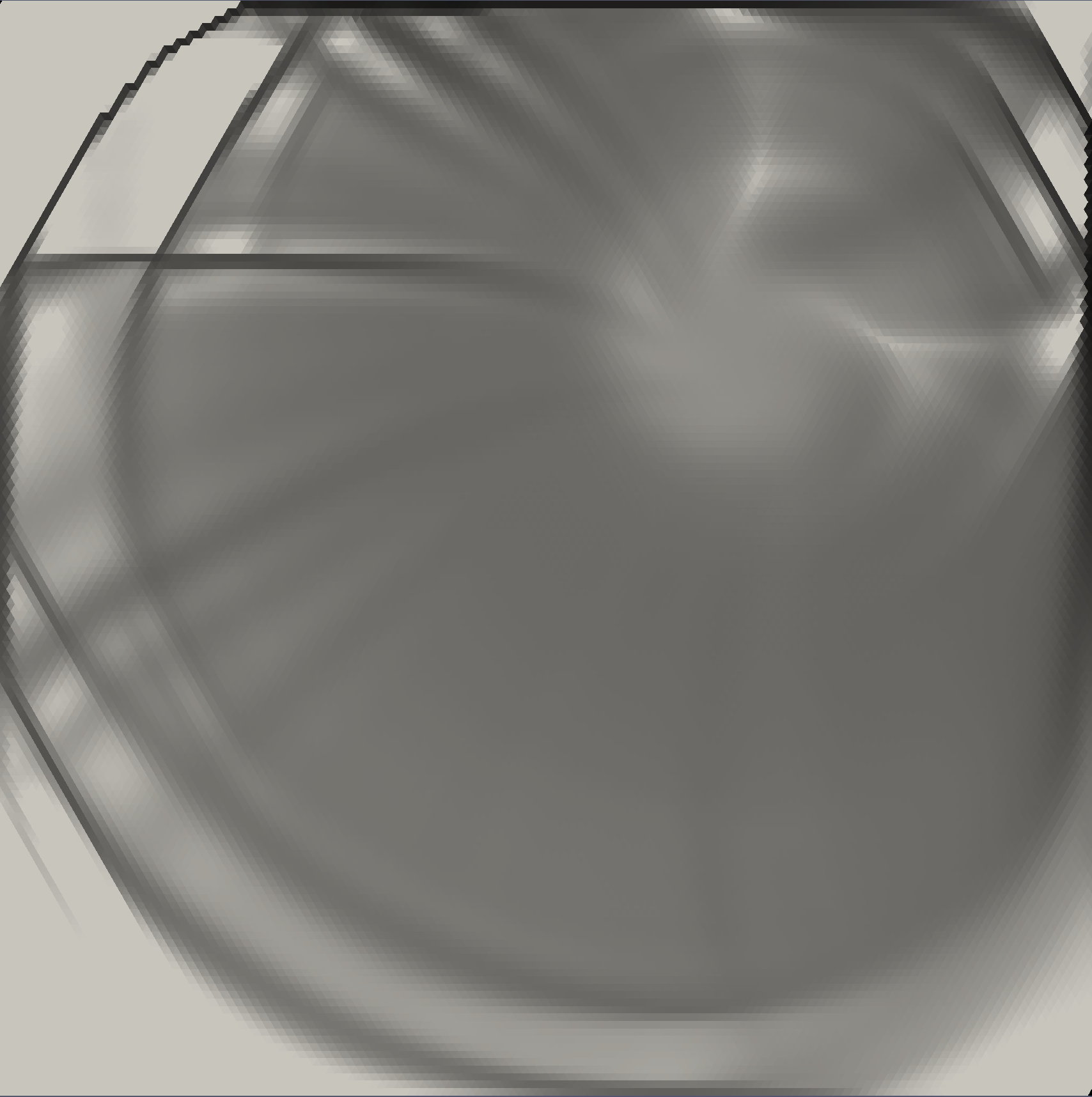}
 &\includegraphics[scale=0.05]{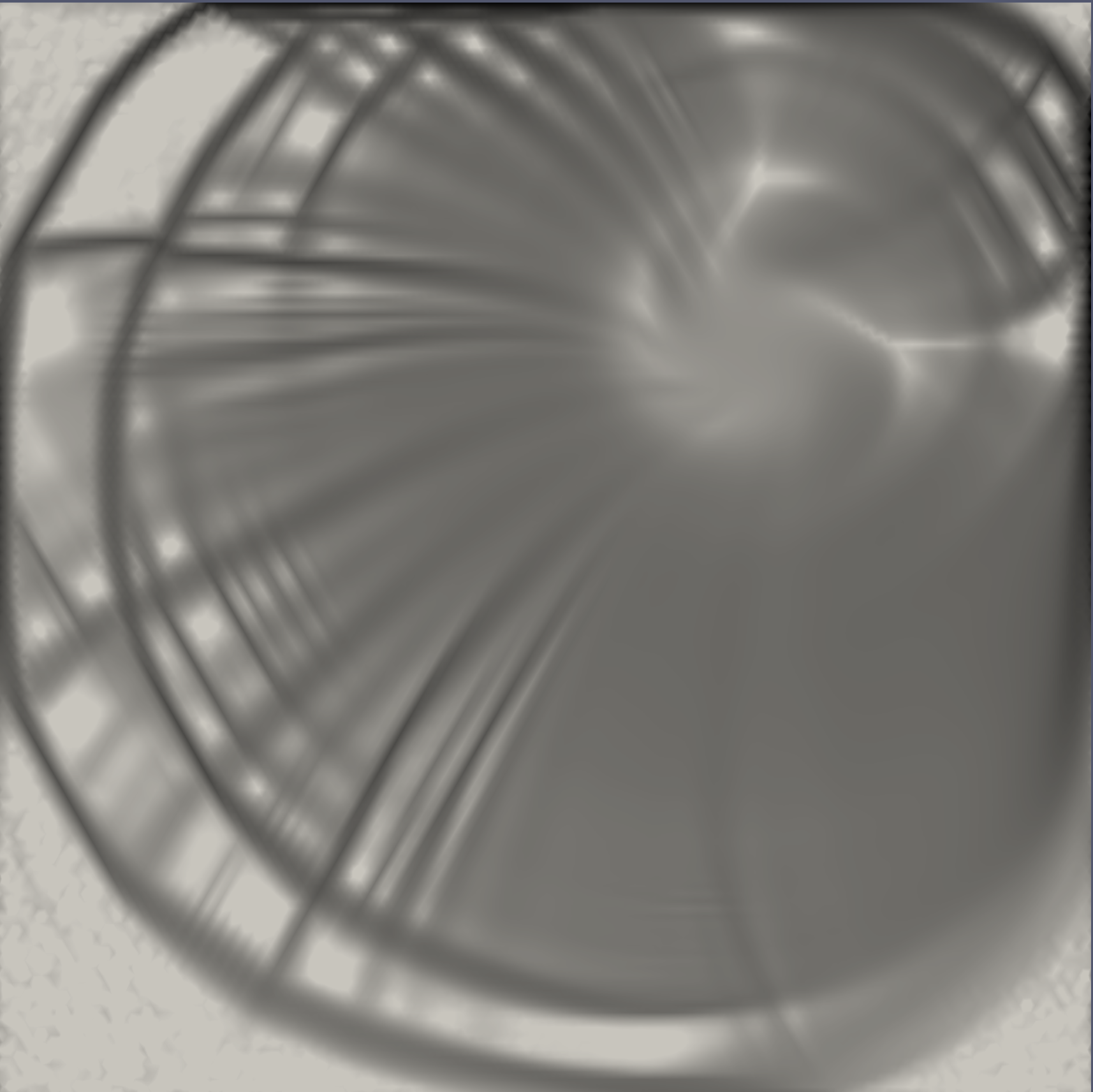}
 &\includegraphics[scale=0.06]{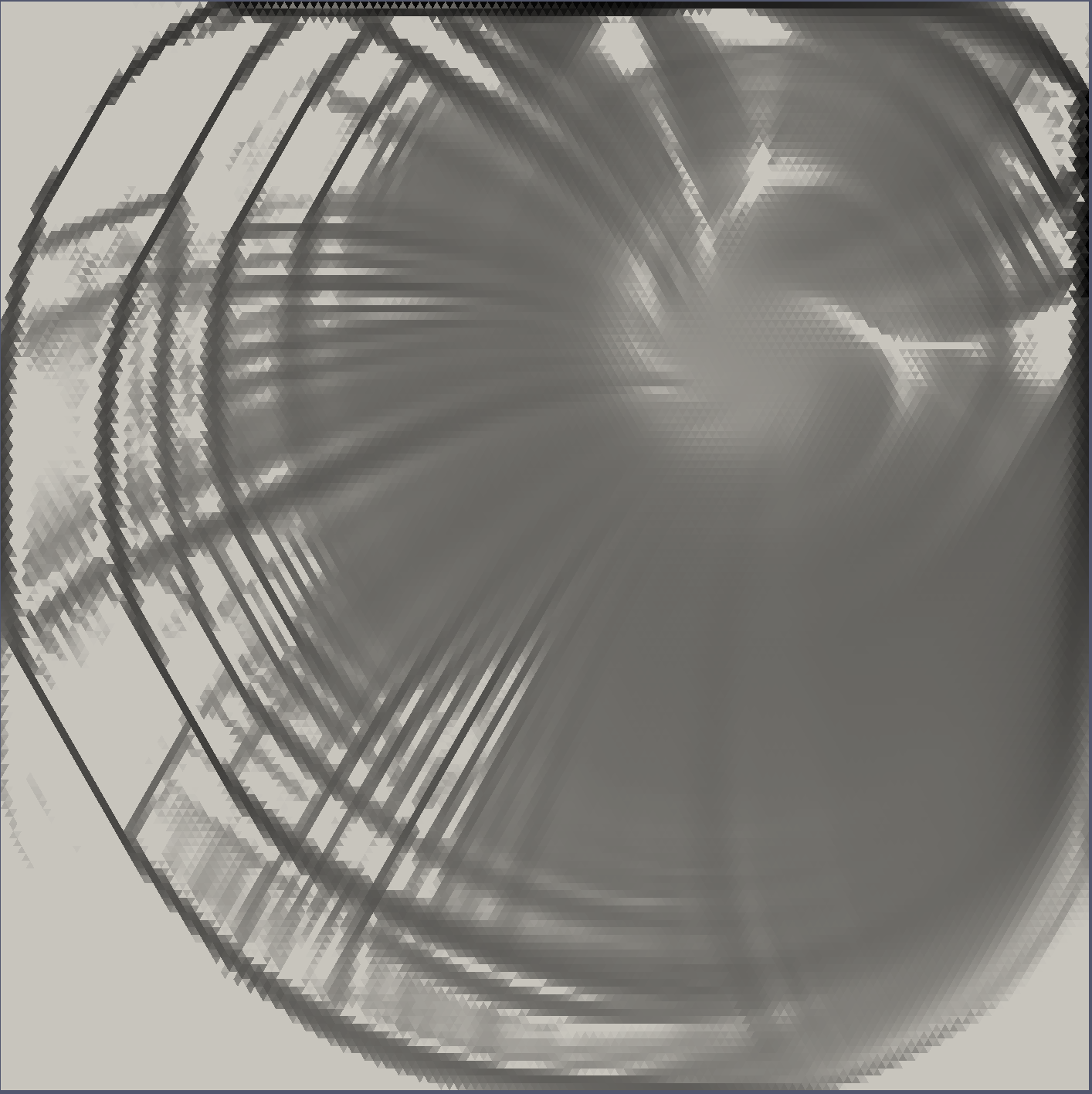}
 &\includegraphics[scale=0.05]{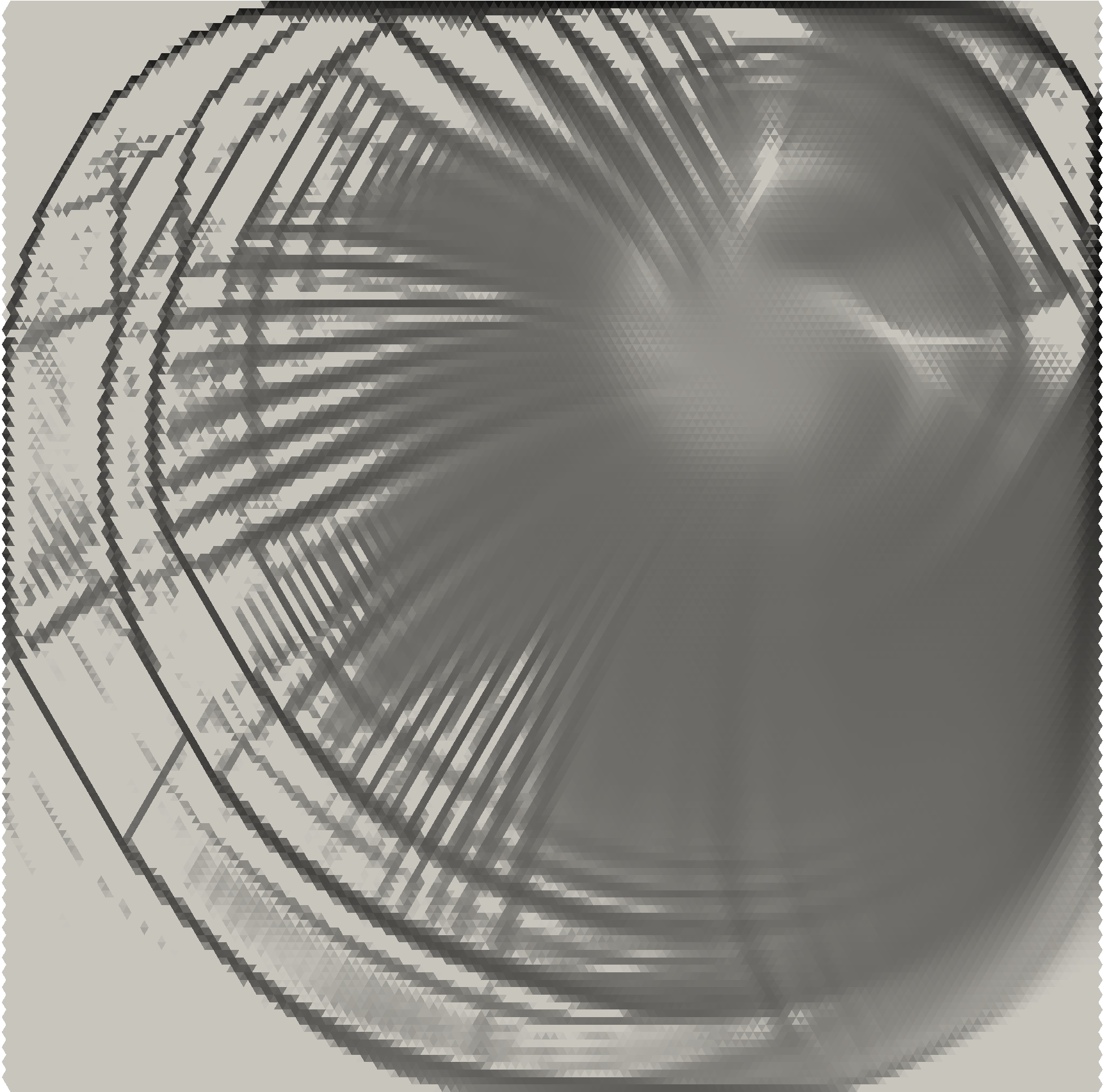}
 \\
 \rotatebox{90}{\phantom{abcdef}2 km}&\includegraphics[scale=0.05]{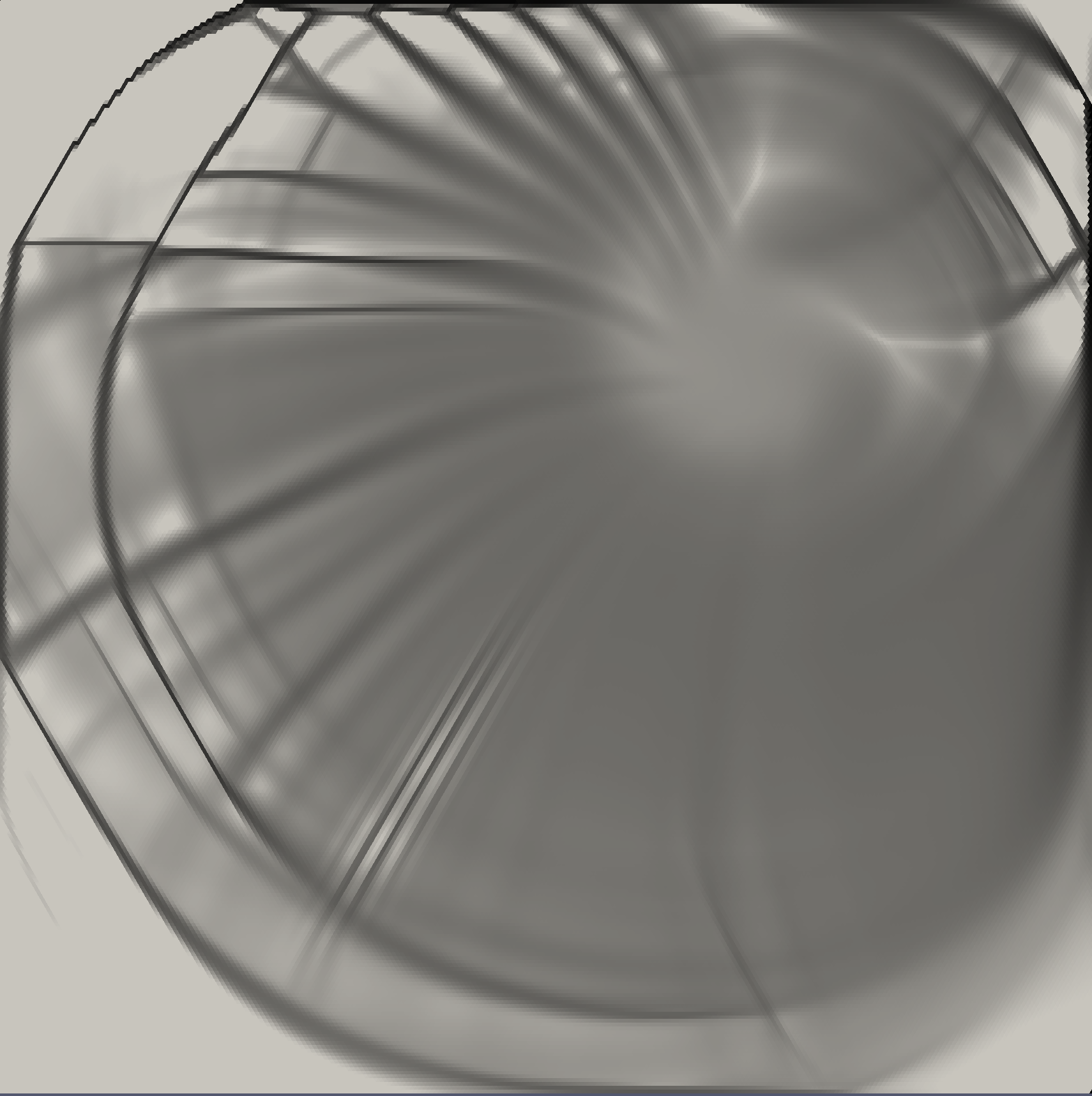}
 &\includegraphics[scale=0.05]{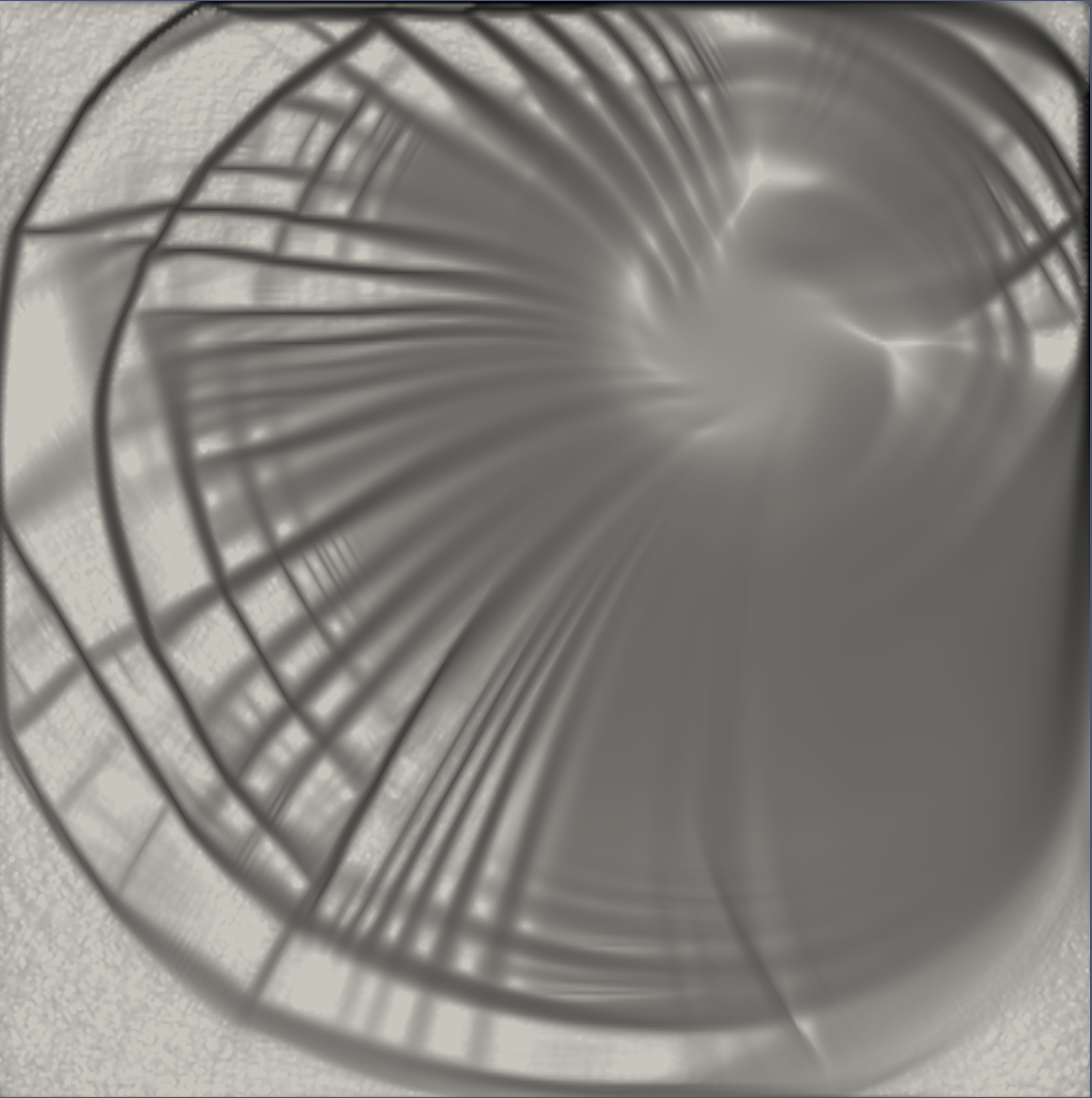}
 &\includegraphics[scale=0.05]{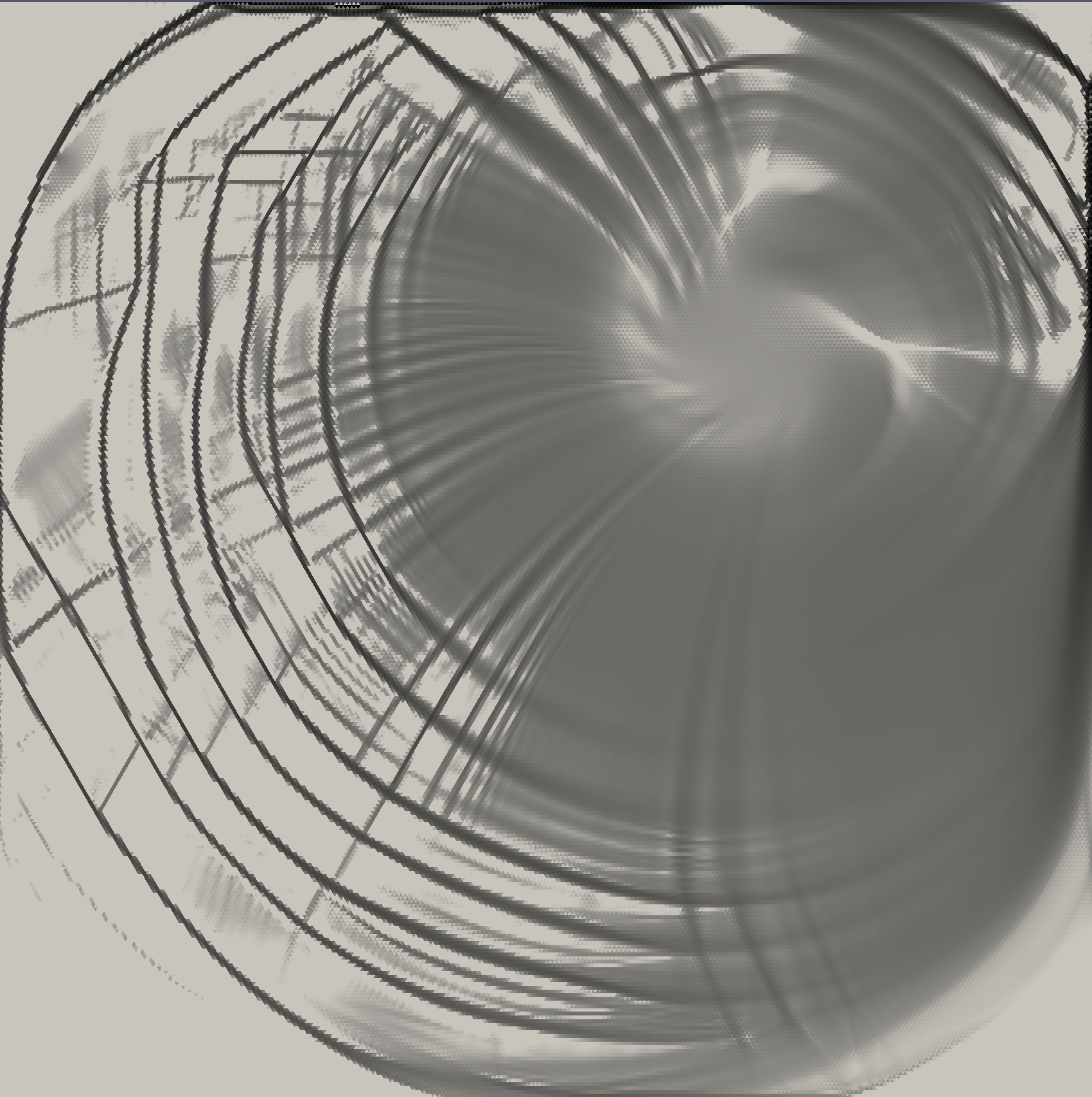}
 &\includegraphics[scale=0.05]{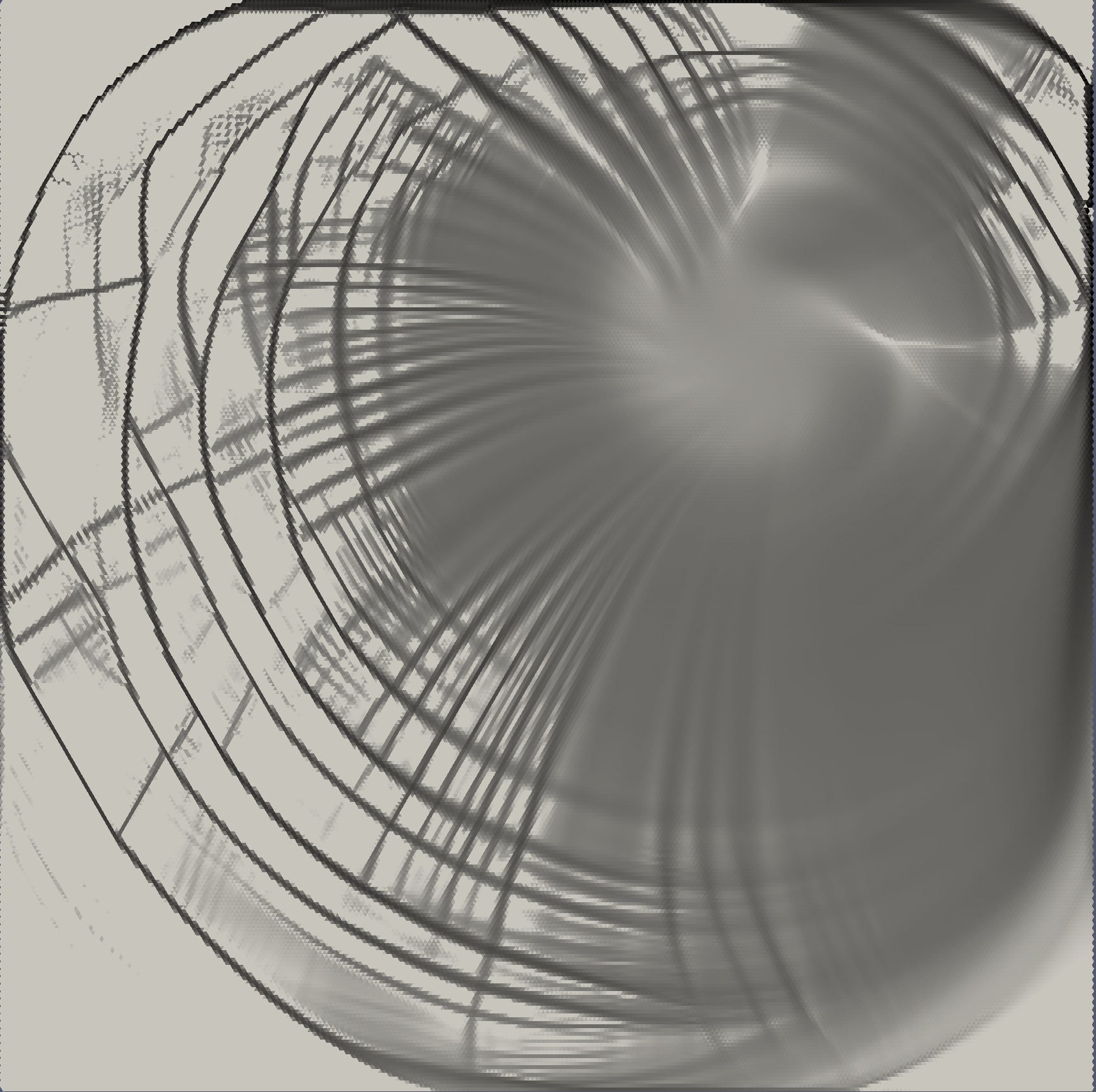}
 \\
 &\includegraphics[scale=0.16]{pics_CR/ICON_sc}&  \includegraphics[scale=0.16]{pics_CR/ICON_sc}&
 \includegraphics[scale=0.16]{pics_CR/ICON_sc}&
 \includegraphics[scale=0.16]{pics_CR/ICON_sc}
 \\
 \end{tabular}
 \caption{Sea ice deformation on a triangular grid. The P1-P1 (A-grid), P1-P0 (B-grid), and the CD-grid analog in form of the CR-P0 pair is computed in FESOM. The CD-grid discretization in ICON refers to the CR-P0 finite element.
 \label{fig:trishear}}
 \end{center}
\end{figure}

\begin{figure}[t]
  \begin{center}
         \includegraphics[scale=0.6]{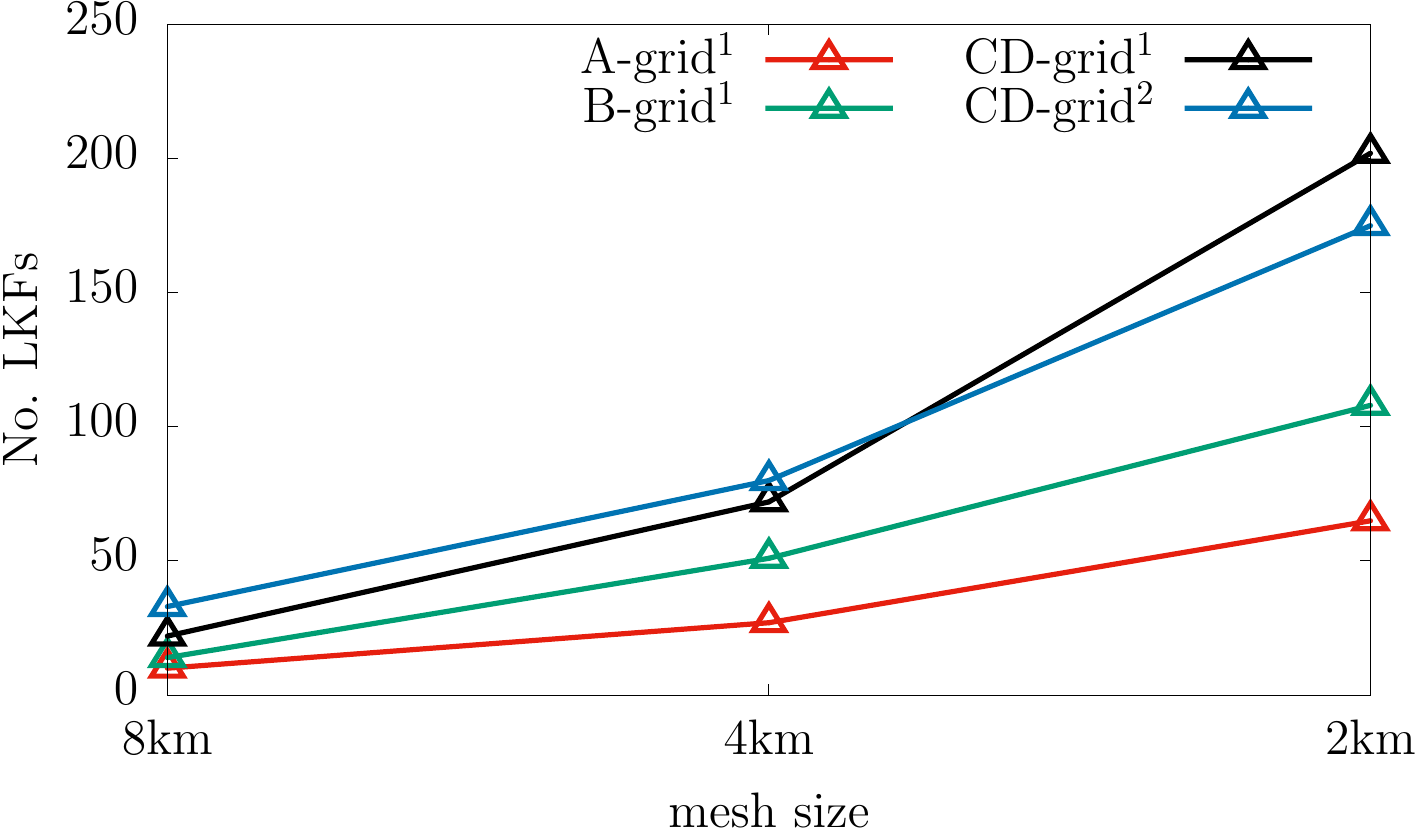}
  \end{center}
  \begin{center}
   \caption{The number of the detected LKFs on triangular meshes from the shear deformation presented in Figure~\ref{fig:trishear}. The subscript 1 and 2 refers to the simulation carried out in the sea ice module of FESOM and ICON respectively.
   \label{fig:vgl_tri}}
  \end{center}
\end{figure}

We start by considering the A-grid like P1-P1 and B-grid type P0-P1 discretization in FESOM. The P0-P1 pair has similarities to the quasi-B grid discretization in FVCOM \citep{Gao2011} and MPAS \citep{Petersen2019}. 
Figure~\ref{fig:triA} and Figure~\ref{fig:trishear} respectively show in the first and second column the sea ice concentration and shear deformation approximated with an A-grid and B-grid type staggering in FESOM. As the only difference between the two approximations is the placement and discretization of the velocity field, we attribute the difference shown in sea ice concentration (Figure~\ref{fig:triA}) and shear deformation (Figure~\ref{fig:trishear}) to the discretization of the sea ice velocity. As the P0-P1 pair has twice as many velocity degrees of freedom compared to the P1-P1 element, it resolves  more LKFs (see Figure~\ref{fig:vgl_tri}). 


Figure~\ref{fig:triA}  shows that the CD-grid approximation in FESOM resolves more LKFs in the sea ice concentration field than the A-grid and B-grid like discretization. This is expected as the CD-grid like CR-P0 staggering has three times more velocity degrees of freedom in the velocity field than the P1-P1 element. Figure~\ref{fig:vgl_tri} shows that we detect up to three times more LKFs in the CD-grid approximation than in the A-grid discretization. The sea ice concentration discretized with the CD-grid approximation in FESOM is advected with an upwind scheme, whereas the A-grid and B-grid type discretizations use a finite element flux correction scheme  based on Taylor-Galerkin method. Changing the transport scheme to a second order finite volume flux-correction scheme does not affect the number of resolved LKFs. However, it does affect the width, the definition and position of LKFs. This will be analyzed in future work.

We proceed with comparing the CD-grid CR-P0 finite element approximation in the framework of ICON and FESOM. In line with the currently used settings, the approximation in ICON is based on equilateral triangles. Thus straight boundaries cannot exactly be resolved and
 the left and right boundaries of the domain consists of outward pointing triangles, see Figure~\ref{fig:triA}. In contrast, the boundaries are straight in FESOM which implies that boundary triangles are not equilateral. Beside this, the discretization is the same. The sea ice concentration shown in Figure~\ref{fig:triA} and the shear deformation presented in Figure~\ref{fig:trishear} are similar in both cases. We conclude that the effect of the different meshes is small. This finding is supported by the results of the LKFs detection algorithm shown in Figure~\ref{fig:vgl_tri}. Both CD-grid like discretizations resolve a similar number of LKFs on the 8km, 4 km and 2 km meshes.

\section{Conclusion}\label{sec:con}
Comparing different viscous-plastic sea ice models, we analyzed the factors that influence the formation of linear kinematic features (LKFs). Naturally, grid resolution has a large effect on the number of simulated LKFs. We also show that the discretization of the velocity field, in particular the staggering of the velocity components on the grid, determines how many LKFs are simulated for a given grid resolution.  The LKFs are triggered when the simulated stress state, which is a function of the velocity, reaches the yield curve in the space of stresses. Both grid resolution and staggering set the underlying degrees of freedom of the velocity field, which in turn are related to the number of LKFs in a simulation.
Other factors such as the geometry of the cell (triangular or quadrilateral), the staggering of the tracers, the advection scheme for the tracers, and the time step are less important for the number of resolved LKFs in our specific test case. These factors mainly affect width and definition of the LKFs and their position in the domain. We expect that details such as the advection scheme would become more important in longer simulations, when the ice is given more time to move over larger distances.

On quadrilateral meshes, the CD-grid discretization doubles the degrees of freedom per grid cell for the velocity field,  
so that with this discretization 
sixteen times fewer
grid cells are required for a similar number of LKFs compared to A, B, and C-grid discretizations. 
The number of LKFs is similar between simulations with the latter three discretizations on quadrilateral meshes. 
The A and CD-grid approximations on quadrilateral meshes resolve as many LKFs as the A and CD-grid equivalents on triangular meshes with the same number of vertices. 

 The CD-grid like finite element discretization of the sea ice velocity on triangular and quadrilateral meshes is numerically instable. The noise is triggered by a nontrivial null space in the strain rate field of the viscous-plastic rheology. 
 The oscillations are controlled by adding a stabilizing term to the momentum equation \citep{Mehlmann2020}.

To summarize, a CD-grid discretization allows us to resolve the same number of LKFs on grids with up to two times coarser resolution than other discretizations. This is an appealing property because simulating realistic LKFs in the viscous-plastic sea ice model requires high spatial resolution.
Otherwise, the type of grid discretization (A, B, or C) and the geometry of the mesh (triangular or quadrilateral cells) have only small effects on the number of simulated LKFs. 


\begin{appendix}

\section{Model configuration}\label{app:conf}

\subsection{MITgcm}
In the MITgcm, the momentum equations and in particular the divergence of the stress tensor are discretized using a finite-volume method on a quadrilateral, curvilinear Arakawa C-grid. The strain rate tensor components are approximated by central differences. For the computation of the viscosities $\zeta$ and $\eta$, the strain rates are combined on a tracer point after averaging $\dot{\epsilon}_{12}^2$, which is naturally defined on corner points, to tracer points (see equation~\ref{deltalimit}). Where necessary, these viscosities are linearly interpolated to corner points. Details of the discretization can be found in the appendix of \citet{LOSCH2010} or at \url{https://mitgcm.readthedocs.io}. In each time step of the simulations in the paper, the nonlinear momentum equation is solved implicitly using a Jacobian-free Newton Krylov (JFNK) solver \citep{Losch2014}. An inexact Newton methods is used, where the linear system is solved with an accuracy depending on the non-linear convergence rate. The non-linear solution required to have a residual $10^{-4}$ times smaller than the initial residual. For details see \citet{Losch2014} or \url{https://mitgcm.readthedocs.io}. Ice thickness and concentration are advected with a second order scheme with a superbee flux limiter.

\subsection{CICE}
CICE is implemented for a quadrilateral, curvilinear B-grid, with velocities located on the grid nodes (corners) of the tracer mesh (``T-grid'').  These nodes are the centers of the dual ``U-grid’' over which the model equations are discretized using a variational approach in order to preserve the energy dissipation properties of the VP equations.  That is, excepting boundary effects, the total rate of internal work equals the dissipation of mechanical energy over the full domain $\Omega$,
\begin{equation}
\int \vt \cdot ( \nabla \cdot \sigma) d\Omega = -\int ( \sigma_{11}\epsilon_{11} + 2\sigma_{12}\epsilon_{12} + \sigma_{22}\epsilon_{22}) d\Omega,
\label{eq:variationalapproachwork}
\end{equation}
where $\vt = (u,v)$ is velocity, $\sigma$ is internal stress and $\dot{\epsilon}_{ij}$ is the strain rate tensor.
Following \citet{Hunke2002}, the stress force components are $F_u = (\nabla\cdot\sigma)_u$ and $F_v = (\nabla\cdot\sigmat)_v$, and the strain rates are derivatives of $\vt$ with general curvilinear metric terms to account for curvature effects of the domain.
Equation~\ref{eq:variationalapproachwork} is differentiated with respect to $u$ and $v$, yielding relationships between the components of $\mathbf{F}$
and the strain rate tensor, via the dissipation rate.  This variational calculation assumes that the stress components are not functions of velocity.  In practice, the differentiation is performed for each U-grid cell, assuming that the velocity is bilinear (that is, linear in $x$ and in $y$) over each cell, and incorporating contributions from all of the surrounding T-grid cells that touch the U-grid centered velocity.  This means that the stress on a given cell boundary (e.g. a T-grid node or U-grid center) may not be the same value as in an adjoining cell, for the same point along the boundary.  For this reason, four values of strain rate (evaluated from the bilinear approximations of $\vt$) and four values of stress and carried for each T-grid cell, representing the four corners of each cell. Derivation of the discretization is presented in detail in \citet{Hunke2002}.

In every time step we solve the momentum equation implicitly by performing 100 nonlinear iterations using a Picard solver. In the implicit solver, the ice velocity $\vt$ appearing in the coefficient $|\vt_{w}-\vt|_2$ in the forcing term \eqref{eq:forcing} is computed as the mean of the two previous nonlinear iterates, following \citet{Lemieux2009}. This linear system is solved using the Flexible Generalized Minimum RESidual (FGMRES, \citep{Saad1996}) method. The FGMRES solver uses a Krylov subspace of dimension 50 and iterates until the relative criterion $10^{-2}$. At each linear iteration, the system is preconditioned by performing five iterations of a Jacobi-preconditioned, GMRES method.  

The tracers in CICE are transported with the incremental remapping scheme \citep{Liscomb2004}. CICE is a multi-thickness category model. Five ice thickness categories were used for the simulations. At initialization, ice of thickness $H$ (following eq. (\ref{initial:H}) is "placed" in the first thickness category. We think that the fact that CICE is based on a multi-thickness category approach has a small impact on the results of this paper. Using only one thickness category results in very similar simulated LKFs (not shown).

\subsection{Gascoigne} 
The model is implemented
in the software library Gascoigne~\citep{Gascoigne}. In the A-grid like Q1-Q1 discretization the velocity as well as the sea ice thickness and the sea ice concentration are co-located at the vertex of a quadrilateral. The B-grid type Q1-Q0 discretization places the velocities at the vertices and the ice concentration and thickness at the cell centers. In the CD-grid like CR-Q0 discretization the velocity is computed at the edge midpoints, whereas the concentration and thickness are calculated at the cell center. The formulation is stabilized as described by \citep{Mehlmann2020}. In the A-grid, B-grid and CD-grid cases the sea ice momentum equation is solved implicitly by a damped Newton method \citep{MehlmannRichter2016newton} and the resulting linear systems of equations are approximated with the GMRES method preconditioned by a parallel multigrid method \citep{MehlmannRichter2016mg,FailerRichter2020}. The nonlinear system is solved to a tolerance of $10^{-13}$, whereas the linear problems in each Newton step are computed with an accuracy of $10^{-2}$. 

In the case of the Q1-Q1 discretization, the advection of the sea ice thickness and concentration is realized with a second order flux-correcting Taylor Galerkin scheme which is described in \citep{Mehlmann2019}. In contrast to the A-grid like Q1-Q1 discretization, a first order upwind scheme is used for the B-grid like Q1-Q0 and CR-Q0 approximations.

\subsection{ICON}
ICON  uses a triangular mesh that is generated by inscribing an icosahedron into a sphere and then subdividing and projecting the subdivided mesh onto the surface of the sphere. The sea ice module in ICON applies the nonconforming Crouzeix-Raviart element to discretize the velocity and uses an upwind scheme to advect the tracers. The discretization corresponds to CD-grid like staggering of the variables. The discretization of the strain rate tensor with the Crouzeix-Raviart element creates a non-trivial  null space. To overcome this shortcoming, the momentum equation needs to be stabilized. A detailed description can be found in \citet{Mehlmann2020}. The nonlinear momentum equation is solved with the mEVP solver of \citet{Kimmritz2015}. To reduce the numerical cost we use 100 sub-iterations per time step on every mesh level. According to the stability criterion, we adjust the parameters $\alpha$ and $\beta$ of the mEVP solver on each mesh level. We use $\alpha=\beta=800$, $\alpha=\beta=2000$, $\alpha=\beta=1200$ on the 8\,km, 4\,km, and 2\,km mesh, respectively.

\subsection{FESOM}
The P1-P1 finite element discretization corresponds to an A-grid staggering, where the velocities and the tracers are stored at the vertices of a triangle. The details of P1-P1 implementation are given in \citet{Danilov2015}.

In the P0-P1 case, the discrete ice velocities are placed at triangles, and scalars are at vertices as in the P1-P1 case. Despite finite-element notation, the implementation of the momentum part follows the finite-volume method. The strain rates are computed at vertices by using median-dual control volumes and replacing the area integral of strain rates by the integral of velocity multiplied with component of outer normal vector over the boundary of median-dual control volumes. The vertex-based strain rates are averaged to edge midpoints. These estimates are further corrected in a least square procedure that increases the weight of across-edge velocity differences in the strain rates.  The ice strength is averaged to edges, and viscosities and stresses are computed at edges too. The stress divergence is computed using the divergence theorem on triangles. Computations of strain rates on vertices with subsequent averaging to edges and least square correction that increases the weight of across-edge velocity difference serve to remove the kernel in the strain rate operator that will be created otherwise. 

The FESOM CP-P0 implementation follows the implementation of the sea ice module in ICON \citep{Mehlmann2020} and differs only by using longitude-latitude coordinates, whereas coordinate systems local to triangles are used in ICON-O. Since the test case is done in plane geometry, both approaches are expected to lead to the same result.

The momentum equation in FESOM is solved with the mEVP solver \citep{Kimmritz2015} using 100 sub-iterations. To ensure stability the  parameters $\alpha$ and $\beta$ of the solver need to adjusted at each mesh resolution. The following triples of $\alpha=\beta$ were used on 8, 4, and 2\,km meshes: (300, 500, 800) for P1 velocities, (500, 800, 1000) for P0 velocities and (800, 1000 and 1500) for the CR velocities. 

In  case of the P1 and P0 velocity discretization the scalar quantities (concentration and thickness) are advanced in time with the FEM-FCT method \citep{Loehner1987}. 
The scalar quantities in the CR-P0 case are constant on triangles and advected with the first-order upwind advection scheme.

\end{appendix}

\bibliographystyle{plainnat}

\begin{thebibliography}{43}
\providecommand{\natexlab}[1]{#1}
\providecommand{\url}[1]{\texttt{#1}}
\expandafter\ifx\csname urlstyle\endcsname\relax
  \providecommand{\doi}[1]{doi: #1}\else
  \providecommand{\doi}{doi: \begingroup \urlstyle{rm}\Url}\fi

\bibitem[Acosta et~al.(2011)Acosta, Apel, Duran, and Lombardi]{Acosta2011}
G.~Acosta, T.~Apel, R.~Duran, and A.~Lombardi.
\newblock Error estimates for {R}aviart-{T}homas interpolation of any order on
  anisotropic tetrahedra.
\newblock \emph{Math. Comput.}, 80:\penalty0 141--163, 2011.

\bibitem[Becker et~al.(2021)Becker, Braack, Meidner, Richter, and
  Vexler]{Gascoigne}
R.~Becker, M.~Braack, D.~Meidner, T.~Richter, and B.~Vexler.
\newblock \emph{The finite element toolkit \textsc{Gascoigne 3d}}, 2021.
\newblock \url{https://www.gascoigne.de}.

\bibitem[Blockley et~al.(2020)Blockley, Vancoppenolle, Hunke, Bitz, Feltham,
  Lemieux, Losch, Maisonnave, Notz, Rampal, Tietsche, Tremblay, Turner,
  Massonnet, Ólason, Roberts, Aksenov, Fichefet, Garric, Iovino, Madec,
  Rousset, Salas~y Melia, and Schroeder]{Blockley2020}
E.~Blockley, M.~Vancoppenolle, E.~C. Hunke, C.~Bitz, D.~L. Feltham, J.-F.
  Lemieux, M.~Losch, E.~Maisonnave, D.~Notz, P.~Rampal, S.~Tietsche,
  B.~Tremblay, A.~Turner, F.~Massonnet, E.~Ólason, A.~Roberts, Y.~Aksenov,
  T.~Fichefet, G.~Garric, D.~Iovino, G.~Madec, C.~Rousset, D.~Salas~y Melia,
  and D.~Schroeder.
\newblock {The Future of Sea Ice Modeling: Where Do We Go from Here?}
\newblock \emph{Bulletin of the American Meteorological Society}, 101\penalty0
  (8):\penalty0 1304--1311, 2020.

\bibitem[Bouillon et~al.(2013)Bouillon, Fichefet, Legat, and
  Madec]{Boullion2013}
S.~Bouillon, T.~Fichefet, V.~Legat, and G.~Madec.
\newblock The elastic-viscous-plastic method revisited.
\newblock \emph{Ocean Modelling}, 71:\penalty0 2--12, 2013.

\bibitem[Coon et~al.(2007)Coon, Kwok, Levy, Pruis, Schreyer, and
  Sulsky]{Coon2007}
M.~Coon, R.~Kwok, G.~Levy, M.~Pruis, H.~Schreyer, and D.~Sulsky.
\newblock {A}rctic {I}ce {D}ynamics {J}oint {E}xperiment {(AIDJEX)} assumptions
  revisited and found inadequate.
\newblock \emph{Journal of Geophysical Research: Oceans}, 112\penalty0 (C11),
  2007.

\bibitem[Crouzeix and Raviart(1973)]{CrouzeixRaviart1973}
M.~Crouzeix and P.-A. Raviart.
\newblock Conforming and nonconforming finite element methods for solving the
  stationary stokes equations. i.
\newblock \emph{Rev. Française Automat. Informat. Recherche Op\'erationnelle
  S\'er. Rouge}, 7:\penalty0 33–75, 1973.

\bibitem[Danilov et~al.(2015)Danilov, Wang, Timmermann, Iakovlev, Sidorenko,
  Kimmritz, Jung, and Schr\"oter]{Danilov2015}
S.~Danilov, Q.~Wang, R.~Timmermann, N.~Iakovlev, D.~Sidorenko, M.~Kimmritz,
  T.~Jung, and J.~Schr\"oter.
\newblock Finite-{E}lement {S}ea {I}ce {M}odel ({FESIM}), version 2.
\newblock \emph{Geosci. Model Dev.}, 8:\penalty0 1747--1761, 2015.

\bibitem[Failer and Richter(2020)]{FailerRichter2020}
L.~Failer and T.~Richter.
\newblock A parallel newton multigrid framework for monolithic fluid-structure
  interactions.
\newblock \emph{Journal of Scientific Computing}, 82\penalty0 (2), 2020.

\bibitem[Feltham(2008)]{Feltham2008}
D.~L. Feltham.
\newblock Sea ice rheology.
\newblock \emph{Annual Review of Fluid Mechanics}, 40\penalty0 (1):\penalty0
  91--112, 2008.

\bibitem[Gao et~al.(2011)Gao, Chen, Qi, and Beardsley]{Gao2011}
G.~Gao, C.~Chen, J.~Qi, and R.~C. Beardsley.
\newblock An unstructured-grid, finite-volume sea ice model: {D}evelopment,
  validation, and application.
\newblock \emph{Journal of Geophysical Research: Oceans}, 116\penalty0 (C8),
  2011.

\bibitem[Girard et~al.(2011)Girard, Bouillon, Weiss, Amitrano, Fichefet,
  Thierry, and Legat]{Girard2011}
L.~Girard, S.~Bouillon, J.~Weiss, D.~Amitrano, T.~Fichefet, L.~Thierry, and
  V.~Legat.
\newblock A new modeling framework for sea-ice mechanics based on
  elasto-brittle rheology.
\newblock \emph{Annals of Glaciology}, 52:\penalty0 123--132, 2011.

\bibitem[Gray and Morland(1994)]{Gray1994}
J.~M. N.~T. Gray and L.~W. Morland.
\newblock A two-dimensional model for the dynamics of sea ice.
\newblock \emph{Philosophical Transactions of the Royal Society of London.
  Series A: Physical and Engineering Sciences}, 347\penalty0 (1682):\penalty0
  219--290, 1994.

\bibitem[Herman(2013)]{Herman2013}
A.~Herman.
\newblock Shear-jamming in two-dimensional granular materials with power-law
  grain-size distribution.
\newblock \emph{Entropy 2013}, 15:\penalty0 4802--4821, 2013.

\bibitem[Hibler(1979)]{Hibler1979}
W.~D. Hibler.
\newblock A dynamic thermodynamic sea ice model.
\newblock \emph{J. Phys. Oceanogr}, 9:\penalty0 815--846, 1979.

\bibitem[Hunke and Dukowicz(1997)]{Hunke1997}
E.~C. Hunke and J.~K. Dukowicz.
\newblock An elastic-viscous-plastic model for sea ice dynamics.
\newblock \emph{J. Phys. Oceanogr.}, 27:\penalty0 1849--1867, 1997.

\bibitem[Hunke and Dukowicz(2002)]{Hunke2002}
E.~C. Hunke and J.~K. Dukowicz.
\newblock {The Elastic–Viscous–Plastic Sea Ice Dynamics Model in General
  Orthogonal Curvilinear Coordinates on a Sphere-Incorporation of Metric
  Terms}.
\newblock \emph{Monthly Weather Review}, 130\penalty0 (7):\penalty0 1848--1865,
  2002.

\bibitem[Hunke et~al.(2015)Hunke, Lipscomb, Turner, Jeffery, and
  Elliott]{Hunke2015}
E.~C. Hunke, W.H. Lipscomb, A.K. Turner, N.~Jeffery, and S.~Elliott.
\newblock \emph{{CICE}: the {L}os {A}lamos {S}ea {I}ce {M}odel {D}ocumentation
  and {S}oftware {U}ser’s {M}anual {V}ersion 5.1 LA-CC-06-012}, 2015.

\bibitem[Hutter and Losch(2020)]{Hutter2020}
N.~Hutter and M.~Losch.
\newblock Feature-based comparison of sea ice deformation in lead-permitting
  sea ice simulations.
\newblock \emph{The Cryosphere}, 14\penalty0 (1):\penalty0 93--113, 2020.

\bibitem[Hutter et~al.(2018)Hutter, Losch, and Menemenlis]{Hutter2018}
N.~Hutter, M.~Losch, and D.~Menemenlis.
\newblock {S}caling properties of {A}rctic sea ice deformation in a
  high-resolution viscous-plastic sea ice model and in satellite observations.
\newblock \emph{J. Geophys. Res.}, 170:\penalty0 18--38, 2018.

\bibitem[Hutter et~al.(2019)Hutter, Zampieri, and Losch]{Hutter2019}
N.~Hutter, L.~Zampieri, and M.~Losch.
\newblock Leads and ridges in {A}rctic sea ice from {RGPS} data and a new
  tracking algorithm.
\newblock \emph{The Cryosphere}, 13\penalty0 (2):\penalty0 627--645, 2019.

\bibitem[Ip et~al.(1991)Ip, Hibler, and Flato]{Hibler1991}
C.~F. Ip, W.~D. Hibler, and G.~M. Flato.
\newblock On the effect of rheology on seasonal sea-ice simulations.
\newblock \emph{Annals of Glaciology}, 15:\penalty0 17--25, 1991.

\bibitem[Kimmritz et~al.(2015)Kimmritz, Danilov, and Losch]{Kimmritz2015}
M.~Kimmritz, S.~Danilov, and M.~Losch.
\newblock On the convergence of the modified elastic-viscous-plastic method for
  solving the sea ice momentum equation.
\newblock \emph{J. Comp. Phys.}, 296:\penalty0 90--100, 2015.

\bibitem[Koldunov et~al.(2019)Koldunov, Danilov, Sidorenko, Hutter, Losch,
  Goessling, Rakowsky, Scholz, Sein, Wang, and Jung]{Koldunov}
N.~Koldunov, S.~Danilov, D.~Sidorenko, N.~Hutter, M.~Losch, H.~Goessling,
  N.~Rakowsky, P.~Scholz, D.~Sein, Q.~Wang, and T.~Jung.
\newblock Fast {EVP} {S}olutions in a {H}igh-{R}esolution {S}ea {I}ce {M}odel.
\newblock \emph{Journal of Advances in Modeling Earth Systems}, 11\penalty0
  (5):\penalty0 1269--1284, 2019.

\bibitem[Kwok et~al.(2008)Kwok, Hunke, Maslowski, and Zhang]{Kwok2008}
R.~Kwok, E.~C. Hunke, D.~Maslowski, and J.~Zhang.
\newblock Variability of sea ice simulations assessed with {RGPS} kinematics.
\newblock \emph{J. Geophys. Res.}, 113\penalty0 (C11), 2008.

\bibitem[Lemieux and Tremblay(2009)]{Lemieux2009}
J.-F. Lemieux and B.~Tremblay.
\newblock Numerical convergence of viscous-plastic sea ice models.
\newblock \emph{J. Geophys. Res.}, 114\penalty0 (C5), 2009.

\bibitem[Lemieux et~al.(2010)Lemieux, Tremblay, Sedl\'{a}\v{c}ek, Tupper,
  Thomas, Huard, and Auclair]{Lemieux2010}
J.-F. Lemieux, B.~Tremblay, J.~Sedl\'{a}\v{c}ek, P.~Tupper, S.~Thomas,
  D.~Huard, and J.P. Auclair.
\newblock Improving the numerical convergence of viscous-plastic sea ice models
  with the {J}acobian-free {N}ewton-{K}rylov method.
\newblock \emph{J. Comp. Phys.}, 229:\penalty0 2840--2852, 2010.

\bibitem[Lemieux et~al.(2012)Lemieux, Knoll, Tremblay, Holland, and
  Losch]{Lemieux2012}
J.-F. Lemieux, D.~Knoll, B.~Tremblay, D.~Holland, and M.~Losch.
\newblock A comparison of the {J}acobian-free {N}ewton-{K}rylov method and the
  {EVP} model for solving the sea ice momentum equation with a viscous-plastic
  formulation: a serial algorithm study.
\newblock \emph{J. Comp. Phys.}, 231:\penalty0 5926--5944, 2012.

\bibitem[Lemieux et~al.(2014)Lemieux, Knoll, Losch, and Girard]{Lemieux2014}
J.-F. Lemieux, D.~Knoll, M.~Losch, and C.~Girard.
\newblock A second-order accurate in time {IM}plicit--{EX}plicit ({IMEX})
  integration scheme for sea ice dynamics.
\newblock \emph{J. Comp. Phys.}, 263:\penalty0 375--392, 2014.

\bibitem[Lipscomb and Hunke(2004)]{Liscomb2004}
W.~H. Lipscomb and E.~C. Hunke.
\newblock {M}odeling {S}ea {I}ce {T}ransport {U}sing {I}ncremental {R}emapping.
\newblock \emph{Monthly Weather Review}, 132\penalty0 (6):\penalty0 1341 --
  1354, 2004.

\bibitem[L\"ohner et~al.(1987)L\"ohner, Morgan, Peraire, and
  Vahdati]{Loehner1987}
R.~L\"ohner, K.~Morgan, J.~Peraire, and M.~Vahdati.
\newblock Finite-element flux-corrected transport ({FEM-FCT}) for the {E}uler
  and {N}avier--{S}tokes equations.
\newblock \emph{Internat. J. Numer. Methods Fluids}, 7:\penalty0 1093--1109,
  1987.

\bibitem[Losch et~al.(2010)Losch, Menemenlis, Campin, Heimbach, and
  Hill]{LOSCH2010}
M.~Losch, D.~Menemenlis, J.-M. Campin, P.~Heimbach, and C.~Hill.
\newblock On the formulation of sea-ice models. {P}art 1: Effects of different
  solver implementations and parameterizations.
\newblock \emph{Ocean Modelling}, 33\penalty0 (1):\penalty0 129 -- 144, 2010.

\bibitem[Losch et~al.(2014)Losch, Fuchs, Lemieux, and Vanselow]{Losch2014}
M.~Losch, A~Fuchs, J.-F. Lemieux, and A.~Vanselow.
\newblock A parallel {J}acobian-free {N}ewton-{K}rylov solver for a coupled sea
  ice-ocean model.
\newblock \emph{J. Comp. Phys.}, 257:\penalty0 901--911, 2014.

\bibitem[Mehlmann and Richter(2017{\natexlab{a}})]{MehlmannRichter2016mg}
C.~Mehlmann and T.~Richter.
\newblock A finite element multigrid-framework to solve the sea ice momentum
  equation.
\newblock \emph{J. Comp. Phys.}, 348:\penalty0 847--861, 2017{\natexlab{a}}.

\bibitem[Mehlmann and Richter(2017{\natexlab{b}})]{MehlmannRichter2016newton}
C.~Mehlmann and T.~Richter.
\newblock A modified global {N}ewton solver for viscous-plastic sea ice models.
\newblock \emph{Ocean Modeling}, 116:\penalty0 96--107, 2017{\natexlab{b}}.

\bibitem[Mehlmann(2019)]{Mehlmann2019}
Carolin Mehlmann.
\newblock \emph{Efficient numerical methods to solve the viscous-plastic sea
  ice model at high spatial resolutions}.
\newblock PhD thesis, Otto-von-Guericke-Universität Magdeburg, Fakultät für
  Mathematik, 2019.
\newblock URL \url{http://dx.doi.org/10.25673/14011}.

\bibitem[Mehlmann and Korn(2021)]{Mehlmann2020}
Carolin Mehlmann and Peter Korn.
\newblock Sea-ice dynamics on triangular grids.
\newblock \emph{J. Comp. Phys.}, 428:\penalty0 110086, 2021.

\bibitem[Petersen et~al.(2019)Petersen, Asay-Davis, Berres, Chen, Feige,
  Hoffman, Jacobsen, Jones, Maltrud, Price, Ringler, Streletz, Turner,
  Van~Roekel, Veneziani, Wolfe, Wolfram, and Woodring]{Petersen2019}
M.~Petersen, X.~Asay-Davis, A.~Berres, Q.~Chen, N.~Feige, M.~Hoffman,
  D.~Jacobsen, P.~Jones, M.~Maltrud, S.~Price, T.~Ringler, G.~Streletz,
  A.~Turner, L.~Van~Roekel, M.~Veneziani, J.~Wolfe, P.~Wolfram, and
  J.~Woodring.
\newblock An {E}valuation of the {O}cean and {S}ea {I}ce {C}limate of {E3SM}
  {U}sing {MPAS} and {I}nterannual {CORE-II} {F}orcing.
\newblock \emph{Journal of Advances in Modeling Earth Systems}, 11\penalty0
  (5):\penalty0 1438--1458, 2019.

\bibitem[Rampal et~al.(2016)Rampal, Bouillon, Olason, and
  Morlighem]{Rampal2016}
P.~Rampal, S.~Bouillon, E.~Olason, and M.~Morlighem.
\newblock ne{XtSIM}: a new lagrangian sea ice model.
\newblock \emph{The Cryosphere}, 10:\penalty0 1055--1073, 2016.

\bibitem[Rannacher and Turek(1992)]{RannacherTurek1992}
R.~Rannacher and T.~Turek.
\newblock Simple nonconforming quadrilateral stokes element.
\newblock \emph{Numer. Meth. Partial Diff. Equ.}, 8:\penalty0 97--111, 1992.

\bibitem[Saad(1996)]{Saad1996}
Y.~Saad.
\newblock \emph{Iterative {M}ethods for {S}parse {L}inear {S}ystems}.
\newblock PWS Publishing Company, 1996.

\bibitem[Tsamados et~al.(2013)Tsamados, Feltham, and Wilchinsky]{Tsamados2013}
M.~L. Tsamados, D.~L. Feltham, and A.~V. Wilchinsky.
\newblock Impact of a new anisotropic rheology on simulations of {A}rctic sea
  ice.
\newblock \emph{Journal of Geophysical Research: Oceans}, 118\penalty0
  (1):\penalty0 91--107, 2013.

\bibitem[Wilchinsky and Feltham(2011)]{Wilchinsky2012}
A.~V. Wilchinsky and D.~L. Feltham.
\newblock Modeling coulombic failure of sea ice with leads.
\newblock \emph{Journal of Geophysical Research: Oceans}, 116\penalty0 (C8),
  2011.

\bibitem[Zhang and Hibler(1991)]{Hibler1997}
J.~L. Zhang and W.~D. Hibler.
\newblock On an efficient numerical method for modeling sea ice dynamics.
\newblock \emph{J. Geophys. Res.}, 102:\penalty0 8691--8702, 1991.

\end{thebibliography}

\end{document}